\documentclass[reqno,12pt,a4paper]{amsart}

\usepackage{mathrsfs}
\usepackage{amssymb}
\usepackage{amsfonts}
\usepackage{latexsym}
\usepackage{amsthm}
\usepackage{graphicx}
\def\lb{\label}

\newcommand{\er}[1]{\textrm{(\ref{#1})}}

\begin{document}

%%%%%%%%%% Some definitions %%%%%%%%%%

%%%%%%%% Equations, theorems %%%%%%%%%
\renewcommand{\theequation}{\arabic{section}.\arabic{equation}}
\theoremstyle{plain}
\newtheorem{theorem}{\bf Theorem}[section]
\newtheorem{lemma}[theorem]{\bf Lemma}
\newtheorem{corollary}[theorem]{\bf Corollary}
\newtheorem{proposition}[theorem]{\bf Proposition}
\newtheorem{definition}[theorem]{\bf Definition}
\theoremstyle{remark}
\newtheorem{remark}[theorem]{\bf Remark}
\newtheorem{example}[theorem]{\bf Example}

%%%%% Alphabet %%%%%
\def\a{\alpha}  \def\cA{{\mathcal A}}     \def\bA{{\bf A}}  \def\mA{{\mathscr A}}
\def\b{\beta}   \def\cB{{\mathcal B}}     \def\bB{{\bf B}}  \def\mB{{\mathscr B}}
\def\g{\gamma}  \def\cC{{\mathcal C}}     \def\bC{{\bf C}}  \def\mC{{\mathscr C}}
\def\G{\Gamma}  \def\cD{{\mathcal D}}     \def\bD{{\bf D}}  \def\mD{{\mathscr D}}
\def\d{\delta}  \def\cE{{\mathcal E}}     \def\bE{{\bf E}}  \def\mE{{\mathscr E}}
\def\D{\Delta}  \def\cF{{\mathcal F}}     \def\bF{{\bf F}}  \def\mF{{\mathscr F}}
\def\c{\chi}    \def\cG{{\mathcal G}}     \def\bG{{\bf G}}  \def\mG{{\mathscr G}}
\def\z{\zeta}   \def\cH{{\mathcal H}}     \def\bH{{\bf H}}  \def\mH{{\mathscr H}}
\def\e{\eta}    \def\cI{{\mathcal I}}     \def\bI{{\bf I}}  \def\mI{{\mathscr I}}
\def\p{\psi}    \def\cJ{{\mathcal J}}     \def\bJ{{\bf J}}  \def\mJ{{\mathscr J}}
\def\vT{\Theta} \def\cK{{\mathcal K}}     \def\bK{{\bf K}}  \def\mK{{\mathscr K}}
\def\k{\kappa}  \def\cL{{\mathcal L}}     \def\bL{{\bf L}}  \def\mL{{\mathscr L}}
\def\l{\lambda} \def\cM{{\mathcal M}}     \def\bM{{\bf M}}  \def\mM{{\mathscr M}}
\def\L{\Lambda} \def\cN{{\mathcal N}}     \def\bN{{\bf N}}  \def\mN{{\mathscr N}}
\def\m{\mu}     \def\cO{{\mathcal O}}     \def\bO{{\bf O}}  \def\mO{{\mathscr O}}
\def\n{\nu}     \def\cP{{\mathcal P}}     \def\bP{{\bf P}}  \def\mP{{\mathscr P}}
\def\r{\rho}    \def\cQ{{\mathcal Q}}     \def\bQ{{\bf Q}}  \def\mQ{{\mathscr Q}}
\def\s{\sigma}  \def\cR{{\mathcal R}}     \def\bR{{\bf R}}  \def\mR{{\mathscr R}}
\def\S{\Sigma}  \def\cS{{\mathcal S}}     \def\bS{{\bf S}}  \def\mS{{\mathscr S}}
\def\t{\tau}    \def\cT{{\mathcal T}}     \def\bT{{\bf T}}  \def\mT{{\mathscr T}}
\def\f{\phi}    \def\cU{{\mathcal U}}     \def\bU{{\bf U}}  \def\mU{{\mathscr U}}
\def\F{\Phi}    \def\cV{{\mathcal V}}     \def\bV{{\bf V}}  \def\mV{{\mathscr V}}
\def\P{\Psi}    \def\cW{{\mathcal W}}     \def\bW{{\bf W}}  \def\mW{{\mathscr W}}
\def\cX{{\mathcal X}}     \def\bX{{\bf X}}  \def\mX{{\mathscr X}}
\def\x{\xi}     \def\cY{{\mathcal Y}}     \def\bY{{\bf Y}}  \def\mY{{\mathscr Y}}
\def\X{\Xi}     \def\cZ{{\mathcal Z}}     \def\bZ{{\bf Z}}  \def\mZ{{\mathscr Z}}
\def\O{\Omega}

%********************* Saburova
\def\be{{\bf e}} \def\bc{{\bf c}}
\def\bv{{\bf v}} \def\bu{{\bf u}}
 \def\mn{\mathrm n}
\def\mm{\mathrm m}
\def\slim{\mathrm s\textrm{-}\lim}
%************************

\newcommand{\mc}{\mathscr {c}}

\newcommand{\gA}{\mathfrak{A}}          \newcommand{\ga}{\mathfrak{a}}
\newcommand{\gB}{\mathfrak{B}}          \newcommand{\gb}{\mathfrak{b}}
\newcommand{\gC}{\mathfrak{C}}          \newcommand{\gc}{\mathfrak{c}}
\newcommand{\gD}{\mathfrak{D}}          \newcommand{\gd}{\mathfrak{d}}
\newcommand{\gE}{\mathfrak{E}}
\newcommand{\gF}{\mathfrak{F}}           \newcommand{\gf}{\mathfrak{f}}
\newcommand{\gG}{\mathfrak{G}}           %\newcommand{\gg}{\mathfrak{g}}
\newcommand{\gH}{\mathfrak{H}}           \newcommand{\gh}{\mathfrak{h}}
\newcommand{\gI}{\mathfrak{I}}           \newcommand{\gi}{\mathfrak{i}}
\newcommand{\gJ}{\mathfrak{J}}           \newcommand{\gj}{\mathfrak{j}}
\newcommand{\gK}{\mathfrak{K}}            \newcommand{\gk}{\mathfrak{k}}
\newcommand{\gL}{\mathfrak{L}}            \newcommand{\gl}{\mathfrak{l}}
\newcommand{\gM}{\mathfrak{M}}            \newcommand{\gm}{\mathfrak{m}}
\newcommand{\gN}{\mathfrak{N}}            \newcommand{\gn}{\mathfrak{n}}
\newcommand{\gO}{\mathfrak{O}}
\newcommand{\gP}{\mathfrak{P}}             \newcommand{\gp}{\mathfrak{p}}
\newcommand{\gQ}{\mathfrak{Q}}             \newcommand{\gq}{\mathfrak{q}}
\newcommand{\gR}{\mathfrak{R}}             \newcommand{\gr}{\mathfrak{r}}
\newcommand{\gS}{\mathfrak{S}}              \newcommand{\gs}{\mathfrak{s}}
\newcommand{\gT}{\mathfrak{T}}             \newcommand{\gt}{\mathfrak{t}}
\newcommand{\gU}{\mathfrak{U}}             \newcommand{\gu}{\mathfrak{u}}
\newcommand{\gV}{\mathfrak{V}}             \newcommand{\gv}{\mathfrak{v}}
\newcommand{\gW}{\mathfrak{W}}             \newcommand{\gw}{\mathfrak{w}}
\newcommand{\gX}{\mathfrak{X}}               \newcommand{\gx}{\mathfrak{x}}
\newcommand{\gY}{\mathfrak{Y}}              \newcommand{\gy}{\mathfrak{y}}
\newcommand{\gZ}{\mathfrak{Z}}             \newcommand{\gz}{\mathfrak{z}}

\def\ve{\varepsilon}   \def\vt{\vartheta}    \def\vp{\varphi}    \def\vk{\varkappa}

\def\A{{\mathbb A}} \def\B{{\mathbb B}} \def\C{{\mathbb C}}
\def\dD{{\mathbb D}} \def\E{{\mathbb E}} \def\dF{{\mathbb F}} \def\dG{{\mathbb G}} \def\H{{\mathbb H}}\def\I{{\mathbb I}} \def\J{{\mathbb J}} \def\K{{\mathbb K}} \def\dL{{\mathbb L}}\def\M{{\mathbb M}} \def\N{{\mathbb N}} \def\dO{{\mathbb O}} \def\dP{{\mathbb P}} \def\R{{\mathbb R}}\def\S{{\mathbb S}} \def\T{{\mathbb T}} \def\U{{\mathbb U}} \def\V{{\mathbb V}}\def\W{{\mathbb W}} \def\X{{\mathbb X}} \def\Y{{\mathbb Y}} \def\Z{{\mathbb Z}}

%%%%% Arrows %%%%%

\def\la{\leftarrow}              \def\ra{\rightarrow}            \def\Ra{\Rightarrow}
\def\ua{\uparrow}                \def\da{\downarrow}
\def\lra{\leftrightarrow}        \def\Lra{\Leftrightarrow}

%%%%% Typography %%%%%

\def\lt{\biggl}                  \def\rt{\biggr}
\def\ol{\overline}               \def\wt{\widetilde}
\def\no{\noindent}

%%%%% Math signs %%%%%

\newcommand{\fr}{\frac}
\newcommand{\tf}{\tfrac}

\def\ul{\underline}

\let\ge\geqslant                 \let\le\leqslant
\def\lan{\langle}                \def\ran{\rangle}
\def\/{\over}                    \def\iy{\infty}
\def\sm{\setminus}               \def\es{\emptyset}
\def\ss{\subset}                 \def\ts{\times}
\def\pa{\partial}                \def\os{\oplus}
\def\om{\ominus}                 \def\ev{\equiv}
\def\iint{\int\!\!\!\int}        \def\iintt{\mathop{\int\!\!\int\!\!\dots\!\!\int}\limits}
\def\el2{\ell^{\,2}}             \def\1{1\!\!1}
\def\sh{\sharp}
\def\wh{\widehat}
\def\bs{\backslash}
\def\intl{\int\limits}
%%%%% Math operations %%%%%

\def\na{\mathop{\mathrm{\nabla}}\nolimits}
\def\sh{\mathop{\mathrm{sh}}\nolimits}
\def\ch{\mathop{\mathrm{ch}}\nolimits}
\def\where{\mathop{\mathrm{where}}\nolimits}
\def\all{\mathop{\mathrm{all}}\nolimits}
\def\as{\mathop{\mathrm{as}}\nolimits}
\def\Area{\mathop{\mathrm{Area}}\nolimits}
\def\arg{\mathop{\mathrm{arg}}\nolimits}
\def\const{\mathop{\mathrm{const}}\nolimits}
\def\det{\mathop{\mathrm{det}}\nolimits}
\def\diag{\mathop{\mathrm{diag}}\nolimits}
\def\diam{\mathop{\mathrm{diam}}\nolimits}
\def\dim{\mathop{\mathrm{dim}}\nolimits}
\def\dist{\mathop{\mathrm{dist}}\nolimits}
\def\Im{\mathop{\mathrm{Im}}\nolimits}
\def\Iso{\mathop{\mathrm{Iso}}\nolimits}
\def\Ker{\mathop{\mathrm{Ker}}\nolimits}
\def\Lip{\mathop{\mathrm{Lip}}\nolimits}
\def\rank{\mathop{\mathrm{rank}}\limits}
\def\Ran{\mathop{\mathrm{Ran}}\nolimits}
\def\Re{\mathop{\mathrm{Re}}\nolimits}
\def\Res{\mathop{\mathrm{Res}}\nolimits}
\def\res{\mathop{\mathrm{res}}\limits}
\def\sign{\mathop{\mathrm{sign}}\nolimits}
\def\span{\mathop{\mathrm{span}}\nolimits}
\def\supp{\mathop{\mathrm{supp}}\nolimits}
\def\Tr{\mathop{\mathrm{Tr}}\nolimits}
\def\BBox{\hspace{1mm}\vrule height6pt width5.5pt depth0pt \hspace{6pt}}

%%%%%%%%%%%%% specialities %%%%%%%%%%%%%%

\newcommand\nh[2]{\widehat{#1}\vphantom{#1}^{(#2)}}
%{{\mathop{#1}\limits^\wedge}\vphantom{#1}^{(#2)}}
\def\dia{\diamond}

\def\Oplus{\bigoplus\nolimits}

%%%%%%%%%%% End of definitions %%%%%%%%%%

%%%%% OLD OLD OLD

\def\qqq{\qquad}
\def\qq{\quad}
\let\ge\geqslant
\let\le\leqslant
\let\geq\geqslant
\let\leq\leqslant

\newcommand{\bea}{\begin{aligned}}
\newcommand{\ena}{\end{aligned}}

\newcommand{\ca}{\begin{cases}}
\newcommand{\ac}{\end{cases}}
\newcommand{\ma}{\begin{pmatrix}}
\newcommand{\am}{\end{pmatrix}}
\renewcommand{\[}{\begin{equation}}
\renewcommand{\]}{\end{equation}}
\def\eq{\begin{equation}}
\def\qe{\end{equation}}
\def\[{\begin{equation}}
\def\bu{\bullet}
\def\bq{\mathbf q}

\title[{Scattering on periodic metric graphs}]
{Scattering on periodic metric graphs}

\date{\today}

\author[Evgeny Korotyaev]{Evgeny Korotyaev}
\address{Dep. of Mathematical Analysis, Saint-Petersburg State University,
 Universitetskaya nab. 7/9, St. Petersburg, 199034, Russia,
\ korotyaev@gmail.com, \
e.korotyaev@spbu.ru,}
\author[Natalia Saburova]{Natalia Saburova}
\address{Dep. of Mathematical Analysis, Algebra and Geometry, Northern (Arctic)
 Federal University, Severnaya Dvina Emb. 17, Arkhangelsk, 163002, Russia,
 \ n.saburova@gmail.com, \ n.saburova@narfu.ru}


\begin{abstract}
We consider the Laplacian on a periodic metric graph and
obtain its decomposition into a direct fiber integral in terms of
the corresponding discrete Laplacian. Eigenfunctions and eigenvalues
of the fiber metric Laplacian are expressed explicitly in terms of
eigenfunctions  and eigenvalues  of the corresponding fiber discrete
Laplacian and eigenfunctions of the Dirichlet problem on the unit
interval. We show that all these eigenfunctions are uniformly
bounded. We apply these results to the periodic metric Laplacian
perturbed by real integrable potentials. We prove the following: a) the wave operators exist and are complete, b) the standard Fredholm determinant is well-defined and is analytic in the
upper half-plane without any modification for any dimension, c) the
determinant  and the corresponding S-matrix satisfy the Birman-Krein
identity.
\end{abstract}

\maketitle

\begin{quotation}

\begin{center}
{\bf Table of Contents}
\end{center}

\vskip 6pt

{\footnotesize

1. Introduction and main results\hfill \pageref{Sec1}\ \ \ \ \

2.  Direct integral for metric Laplacians \hfill
\pageref{Sec2}\ \ \ \ \

3.  Eigenfunctions of metric Laplacians \hfill \pageref{Sec3}\ \ \ \
\

4.  Schr\"odinger operators on periodic metric graphs \hfill
\pageref{Sec4}\ \ \ \ \

5. $d$-Dimensional lattice \hfill \pageref{Sec5}\ \ \ \ \

6. Graphene lattice \hfill \pageref{Sec6}\ \ \ \ \

7. Stanene  lattice \hfill \pageref{Sec7}\ \ \ \ \

Acknowledgments \hfill \pageref{Ack}\ \ \ \ \

References \hfill \pageref{Ref}\ \ \ \ \

\ \ \ \ \ }
\end{quotation}

%%%%%%%%%%%%%%%%%%%%%%%%%%%%%%%%%%%%%%%%%%%%%%%%%%%%%

\vskip 0.25cm

%%%%%%%%%%%%%%%%%%%%%%%%%%%%%%%%%%%%%%%%%%%%%%%%%%%%%

\vskip 0.25cm

\section {\lb{Sec1}Introduction and main results}
\setcounter{equation}{0}

\subsection{Introduction}
We consider a Schr\"odinger operator $H=H_0+Q$ on a periodic
metric graph $\cG$ with the same edge lengths, i.e., on the so-called
periodic equilateral graph. Here $H_0=\D_M$ is a
positive metric Laplace operator with Kirchhoff vertex conditions
and $Q$  is a real integrable  potential. Differential
operators on metric graphs arise naturally as simplified models in
mathematics, physics, chemistry, and engineering.

It is well-known that $\D_M$ is self-adjoint and  its spectrum
consists of an absolutely continuous part and an infinite number of
flat bands (i.e., eigenvalues with infinite multiplicity). The
absolutely continuous spectrum is a union of an infinite number of
spectral bands separated by gaps. Note that the flat bands always exist, mainly
due to the existence of the so-called "Dirichlet-eigenvalues" (the corresponding eigenfunctions vanish on all vertices). The spectrum of metric  Laplacians
has deep relationships with the spectrum of discrete Laplacians $\D$ on
the corresponding discrete graphs, see  \cite{B85},
\cite{BGP08}, \cite{C97} and references therein.

There are a lot of papers, and even books, on the spectrum of
discrete and metric Laplacians  on finite and infinite graphs (see
\cite{BK13}, \cite{Ch97}, \cite{CDS95}, \cite{CDGT88}, \cite{P12}
and references therein).

It is well known that periodic operators can be decomposed into a
direct integral.  The existence of the direct integral for the
discrete Schr\"odinger operator with periodic potentials on periodic graphs was discussed in
many papers (see, e.g., \cite{HS99}, \cite{KSS98}, \cite{SS92}). In the case of a
concrete periodic graph it is not difficult to write down an
explicit expression for the fiber operator. In the case of an
arbitrary periodic graph explicit forms of the fiber operators were
given in \cite{KS14a}, \cite{S13}. In particular, from this direct
integral decomposition it follows that the spectrum of the Schr\"odinger
operator on periodic discrete graphs consists  of an absolutely
continuous part and a finite number of flat bands (i.e., eigenvalues
with infinite multiplicity). The absolutely continuous spectrum
consists of a finite number of bands (intervals) separated by gaps.
In \cite{KS14a}, \cite{KS19} the authors estimated the Lebesgue measure of the
spectrum of the discrete Schr\"odinger operator in terms of
geometric parameters of the graph and the potentials. Furthermore,
they estimated a global variation of the Lebesgue measure of the
spectrum and a global variation of gap-length in terms of the potentials
and geometric parameters of the graph. In \cite{KS16a} the authors considered Laplacians on periodic equilateral metric graphs and estimated the Lebesgue measure of the bands and gaps on a finite interval in terms of geometric parameters
of the graph.

A localization of spectral bands of Laplacians
both on metric and discrete periodic graphs in terms of eigenvalues
of the operator on finite graphs (the so-called eigenvalue
bracketing) was described in \cite{KS14b}, \cite{KS15}, \cite{KS19}, \cite{LP08}.

Effective masses for Laplacians on periodic discrete and metric
equilateral graphs were studied in \cite{KS16b}. The authors
estimated effective masses associated with the ends of each
spectral band in terms of geometric parameters of the graphs. Moreover,
at the beginning of the spectrum they obtained two-sided estimates of
the effective mass in terms of geometric parameters of the graphs.

Results about Laplacians on periodic graphs are used in spectral analysis of the Schr\"odinger operator with a decaying potential and also for the study of the Laplacians on periodic graphs with various defects. We briefly describe these works. The scattering problem for the Schr\"odinger operator with a decaying potential on the lattice $\Z^d$, $d>1$, was considered in the papers \cite{BS99,IK12,IM14,Ko10,KM17,Na14,RoS09,SV01}, see also
references therein. Inverse scattering theory for the discrete Schr\"odinger operators with finitely supported potentials was considered
in \cite{IK12} for the case of the lattice $\Z^d$ and in \cite{A12}
for the case of the hexagonal lattice. The discrete Schr\"odinger operator with decaying potentials on arbitrary periodic graphs was studied in \cite{KSl20}, \cite{PR18}. The Schr\"odinger operator with a potential periodic in some directions and finitely supported in others on arbitrary periodic graphs was investigated in the article \cite{KS17}. In papers \cite{AIM16,AIM18,KS17a,Ku14,SS17} the Laplace and Schr\"odinger operators on periodic graphs with different defects were considered.

\smallskip

We describe the main results of this paper:

{\it i) The direct integral for the metric Laplacian $\D_M$ is described in terms
of the discrete Laplacian $\D$.

ii)  Eigenfunctions of fiber operators for the metric Laplacian
with eigenvalues from the set $\s(\D_M)\sm\s_D$, where $\s_D=\{(\pi j)^2 : j\in\N\}$ is the so-called
Dirichlet spectrum, are described in
terms of eigenfunctions of fiber operators for the discrete
Laplacian~$\D$.

iii) We compute eigenfunctions of fiber operators for the metric Laplacian
corresponding to the Dirichlet spectrum $\s_D$ and determine the multiplicity of this spectrum.

iv) We describe eigenfunctions of the absolutely continuous spectrum of the metric Laplacian $\D_M$ and show that all eigenfunctions of $\D_M$ are uniformly bounded.

v) We consider scattering for the Schr\"odinger operator
$H=\D_M+Q$, where the potential $Q\in L^1(\cG)$   is real.
In particular, we obtain

$\bu $ the existence and completeness of the wave operators;

$\bu $ the standard  Fredholm determinant (without any modification
for any dimension) is well defined and is analytic in the upper
half-plane and its main properties are discussed;

$\bu $ the difference of the resolvents for the Schr\"odinger operator
and for the Laplace operator  belongs to the trace class for any
dimension. Note that for the corresponding Schr\"odinger operators on
$\R^d$, $d\geq2$, it does not hold true.}

\medskip

The proof of these results is based on a detailed analysis of
eigenfunctions of fiber operators for the metric Laplacian  obtained
in the present paper. Note that for short range potentials one needs
to develop Mourre's or Enss's  approaches on periodic metric
graphs.

\subsection{Metric Laplacians} Let $\cG=(\cV,\cE)$ be a connected infinite graph, possibly  having loops and multiple edges and embedded into the space $\R^d$. Here $\cV$ is the set of its vertices and $\cE$ is the set of its unoriented edges. We identify each edge $\be\in\cE$ with the segment $[0,1]$. This introduces an orientation on the edge set $\cE$. An edge starting at a vertex $u$ and ending at a vertex $v$ from $\cV$ will be denoted as the ordered pair $(u,v)\in\cE$ and is said to be \emph{incident} to the vertices. Vertices $u,v\in\cV$ will be called \emph{adjacent} and denoted by $u\sim v$, if $(u,v)\in \cE$. Denote by $I(v)$ the set of all edges from $\cE$ incident to the vertex $v\in\cV$. We define the \emph{degree} ${\vk}_v$ of the vertex $v\in\cV$ as the number of all edges from $I(v)$ with loops counted twice.

Let $\G$ be a lattice of rank $d$ in $\R^d$ with a basis $\A=\{a_1,\ldots,a_d\}$, i.e.,
$$
\G=\Big\{a : a=\sum_{s=1}^dn_sa_s, \; n_s\in\Z,\; s\in\N_d\Big\}, \qqq \N_d=\{1,\ldots,d\},
$$
and let
\[\lb{fuce}
\Omega=\Big\{x\in\R^d : x=\sum_{s=1}^dt_sa_s, \; 0\leq t_s<1,\; s\in\N_d\Big\}
\]
be the \emph{fundamental cell} of the lattice $\G$. We define the equivalence relation on $\R^d$:
$$
x\equiv y \; (\hspace{-4mm}\mod \G) \qq\Leftrightarrow\qq x-y\in\G \qqq
\forall\, x,y\in\R^d.
$$

Below we consider locally finite $\G$-periodic metric equilateral graphs $\cG$, i.e., graphs satisfying the following conditions:
\begin{itemize}
  \item[1)] $\cG=\cG+a$ for any $a\in\G$;
  \item[2)] the quotient graph  $\cG_*=\cG/\G$ is compact;
  \item[3)] all edges of the graph $\cG$ have the same length.
\end{itemize}
The basis $a_1,\ldots,a_d$ of the lattice $\G$ is called the {\it periods} of $\cG$. We also call the quotient graph $\cG_*=\cG/\G$ the \emph{fundamental graph}
of the periodic graph $\cG$. The fundamental graph $\cG_*$ is a graph on the $d$-dimensional torus $\R^d/\G$. The graph $\cG_*=(\cV_*,\cE_*)$ has the vertex set $\cV_*=\cV/\G$ and the set $\cE_*=\cE/\G$ of oriented edges which are finite. Denote by $I_*(v)$ the set of all edges from $\cE_*$ incident to the vertex $v\in\cV_*$.

\setlength{\unitlength}{1.0mm}
\begin{figure}[h]
\centering
\unitlength 1.2mm % = 2.845pt
\linethickness{0.4pt}
\ifx\plotpoint\undefined\newsavebox{\plotpoint}\fi % GNUPLOT compatibility

\begin{picture}(90,50)(0,0)
\put(5,10){\emph{a})}
%\put(7,38){$\cG$}
\put(22.5,22.5){\vector(1,0){15.00}}
\put(22.5,22.5){\vector(0,1){15.00}}

\put(34,20.5){$\scriptstyle a_1$}
\put(19.5,35){$\scriptstyle a_2$}

\multiput(22.5,37.5)(4,0){4}{\line(1,0){2}}
\multiput(37.5,22.5)(0,4){4}{\line(0,1){2}}

\bezier{25}(23.0,22.5)(23.0,30.0)(23.0,37.5)
\bezier{25}(23.5,22.5)(23.5,30.0)(23.5,37.5)
\bezier{25}(24.0,22.5)(24.0,30.0)(24.0,37.5)
\bezier{25}(24.5,22.5)(24.5,30.0)(24.5,37.5)
\bezier{25}(25.0,22.5)(25.0,30.0)(25.0,37.5)
\bezier{25}(25.5,22.5)(25.5,30.0)(25.5,37.5)
\bezier{25}(26.0,22.5)(26.0,30.0)(26.0,37.5)
\bezier{25}(26.5,22.5)(26.5,30.0)(26.5,37.5)
\bezier{25}(27.0,22.5)(27.0,30.0)(27.0,37.5)
\bezier{25}(27.5,22.5)(27.5,30.0)(27.5,37.5)
\bezier{25}(28.0,22.5)(28.0,30.0)(28.0,37.5)
\bezier{25}(28.5,22.5)(28.5,30.0)(28.5,37.5)
\bezier{25}(29.0,22.5)(29.0,30.0)(29.0,37.5)
\bezier{25}(29.5,22.5)(29.5,30.0)(29.5,37.5)
\bezier{25}(30.0,22.5)(30.0,30.0)(30.0,37.5)
\bezier{25}(30.5,22.5)(30.5,30.0)(30.5,37.5)
\bezier{25}(31.0,22.5)(31.0,30.0)(31.0,37.5)
\bezier{25}(31.5,22.5)(31.5,30.0)(31.5,37.5)
\bezier{25}(32.0,22.5)(32.0,30.0)(32.0,37.5)
\bezier{25}(32.5,22.5)(32.5,30.0)(32.5,37.5)
\bezier{25}(33.0,22.5)(33.0,30.0)(33.0,37.5)
\bezier{25}(33.5,22.5)(33.5,30.0)(33.5,37.5)
\bezier{25}(34.0,22.5)(34.0,30.0)(34.0,37.5)
\bezier{25}(34.5,22.5)(34.5,30.0)(34.5,37.5)
\bezier{25}(35.0,22.5)(35.0,30.0)(35.0,37.5)
\bezier{25}(35.5,22.5)(35.5,30.0)(35.5,37.5)
\bezier{25}(36.0,22.5)(36.0,30.0)(36.0,37.5)
\bezier{25}(36.5,22.5)(36.5,30.0)(36.5,37.5)
\bezier{25}(37.0,22.5)(37.0,30.0)(37.0,37.5)
\bezier{25}(37.5,22.5)(37.5,30.0)(37.5,37.5)

\put(29.0,29){$\Omega$}

%***********************
\put(10,15){\line(1,-1){5.00}}
\put(10,15){\line(1,1){5.00}}
\put(20,15){\line(-1,-1){5.00}}
\put(20,15){\line(-1,1){5.00}}
\put(20,15){\line(1,0){5.00}}
\put(15,20){\line(0,1){5.00}}
\put(10,15){\circle{1.1}}
\put(15,10){\circle{1.1}}
\put(20,15){\circle{1.1}}
\put(15,20){\circle{1.1}}
%****************************
\put(25,15){\line(1,-1){5.00}}
\put(25,15){\line(1,1){5.00}}
\put(35,15){\line(-1,-1){5.00}}
\put(35,15){\line(-1,1){5.00}}
\put(35,15){\line(1,0){5.00}}
\put(30,20){\line(0,1){5.00}}
\put(25,15){\circle{1.1}}
\put(30,10){\circle{1.1}}
\put(35,15){\circle{1.1}}
\put(30,20){\circle{1.1}}
%****************************
\put(40,15){\line(1,-1){5.00}}
\put(40,15){\line(1,1){5.00}}
\put(50,15){\line(-1,-1){5.00}}
\put(50,15){\line(-1,1){5.00}}
\put(45,20){\line(0,1){5.00}}
\put(40,15){\circle{1.1}}
\put(45,10){\circle{1.1}}
\put(50,15){\circle{1.1}}
\put(45,20){\circle{1.1}}

%***********************
\put(10,30){\line(1,-1){5.00}}
\put(10,30){\line(1,1){5.00}}
\put(20,30){\line(-1,-1){5.00}}
\put(20,30){\line(-1,1){5.00}}
\put(20,30){\line(1,0){5.00}}
\put(15,35){\line(0,1){5.00}}
\put(10,30){\circle{1.1}}
\put(15,25){\circle{1.1}}
\put(20,30){\circle{1.1}}
\put(15,35){\circle{1.1}}
%****************************
\put(25,30){\line(1,-1){5.00}}
\put(25,30){\line(1,1){5.00}}
\put(35,30){\line(-1,-1){5.00}}
\put(35,30){\line(-1,1){5.00}}
\put(35,30){\line(1,0){5.00}}
\put(30,35){\line(0,1){5.00}}
\put(25,30){\circle*{1.1}}
\put(30,25){\circle*{1.1}}
\put(35,30){\circle*{1.1}}
\put(30,35){\circle*{1.1}}
%****************************
\put(40,30){\line(1,-1){5.00}}
\put(40,30){\line(1,1){5.00}}
\put(50,30){\line(-1,-1){5.00}}
\put(50,30){\line(-1,1){5.00}}
\put(45,35){\line(0,1){5.00}}
\put(40,30){\circle{1.1}}
\put(45,25){\circle{1.1}}
\put(50,30){\circle{1.1}}
\put(45,35){\circle{1.1}}

%***********************
\put(10,45){\line(1,-1){5.00}}
\put(10,45){\line(1,1){5.00}}
\put(20,45){\line(-1,-1){5.00}}
\put(20,45){\line(-1,1){5.00}}
\put(20,45){\line(1,0){5.00}}
\put(10,45){\circle{1.1}}
\put(15,40){\circle{1.1}}
\put(20,45){\circle{1.1}}
\put(15,50){\circle{1.1}}
%****************************
\put(25,45){\line(1,-1){5.00}}
\put(25,45){\line(1,1){5.00}}
\put(35,45){\line(-1,-1){5.00}}
\put(35,45){\line(-1,1){5.00}}
\put(35,45){\line(1,0){5.00}}
\put(25,45){\circle{1.1}}
\put(30,40){\circle{1.1}}
\put(35,45){\circle{1.1}}
\put(30,50){\circle{1.1}}
%****************************
\put(40,45){\line(1,-1){5.00}}
\put(40,45){\line(1,1){5.00}}
\put(50,45){\line(-1,-1){5.00}}
\put(50,45){\line(-1,1){5.00}}
\put(40,45){\circle{1.1}}
\put(45,40){\circle{1.1}}
\put(50,45){\circle{1.1}}
\put(45,50){\circle{1.1}}
%******************************
% ends
\put(10,15){\line(-1,0){2.00}}
\put(10,30){\line(-1,0){2.00}}
\put(10,45){\line(-1,0){2.00}}
\put(50,15){\line(1,0){2.00}}
\put(50,30){\line(1,0){2.00}}
\put(50,45){\line(1,0){2.00}}
\put(15,10){\line(0,-1){2.00}}
\put(30,10){\line(0,-1){2.00}}
\put(45,10){\line(0,-1){2.00}}
\put(15,50){\line(0,1){2.00}}
\put(30,50){\line(0,1){2.00}}
\put(45,50){\line(0,1){2.00}}
%*******************************
\put(23.0,31.5){$\scriptstyle v_3$}
\put(31,24){$\scriptstyle v_2$}
\put(34.5,31.5){$\scriptstyle v_1$}
\put(31,35){$\scriptstyle v_4$}

\put(22.5,20){$\scriptstyle v_4-a_2$}
\put(41,29.5){$\scriptstyle v_3+a_1$}
\put(12.6,29.5){$\scriptstyle v_1-a_1$}
\put(31.0,39.5){$\scriptstyle v_2+a_2$}

%********************************************
\put(62,10){\emph{b})}
%\put(60,38){$\cG_*$}
\put(70,20){\vector(1,0){20.00}}
\put(70,20){\vector(0,1){20.00}}

\put(86,17.5){$a_1$}
\put(66,38){$a_2$}

\put(72.5,31.5){$v_3$}
\put(81,24){$v_2$}
\put(84.5,31.5){$v_1$}
\put(81,35){$v_4$}

\put(75,30){\line(1,-1){5.00}}
\put(75,30){\line(1,1){5.00}}
\put(85,30){\line(-1,-1){5.00}}
\put(85,30){\line(-1,1){5.00}}
\put(85,30){\line(1,0){5.00}}
\put(80,35){\line(0,1){5.00}}
\put(70,30){\line(1,0){5.00}}
\put(80,20){\line(0,1){5.00}}
\put(75,30){\circle*{1.1}}
\put(80,25){\circle*{1.1}}
\put(85,30){\circle*{1.1}}
\put(80,35){\circle*{1.1}}

\multiput(71,40)(4,0){5}{\line(1,0){2}}
\multiput(90,21)(0,4){5}{\line(0,1){2}}

\bezier{30}(70.5,20.0)(70.5,30.0)(70.5,40.0)
\bezier{30}(71.0,20.0)(71.0,30.0)(71.0,40.0)
\bezier{30}(71.5,20.0)(71.5,30.0)(71.5,40.0)
\bezier{30}(72.0,20.0)(72.0,30.0)(72.0,40.0)
\bezier{30}(72.5,20.0)(72.5,30.0)(72.5,40.0)
\bezier{30}(73.0,20.0)(73.0,30.0)(73.0,40.0)
\bezier{30}(73.5,20.0)(73.5,30.0)(73.5,40.0)
\bezier{30}(74.0,20.0)(74.0,30.0)(74.0,40.0)
\bezier{30}(74.5,20.0)(74.5,30.0)(74.5,40.0)
\bezier{30}(75.0,20.0)(75.0,30.0)(75.0,40.0)
\bezier{30}(75.5,20.0)(75.5,30.0)(75.5,40.0)
\bezier{30}(76.0,20.0)(76.0,30.0)(76.0,40.0)
\bezier{30}(76.5,20.0)(76.5,30.0)(76.5,40.0)
\bezier{30}(77.0,20.0)(77.0,30.0)(77.0,40.0)
\bezier{30}(77.5,20.0)(77.5,30.0)(77.5,40.0)
\bezier{30}(78.0,20.0)(78.0,30.0)(78.0,40.0)
\bezier{30}(78.5,20.0)(78.5,30.0)(78.5,40.0)
\bezier{30}(79.0,20.0)(79.0,30.0)(79.0,40.0)
\bezier{30}(79.5,20.0)(79.5,30.0)(79.5,40.0)
\bezier{30}(80.0,20.0)(80.0,30.0)(80.0,40.0)
\bezier{30}(80.5,20.0)(80.5,30.0)(80.5,40.0)
\bezier{30}(81.0,20.0)(81.0,30.0)(81.0,40.0)
\bezier{30}(81.5,20.0)(81.5,30.0)(81.5,40.0)
\bezier{30}(82.0,20.0)(82.0,30.0)(82.0,40.0)
\bezier{30}(82.5,20.0)(82.5,30.0)(82.5,40.0)
\bezier{30}(83.0,20.0)(83.0,30.0)(83.0,40.0)
\bezier{30}(83.5,20.0)(83.5,30.0)(83.5,40.0)
\bezier{30}(84.0,20.0)(84.0,30.0)(84.0,40.0)
\bezier{30}(84.5,20.0)(84.5,30.0)(84.5,40.0)
\bezier{30}(85.0,20.0)(85.0,30.0)(85.0,40.0)
\bezier{30}(85.5,20.0)(85.5,30.0)(85.5,40.0)
\bezier{30}(86.0,20.0)(86.0,30.0)(86.0,40.0)
\bezier{30}(86.5,20.0)(86.5,30.0)(86.5,40.0)
\bezier{30}(87.0,20.0)(87.0,30.0)(87.0,40.0)
\bezier{30}(87.5,20.0)(87.5,30.0)(87.5,40.0)
\bezier{30}(88.0,20.0)(88.0,30.0)(88.0,40.0)
\bezier{30}(88.5,20.0)(88.5,30.0)(88.5,40.0)
\bezier{30}(89.0,20.0)(89.0,30.0)(89.0,40.0)
\bezier{30}(89.5,20.0)(89.5,30.0)(89.5,40.0)
\bezier{30}(90.0,20.0)(90.0,30.0)(90.0,40.0)

\put(79.0,29){$\Omega$}
\end{picture}
\vspace{-10mm} \caption{ \footnotesize \emph{a}) A periodic graph $\cG$, the vectors $a_1,a_2$ are the periods of the graph, the fundamental cell $\Omega$ is shaded; \; \emph{b}) the fundamental graph $\cG_*$, the opposite edges of $\Omega$ are identified.} \label{ff.0.11}
\end{figure}

\begin{example}
We consider the periodic graph $\cG$ shown in Fig.\ref{ff.0.11}\emph{a}. The periods $a_1,a_2$ of the graph and the fundamental cell $\Omega$ are also shown in the figure. The fundamental graph $\cG_*=\cG/\G$ is a graph on the two-dimensional torus $\R^2/\G$, where $\G$ is the lattice generated by the vectors $a_1,a_2$. The torus is obtained from the fundamental cell $\Omega$ by identification of its opposite edges. The fundamental graph $\cG_*$ (Fig.\ref{ff.0.11}\emph{b}) consists of four vertices $v_1,v_2,v_3,v_4$ and six edges
$$
\be_1=(v_1,v_2),\qq \be_2=(v_2,v_3),\qq \be_3=(v_3,v_4),\qq \be_4=(v_4,v_1),\qq \be_5=(v_1,v_3),\qq \be_6=(v_2,v_4).
$$
\end{example}

Recall that each edge $\be\in\cE$ is identified with the segment $[0,1]$. This identification introduces a local coordinate $t\in[0,1]$ along $\be$. For each function $y$ on $\cG$ we define a function $y_{\be}=y\big|_{\be}$, $\be\in\cE$. We identify each function $y_{\be}$ on $\be$ with a function on $[0,1]$ by using the local coordinate $t\in[0,1]$. Let $L^2(\cG)=\os_{\be\in\cE} L^2(\be)$ be
the Hilbert space of all functions $y=(y_\be)_{\be\in\cE}$, where
each $y_\be\in L^2(\be)=L^2(0,1)$, equipped with the norm
$$
\|y\|^2_{L^2(\cG)}=\sum_{\be\in\cE}\|y_\be\|^2_{L^2(0,1)}<\infty.
$$
We define the positive \emph{metric Laplacian} $\D_M$ on
$y=(y_\be)_{\be\in\cE}\in L^2(\cG)$ by
$$
(\D_My)_\be=-y''_\be,\qqq {\rm where}\qqq (y''_\be)_{\be\in\cE}\in
L^2(\cG),
$$
and $y$ satisfies the \emph{Kirchhoff conditions}:
\[
\lb{Dom1} y \textrm{ is continuous on }\cG,\qqq \sum\limits_{\be\in
I(v)}(-1)^{\d(\be,v)}\,y_\be'\big(\d(\be,v)\big)=0,  \qq \forall
v\in\cV,
\]
\[\lb{tev}
\d(\be,v)=\left\{
\begin{array}{rl}
1, \qq &\textrm{if $v$ is the terminal vertex of the edge $\be$} \\
 0, \qq &\textrm{if $v$ is the initial vertex of the edge $\be$}.
\end{array}\right.
\]

\begin{remark}
The spectrum of the metric Laplacian $\D_M$ on
$\cG=(\cV,\cE)$  does not depend on the orientation of the edges
$\be\in\cE$.
\end{remark}

\subsection{Discrete Laplacians.} We introduce the inverse edge
$\ul\be=(v,u)$ for each oriented edge $(u,v)\in \cE$. Let $\cA$ and
$\cA_*$  denote the sets of all edges from $\cE$ and $\cE_*$,
respectively, and their inverse edges. We introduce  the Hilbert
space
\[
\lb{ellV} \ell^2(\cV)=\Big\{f: \cV\ra\C,\
\textstyle\sum\limits_{v\in\cV}\vk_v|f(v)|^2<\infty\Big\},
\]
equipped with the inner product
\[\lb{ipV}
\lan f,g\ran_\cV=\sum_{v\in \cV}\vk_vf(v)\ol{g(v)}.
\]
We define the discrete \emph{normalized Laplace operator} $\D$ on $\ell^2(\cV)$ by
\[\lb{DOL}
\big(\D f\big)(v)=-\frac1{\vk_v}\sum\limits_{(v,\,u)\in\cA}f(u), \qqq f=(f(v))_{v\in\cV}\in\ell^2(\cV),
\]
where ${\vk}_v$ is the degree of the vertex $v\in\cV$ and the sum is
taken over  all oriented edges from $\cA$ starting at the vertex
$v\in\cV$. It is well known (see, e.g., \cite{MW89}) that $\D$ is
self-adjoint and the point $-1$ belongs to its spectrum $\s(\D)$
containing in $[-1,1]$, i.e.,
\[\lb{mp}
-1\in\s(\D)\subset[-1,1].
\]
This Laplacian $\D$ on $\ell^2(\cV)$ has the standard decomposition
into a constant fiber direct integral  for some unitary operator $U:
\ell^2(\cV)\to\gH$:
\[
\lb{raz}
\gH=\int^\oplus_{\T^d}\ell^2(\cV_*)\,{d\vt\/(2\pi)^d}\,,\qqq
U\D U^{-1}=\int^\oplus_{\T^d}\D(\vt)\,{d\vt\/(2\pi)^d}\,,
\]
where $\T^d=\R^d/(2\pi\Z)^d$, the fiber Laplacian $\D(\vt)$ on the fiber space $\ell^2(\cV_*)=\C^\n$ is given by \er{l2.15}, and $\nu=\#\cV_*$ is the number of the fundamental graph vertices. The parameter $\vt$ is called the \emph{quasimomentum}. For each $\vt\in\T^d$ the fiber operator $\D(\vt)$ has $\n$ eigenvalues. It is known (see, e.g., Proposition 4.2 in \cite{HN09}) that  $\l_*$ is an eigenvalue of $\D$ with infinite multiplicity iff $\l_*$ is an eigenvalue of $\D(\vt)$ for any $\vt\in\T^d$. We call $\{\l_*\}$ a \emph{flat band}. We define the multiplicity of flat bands as follows: a flat band $\{\l_*\}$ of $\D$ has multiplicity $m$ iff $\l_*$ is an eigenvalue of $\D(\vt)$ with multiplicity $m$ for almost all $\vt\in\T^d$. Thus, if the operator $\D$ has $r\ge 0$ flat bands, then we denote the corresponding eigenvalues of $\D(\vt)$ (counting multiplicities) by
\[
\lb{fb1}
\l_n=\l_n(\vt)=\const, \qqq \forall\, \n-r< n\le \n=\#\cV_*,\qqq  \forall\,\vt\in \T^d.
\]
All other eigenvalues $\l_1(\vt),\dots,\l_{\n-r}(\vt)$ are not
constant. They can be enumerated in increasing order (counting
multiplicities) by
\[
\label{eq.3.1} \l_1(\vt )\leq\l_2(\vt )\leq\ldots\leq\l_{\nu-r}(\vt),
\qqq \forall\,\vt\in\T^d.
\]
Since $\D(\vt)$ is self-adjoint and analytic in $\vt\in\T^d$, each
$\l_n(\cdot)$,  $n\in\N_\n=\{1,\ldots,\n\}$, is a real and piecewise
analytic function on the torus $\T^d$. Letting $P_n(\vt)$ be the associated
projection, we have
\[
\D(\vt)=\sum_{n=1}^\n \l_n(\vt)P_n(\vt),\qqq  \vt\in\T^d.
\]
We define the \emph{spectral bands} $\s_n(\D)$ by
\[
\lb{ban.1}
\s_n(\D)=[\l_n^-,\l_n^+]=\l_n(\T^d),\qqq n\in\N_{\n}.
\]
Thus, the spectrum of the Laplacian $\D$ on  $\cG$ has the form
$\s(\D)=\bigcup_{n=1}^{\nu}\limits\s_n(\D)$ and we get
\[ \lb{r0}
\textstyle\s(\D)=\s_{ac}(\D)\cup \s_{fb}(\D),\qqq
\s_{ac}(\D)=\bigcup\limits_{n=1}^{\nu-r}\s_n(\D),\qqq
\s_{fb}(\D)=\{\l_{\n-r+1},\ldots,\l_\n\}.
\]
Here and below $\s_{ac}(\D)$ is the absolutely continuous spectrum of $\D$, which is a union of non-degenerate intervals from \er{ban.1}, and
$\s_{fb}(\D)$ is the set of all flat bands (eigenvalues of infinite
multiplicity). An open interval between two neighboring non-degenerate spectral bands is called a \emph{spectral gap}. If $\l_{\nu-r}^+<1$, then it is convenient for us to also call an open interval $(\l_{\nu-r}^+,1)$ a gap of the operator $\D$.

\begin{remark}
It is known that the first spectral band
$\s_1(\D)=[-1,\l_1^{+}]$ of the Laplacian~$\D$ is non-degenerate. The proof is quite simple. Indeed, assume that the first spectral band
$\s_1(\D)$ is degenerate, i.e., the point $-1$ (which belongs to the spectrum of $\D$) is an eigenvalue of the operator $\D$ with infinite multiplicity. Then (see, e.g., Theorem 4.5.2 in \cite{BK13}) there exists an eigenfunction $0\neq f\in\ell^2(\cV)$, corresponding to this eigenvalue and having finite support. By a discrete analogue of the maximum principle (see, e.g., Theorem 7.7 in \cite{S13}), we conclude that $f=0$. We get a contradiction.
\end{remark}

\subsection{Direct integral of metric Laplacians.}
We introduce the Hilbert space (a constant fiber direct integral)
$$
\mH=L^2\big(\T^d,{d\vt
\/(2\pi)^d}\,,\cH\big)=\int_{\T^d}^{\os}\cH\,{d\vt \/(2\pi)^d}\,,
\qqq \cH=L^2(\cG_*),
$$
equipped with the norm
$$
\|g\|^2_{\mH}=\int_{\T^d}\|g(\vt,\cdot )\|_{\cH}^2\frac{d\vt
}{(2\pi)^d}\,,
$$
where the function $g(\vt,\cdot)\in  L^2(\cG_*)$ for all
$\vt\in\T^d$. We have the preliminary standard result about a direct
integral, see, e.g., \cite{RS78}. Recall that $\{a_1,\ldots,a_d\}$ is the basis of the lattice $\G$.

\begin{theorem}\label{T1}
The metric Laplacian $\D_M$ on $L^2(\cG)$ has the
following decomposition into a constant fiber direct integral
\[
\lb{Mraz}
\mU\D_M\mU^{-1}=\int^\oplus_{\T^d}\D_M(\vt)\,{d\vt\/(2\pi)^d}\,,
\]
for the unitary operator $\mU:L^2(\cG)\to\mH$ defined by
\[
\lb{5001}
(\mU h)(\vt,x)=\sum\limits_{m=(m_1,\ldots,m_d)\in\Z^d}e^{-i\lan m,\vt\ran}
h(x+m_1a_1+\ldots+m_da_d), \qq
(\vt,x)\in \T^d\ts \cG_*.
\]
Here the fiber operator $\D_M(\vt)$ acts on $y=(y_\be)_{\be\in\cE_*}\in L^2(\cG_*)$ by
\[\lb{di1}
(\D_M(\vt)y)_\be=-y''_\be,\qqq (y''_\be)_{\be\in\cE_*}\in
L^2(\cG_*),
\]
and the function $y$  satisfies the so-called quasi-periodic
conditions at any $v\in\cV_*$:
\[
\lb{FBC} e^{-i\d(\be_1,v)\lan\t
(\be_1),\,\vt\ran}y_{\be_1}\big(\d(\be_1,v)\big)
=e^{-i\d(\be_2,v)\lan\t
(\be_2),\,\vt\ran}y_{\be_2}\big(\d(\be_2,v)\big), \qq \forall\,
\be_1,\be_2\in I_\ast(v),
\]
\[
\lb{di2} \sum\limits_{\be\in
I_\ast(v)}(-1)^{\d(\be,v)}e^{-i\d(\be,v)\lan\t (\be),\,\vt\ran}\,y_\be'
\big(\d(\be,v)\big)=0,
\]
where $\d(\be,v)$ is defined by \er{tev}, $\t
(\be)\in\Z^d$ is the index of the edge $\be\in\cE_*$, defined in subsection 2.1, and
$\lan\cdot\,,\cdot\ran$ denotes the standard inner product in
$\R^d$.
\end{theorem}

\begin{remark}
In \er{5001} we identify an edge $\be_*\in\cE_*$ of the fundamental graph $\cG_*$ with an edge $\be\in\cE$ of the periodic graph $\cG$ (for more details see subsection 2.2).
\end{remark}

\subsection{Spectra of the metric Laplacian $\D_M$ and its fiber operators $\D_M(\vt)$.} Consider the eigenvalues problem on the unit interval with Dirichlet boundary conditions
\[
\lb{Dp}
-y''=E y,\qqq y(0)=y(1)=0.
\]
It is known that $(\pi j)^2$, $j\in \N$, are the so-called Dirichlet
eigenvalues of the problem \er{Dp}. Thus the spectrum of this
problem is given by
\[
\lb{Dspec} \s_D=\{(\pi j)^2 : j\in\N\}.
\]

We denote by $\b$ the Betti number of the fundamental graph $\cG_*=(\cV_*,\cE_*)$:
\[\lb{Benu}
\b=\#\cE_*-\#\cV_*+1,
\]
where $\#A$ is the number of elements in a set~$A$. The Betti number $\b$ of a finite connected graph $\cG_*$ can also be defined as the dimension of the cycle space $\cC$ of the graph $\cG_*$, i.e.,
\[\lb{becc}
\b=\dim\cC.
\]
Due to the connectivity of the $\G$-periodic graph $\cG$, the fundamental graph $\cG_*$ contains at least $d$ linearly independent cycles, where $d$ is the rank of the lattice $\G$. Then, by \er{becc}, the Betti number of the fundamental graph $\cG_*$ satisfies the inequality $\b\geq d$.

\medskip

We describe our first result about eigenfunctions of the fiber operators $\D_M(\vt)$, $\vt\neq0$, corresponding to the eigenvalues $(\pi j)^2$, $j\in
\N$. Note that in Theorem \ref{TT1} we do not consider the fiber operator $\D_M(0)$ which is just the metric Laplacian on the fundamental graph $\cG_*$ and is well studied (see Remarks \ref{mlfg}.3 and \ref{De0}).

\begin{theorem}\label{TT1}
Let $\D_M(\vt)$, $\vt\in\T^d\sm\{0\}$, be the fiber operator defined by \er{di1} -- \er{di2} on the fundamental graph $\cG_*=(\cV_*,\cE_*)$. Then for all $j\in\N$ the following
statements hold true:

i) $(\pi j)^2$ is an eigenvalue of the operator $\D_M(\vt)$.

ii) Each eigenvalue $(\pi j)^2$ has multiplicity $\b-1$, where $\b$ is the Betti number of the fundamental graph $\cG_*$ defined by \er{Benu}.

iii) The corresponding normalized eigenfunctions $\P_{j,s}^0=\big(\P_{j,s,\be}^0\big)_{\be\in\cE_*}$, $s\in\N_{\b-1}$, have the form:
\[
\lb{eif}
\begin{aligned}
\P_{j,s,\be}^0(\vt,t)=X_{j,s,\be}(\vt)\sqrt{2}\,\sin(\pi jt),
\qq t\in[0,1],\\
\big\|\P_{j,s}^0(\vt,\cdot)\big\|_{L^2(\cG_*)}^2=\sum_{\be\in
\cE_*}\big|X_{j,s,\be}(\vt)\big|^2=1,
\end{aligned}
\]
where $\big(X_{j,s,\be}(\vt)\big)_{\be\in\cE_*}$ is a normalized
solution of the system of $\n$ equations
\[
\lb{sq5} \sum_{\substack{\be\in\cE_\ast, \\ \textrm{ starting at
$v$}}}x_\be-(-1)^j \sum_{\substack{\be\in\cE_\ast, \\ \textrm{ ending
at $v$}}}e^{-i\lan\t ({\bf e}),\,\vt\ran}\,x_\be=0, \qqq \forall\, v\in
\cV_*,
\]
where $\n=\#\cV_*$, one equation for each $v\in\cV_*$. This system has rank $\n$ and there exist $\b-1$ linearly independent normalized solutions $\big(X_{j,s,\be}(\vt)\big)_{\be\in\cE_*}$,
$s\in\N_{\b-1}$. Moreover, $X_{j,s,\be}(\vt)=X_{j+2,s,\be}(\vt)$
for all $(s,\be)\in\N_{\b-1}\ts\cE_*$ and the eigenfunctions
$\P_{j,s}^0(\vt,t)$ satisfy
\[\lb{supv}
\sup_{(j,s)\in\N\ts\N_{\b-1}}\big\|\P_{j,s}^0(\vt,\cdot)\big\|_{L^\iy(\cG_*)}\le \sqrt{2}\,.
\]
\end{theorem}

\begin{remark}\lb{mlfg}
1) In order to determine $X_{j,s,\be}$ in \er{eif}, for each $j\in\N$ we have the
linear system \er{sq5} of $\#\cV_*$  equations with $\#\cE_*$ variables.

2) The eigenfunctions $\P_{j,s}^0=\big(\P_{j,s,\be}^0\big)_{\be\in\cE_*}$, $(j,s)\in\N\ts\N_{\b-1}$, defined by \er{eif} vanish at each vertex of the fundamental graph $\cG_*$.

3) The fiber operator $\D_M(0)$ is the metric Laplacian on the fundamental graph $\cG_*$. Moreover, it is known (see e.g., \cite{B85}, \cite{C97}, \cite{LP08}) that

\emph{i}) $(\pi j)^2$, $j\in\N$, is an eigenvalue of the operator $\D_M(0)$;

\emph{ii}) each eigenvalue $(2\pi j)^2$ of $\D_M(0)$ has multiplicity $\b+1$;

\emph{iii}) each eigenvalue $(2j+1)^2\pi^2$ of $\D_M(0)$ has multiplicity $\b+1$, if $\cG_*$ is bipartite and multiplicity $\b-1$, if $\cG_*$ is not bipartite;

\emph{iv}) the corresponding eigenfunctions are also constructed in the papers mentioned above (for more details see Remark \ref{De0}).

From item \emph{i}) and Theorem \ref{TT1}.\emph{i}) it follows that $(\pi j)^2$, $j\in\N$, is an eigenvalue of $\D_M(\vt)$ for all $\vt\in\T^d$.
\end{remark}

Theorem \ref{TT1} allows us to decompose the fiber operator
$\D_M(\vt)$ in the following form
\[
\lb{1S} \D_M(\vt)=\D_{MD}(\vt)\os \D_{MV}(\vt),\qqq \textrm{where} \qqq
\D_{MD}(\vt)=\sum_{j\ge 1}(\pi j)^2\cP_j^0(\vt),
\]
$\cP_j^0(\vt)$ is the projection associated with the Dirichlet
eigenvalue $(\pi j)^2$ of $\D_M(\vt)$. The operator $\D_{MD}(\vt)$
corresponds  to the Dirichlet spectrum $\s_D$. The other part
$\D_{MV}(\vt)$ is associated with the vertex set $\cV_*$ and it is
discussed in the next theorem. Without loss of generality we may
assume that for a bipartite periodic graph $\cG$ the fundamental
graph $\cG_*$ is also bipartite. Indeed, if $\cG_*$ is non-bipartite, then we just need to go to the vectors $2a_1,\ldots,2a_d$ as the periods of the graph $\cG$, gluing together several copies of $\cG_*$ to a new bipartite fundamental
graph. For example, if $\cG=\Z$, then the bipartite fundamental graph $\cG_*$ is obtained by identification of every second vertex of $\cG$.

We describe the "vertex fiber Laplacian" $\D_{MV}(\vt)$ and its
eigenfunctions. These  eigenfunctions are expressed in terms of
eigenfunctions of the discrete fiber Laplacians.

\begin{theorem}\label{TT2}
The operator $\D_{MV}(\vt)$ defined by \er{1S} has the following form
\[
\lb{DVS} \D_{MV}(\vt)=\sum_{(j,n)\in \Z_0\ts\N_\n}
z_{j,n}^2(\vt)\cP_{j,n}(\vt),\qqq \textrm{ for all} \qq
\vt\in\T^d\setminus\{0\},
\]
where $\cP_{j,n}(\vt)$ is the projection associated with the
eigenvalue  $z_{j,n}^2(\vt)$ of $\D_{MV}(\vt)$ given by
\[
\lb{egv1}
\begin{aligned}
z_{j,n}(\vt)=\ca
  z_n(\vt)+\pi j, & \textrm{if $j$ is even} \\
  (\pi-z_n(\vt))+\pi j, & \textrm{if $j$ is odd}
\ac\;, \qq z_n(\vt)=\arccos(-\l_n(\vt))\in [0,\pi],
\end{aligned}
\]
$\Z_0=\{0,1,2,\ldots\}$. The corresponding normalized eigenfunctions
$\P_{j,n}(\vt)=\big(\P_{j,n,\be}(\vt,t)\big)_{\be\in\cE_*}$ satisfy:
\[\lb{ms6}
\begin{aligned}
& \P_{j,n,\be}(\vt,t)=\textstyle {\sqrt{2}\/\sin z_n(\vt)}\,
\big(\p_n(\vt,u)\sin(z_{j,n}(\vt)\,(1-t))+\p_n(\vt,v)e^{i\lan\t
(\be),\,\vt\ran}\,\sin(z_{j,n}(\vt)\,t)\big),\\
\end{aligned}
\]
\[
\lb{ms8}
\sup _{(j,n,\vt)\in \Z_0\ts\N_\n\ts
 \T^d}\|\P_{j,n}(\vt)\|_{L^\iy(\cG_*)}=C<\iy,
\]
where $\be=(u,v)$, and $\p_n(\vt,\cdot)\in\ell^2(\cV_*)$ is the
normalized eigenfunction of the discrete fiber Laplacian $\D(\vt)$
given by \er{l2.15}, corresponding to the eigenvalue $\l_n(\vt)$.
\end{theorem}

\begin{remark}
1) The eigenvalue  $z_{j,n}^2(\vt)\neq (\pi
j)^2$ for all $\vt\in\T^d\setminus\{0\}$ and all $(j,n)\in\Z_0\ts\N_\n$.

2) Recently, Pankrashkin (see Proposition 1 in \cite{Pa12},  Theorem
14 in \cite{Pa13}) showed the following result. Let $\e(z)=-\cos
\sqrt z$ and $\c_\omega$ be the characteristic function of the set
$\omega\ss \R$ and let an interval $J\ss\R\sm \s_D$. Then the operator
$A=\D_{MV} \c_J(\D_{MV})$ is unitarily equivalent to the operator
$B=\e^{-1}(\D \c_{\e(J)}(\D))$, i.e., the identity
$A=UBU^{-1}$ holds true for some
unitary operator $U$.

3) We do not know any paper about \er{ms6} -- \er{ms8}. These
results are important to study the Fredholm determinant in Theorem
\ref{TSD}.

4)  The explicit form of the constant $C$ in \er{ms8} is defined  in
Propositions \ref{T.As}, \ref{T.As1}. Results about the effective
masses for Laplacians \cite{KS16b} are essential in the proof of \er{ms8}.
\end{remark}

Theorems \ref{T1}, \ref{TT1} and \ref{TT2} imply the following well-known results, see,
e.g., \cite{BK13}, \cite{BGP08}, \cite{P12}:

\begin{corollary}\label{TsDM}
i) The spectrum of the Laplacian $\D_M$ on $L^2(\cG)$ is given by
\[\lb{sDM1}
\begin{aligned}
&\s(\D_M)=\s_{ac}(\D_M)\cup \s_{fb}(\D_M),\qq
\s_{ac}(\D_M)=\hspace{-3mm}\bigcup_{(j,n)\in\Z_0\ts\N_{\n-r}}\hspace{-3mm}\s_{j,n}(\D_M),
 \qq \s_{j,n}(\D_M)=[E_{j,n}^-,E_{j,n}^+],
 \\
&E_{j,n}^\pm=\ca\big(z_n^\pm+\pi j\big)^2, & \; \textrm{if $j$ is even}
\\
\big(\pi-z_n^\mp+\pi j\big)^2, & \; \textrm{if $j$ is odd} \ac, \qq
z_n^\pm=\arccos(-\l_n^\pm)\in [0,\pi], \qq n\in\N_{\n-r},
\end{aligned}
\]
and
\[\lb{sDM1a}
\begin{aligned}
&\s_{fb}(\D_M)=\s_{fb}(\D_{MD})\cup\s_{fb}(\D_{MV}),\qqq
\s_{fb}(\D_{MD})=\s_D=\{(\pi j)^2 : j\in\N\},\\
&\s_{fb}(\D_{MV})=\bigcup_{j=0}^{\infty}\{E_{j,\n-r+1},\ldots,E_{j,\n}\},\qqq
E_{j,n}=\ca\big(z_n+\pi j\big)^2, & \; \textrm{if $j$ is even} \\
\big((\pi-z_n)+\pi j\big)^2, & \; \textrm{if $j$ is odd}
\ac,\\
& z_n=\arccos(-\l_n)\in (0,\pi),
\qqq n=\n-r+1,\ldots,\n,\\
\end{aligned}
\]
where $\l_n^\pm$, $n\in\N_{\n-r}$, are defined in \er{ban.1}, and $r\ge 0$ is the number of the flat bands $\{\l_{\n-r+1}\},\ldots,\{\l_\n\}$ of the discrete Laplacian $\D$ (counting multiplicities).

ii) $\s(\D_M)=[0,+\iy)$ iff $\s(\D)=[-1,1]$.

iii) Both the sets $\s_{ac}(\sqrt{\D_M}\,)$ and $\s_{fb}(\sqrt{\D_M}\,)$
 are $2\pi$-periodic on the half-line $(0,+\infty)$ and  are symmetric on the interval $(0,2\pi)$ with respect to the point $\pi$.

iv) The spectrum of the operator $\D_M$ has exactly $k$ gaps on the interval $[0,\pi^2]$ iff $\s(\D)$ has exactly $k$ gaps. The number of flat bands of $\D_M$ on $[0,\pi^2]$ is $r+1$.

\end{corollary}

\setlength{\unitlength}{1.0mm}
\begin{figure}[h]
\centering
\unitlength 0.9mm % = 2.845pt
\linethickness{0.4pt}
\ifx\plotpoint\undefined\newsavebox{\plotpoint}\fi % GNUPLOT compatibility
\begin{picture}(170,80)(0,0)

\put(-5,40){\vector(1,0){170.00}}

\put(0,5){\vector(0,1){74.00}}

\put(0.2,10){\line(0,1){22.20}}
\put(-0.2,10){\line(0,1){22.20}}
\put(0.4,10){\line(0,1){22.2}}
\put(-0.4,10){\line(0,1){22.20}}

\multiput(0,32)(4,0){41}{\line(1,0){2}}
\multiput(0,45)(4,0){41}{\line(1,0){2}}
\multiput(0,53)(4,0){41}{\line(1,0){2}}
\multiput(0,66)(4,0){41}{\line(1,0){2}}

\put(0,45){\circle*{1.0}}

\put(0.4,53){\line(0,1){13.20}}
\put(0.2,53){\line(0,1){13.20}}
\put(-0.2,53){\line(0,1){13.20}}
\put(-0.4,53){\line(0,1){13.20}}

\put(-13.0,9){$\scriptstyle\l_1^-=-1$}
\put(-5.0,31){$\scriptstyle\l_1^+$}
\put(-5.0,65){$\scriptstyle\l_2^+$}
\put(-5.0,53){$\scriptstyle\l_2^-$}
\put(-4.0,37.0){$\scriptstyle 0=z_1^-$}
\put(23.0,36.8){$\scriptstyle z_1^+$}
\put(95.0,41.5){$\scriptstyle 2\pi-z_1^+$}
\put(147.0,41.5){$\scriptstyle 2\pi+z_1^+$}
\put(8.0,41.5){$\scriptstyle \s_1(\D_M^{1/2})$}
\put(39.5,41.5){$\scriptstyle \s_2(\D_M^{1/2})$}
\put(22.0,20.0){$\l=-\cos z$}
\put(-4.5,45){$\scriptstyle\l_3$}
\put(-9.5,59){$\scriptstyle\s_2(\D)$}
\put(-3.5,70){$\scriptstyle 1$}
\put(-9.5,20){$\scriptstyle\s_1(\D)$}
\put(-1,10){\line(1,0){2.00}}
\put(-1,70){\line(1,0){2.00}}

\multiput(-1,10)(4,0){41}{\line(1,0){2}}
\multiput(-1,70)(4,0){41}{\line(1,0){2}}
\bezier{600}(0,10)(15,9)(31.4,40)
\bezier{600}(94.2,40)(62.8,100)(31.4,40)
\bezier{600}(94.2,40)(125.6,-20)(157,40)

\put(34.2,40.0){\circle*{1.0}}
\put(32.5,37.0){$\scriptstyle z_3$}
\put(86.0,37.){$\scriptstyle 2\pi-z_3$}
\multiput(26.9,32.5)(0,2){4}{\line(0,1){1}}
\multiput(34.2,40.0)(0,2){3}{\line(0,1){1}}
\multiput(51.2,40.7)(0,2){13}{\line(0,1){1}}
\multiput(39.2,40.1)(0,2){7}{\line(0,1){1}}

\put(62.8,40.0){\circle*{1.0}}
\put(125.6,40.0){\circle*{1.5}}
\put(62.0,37.5){$\scriptstyle \pi$}
\put(124.0,36.5){$\scriptstyle 2\pi$}

\multiput(74.4,40.7)(0,2){13}{\line(0,1){1}}
\multiput(86.4,40.1)(0,2){7}{\line(0,1){1}}
\multiput(98.7,32.5)(0,2){4}{\line(0,1){1}}
\multiput(91.4,40.0)(0,2){3}{\line(0,1){1}}
\multiput(152.4,32.5)(0,2){4}{\line(0,1){1}}
\put(74.4,40.3){\line(1,0){12}}
\put(74.4,40.2){\line(1,0){12}}
\put(74.4,40.1){\line(1,0){12}}
\put(74.4,39.9){\line(1,0){12}}
\put(74.4,39.8){\line(1,0){12}}
\put(74.4,39.7){\line(1,0){12}}
\put(91.4,40.0){\circle*{1.0}}

\put(98.7,40.3){\line(1,0){53.8}}
\put(98.7,40.2){\line(1,0){53.8}}
\put(98.7,40.1){\line(1,0){53.8}}
\put(98.7,39.9){\line(1,0){53.8}}
\put(98.7,39.8){\line(1,0){53.8}}
\put(98.7,39.7){\line(1,0){53.8}}
\put(39.2,40.3){\line(1,0){12}}
\put(39.2,40.2){\line(1,0){12}}
\put(39.2,40.1){\line(1,0){12}}
\put(39.2,39.9){\line(1,0){12}}
\put(39.2,39.8){\line(1,0){12}}
\put(39.2,39.7){\line(1,0){12}}
\put(39.0,37.0){$\scriptstyle z_2^-$}
\put(50.0,37.0){$\scriptstyle z_2^+$}

\put(0,40.3){\line(1,0){26.9}}
\put(0,40.2){\line(1,0){26.9}}
\put(0,40.1){\line(1,0){26.9}}
\put(0,39.9){\line(1,0){26.9}}
\put(0,39.8){\line(1,0){26.9}}
\put(0,39.7){\line(1,0){26.9}}

\put(162.0,36.0){$z$}

\put(-4,76.0){$\l$}

\end{picture}
%\vspace{-0.5cm}
\caption{\footnotesize Relation between the spectra of $\D$ and $\sqrt{\D_M}$.}
\label{fRel}
\end{figure}

\begin{remark}
1) Recall that Br\"uning-Geyler-Pankrashkin
\cite{BGP08} proved the following identities
\[
\lb{DMDx} \s_\a(\D_M)\sm\s_D=\{E\in\R\sm\s_D: -\cos\sqrt E\in
\s_\a(\D)\},\qqq \a\in\{\rm fb,\ ac\},
\]
where $\s_D=\{(\pi j)^2 : j\in\N\}$. The relation between the spectra of the operators $\D$ and $\sqrt{\D_M}$ is shown in Fig.\ref{fRel}. The spectrum of the discrete Laplacian $\D$ is along the vertical axis $\l$. It consists of an absolutely continuous part which is a union of a finite number of the non-degenerate spectral bands ($\s_1(\D)$ and $\s_2(\D)$ in the figure) and a finite number of the flat bands ($\{\l_3\}$ in the figure). The spectrum of $\sqrt{\D_M}$ is along the horizontal axis $z$.
In order to obtain $\s(\sqrt{\D_M}\,)$ one needs to find the preimage of $\s(\D)$ under the function $\l=-\cos z$. Each spectral band and each flat band of the discrete Laplacian $\D$ generate infinitely many spectral bands and flat bands, respectively, for the operator $\sqrt{\D_M}$. For example, in the figure the flat band $\{\l_3\}$ of the discrete Laplacian $\D$ generates the eigenvalues $z_3,2\pi-z_3,2\pi+z_3,4\pi-z_3,\ldots$ with infinite multiplicity.  Moreover, $\sqrt{\D_M}$ has the additional eigenvalues $\pi,2\pi,3\pi,\ldots$ with infinite multiplicity. We note that the spectrum of $\sqrt{\D_M}$ is $2\pi$-periodic on $(0,+\iy)$ and is symmetric with respect to the point $\pi$ on the interval $(0,2\pi)$.

2) If $\s(\D)\neq[-1,1]$, then in the spectrum of $\D_M$ there exist \textbf{infinitely many} gaps $\g_1,\g_2,\ldots$ and $|\g_n|\rightarrow\infty$ as $n\rightarrow\infty$.

3) The identities \er{sDM1} show that spectral
properties of the discrete Laplacian $\D$ are important to study
spectral properties of the metric  Laplacian $\D_M$.
\end{remark}

Now we discuss the connection between generalized eigenfunctions of the absolutely continuous spectrum of the discrete and metric Laplacians and eigenfunctions of their fiber operators. Let $\ell^2_{loc}(\cV)$ be the set of all sequences on $\cV$, and let $L^2_{loc}(\cG)$ be the set of all functions which are square integrable on any finite subset in $\cE$. Recall that $\{a_1,\ldots,a_d\}$ is the basis of the lattice $\G$.

\begin{proposition}
\label{TLEF}
Let $\vt\in\T^d\sm\{0\}$, $j\in\{0,1,2,\ldots\}$, $n\in\{1,2,\ldots,\n-r\}$. Let $\p_n(\vt,\cdot)\in\ell^2(\cV_*)$ be the normalized eigenfunction of the fiber discrete Laplacian $\D(\vt)$ corresponding to the eigenvalue $\l_n(\vt)\in\s_{ac}(\D)$ and let $\P_{j,n}(\vt,\cdot)\in L^2(\cG_*)$ be the normalized eigenfunction of the fiber metric Laplacian $\D_M(\vt)$ corresponding to the eigenvalue $z_{j,n}^2(\vt)\in\s_{ac}(\D_M)$, where $z_{j,n}^2(\vt)$ is defined by \er{egv1}. Then the following statements hold true.

i) The generalized eigenfunction $\wt\p_n(\vt,\cdot)\in\ell^2_{loc}(\cV)$ of the discrete Laplacian $\D$ corresponding to $\l_n(\vt)\in\s_{ac}(\D)$ has the form
\[\lb{gde}
\begin{aligned}
\wt\p_n(\vt,v+m_1a_1+\ldots+m_da_d)=e^{i\lan m,\,\vt\ran}\p_n(\vt,v), \\ \forall\,v\in \cV_*, \qqq \forall\, m=(m_1,\ldots,m_d)\in\Z^d.
\end{aligned}
\]

ii) The generalized eigenfunction $\wt\P_{j,n}(\vt,\cdot)\in L^2_{loc}(\cG)$ of the metric Laplacian $\D_M$ corresponding to $z_{j,n}^2(\vt)\in\s_{ac}(\D_M)$ has the form
\[\lb{gme}
\begin{aligned}
\wt\P_{j,n}(\vt,x+m_1a_1+\ldots+m_da_d)=e^{i\lan m,\,\vt\ran}\P_{j,n}(\vt,x),\\ \qqq \forall\, x\in\cG_*, \qqq \forall\, m=(m_1,\ldots,m_d)\in\Z^d.
\end{aligned}
\]

iii) The generalized eigenfunctions from \er{gde}, \er{gme} satisfy
\[\lb{gdep}
\wt\p_n(\vt,v)=e^{i\lan v_\A,\,\vt\ran}\vp_n(\vt,v), \qqq \forall\,v\in\cV;
\]
\[\lb{gmep}
\wt\P_{j,n}(\vt,x)=e^{i\lan x_\A,\,\vt\ran}\Phi_{j,n}(\vt,x), \qqq \forall\,x\in \cG;
\]
where the functions $\vp_n(\vt,v)$ and $\Phi_{j,n}(\vt,x)$
are $\G$-periodic with respect to $v$ and $x$, respectively, and $x_\A\in\R^d$ is the coordinate vector of $x\in\R^d$ with respect to the basis $\A=\{a_1,\ldots,a_d\}$ of the lattice $\G$ (see \er{cola}).
\end{proposition}

\begin{remark}
1) In \er{gde}, \er{gme} we identify the vertices and edges of the fundamental graph $\cG_*$ with some vertices and edges of the periodic graph $\cG$. For more details see the proof of this proposition.

2) From \er{gme} and \er{ms8} it follows that all generalized eigenfunctions of the absolutely continuous spectrum of the Laplacian on the periodic metric graph $\cG$ are uniformly bounded in $L^\iy(\cG)$.

3) Formulas \er{ms6}, \er{gde} and \er{gme} give that the eigenfunctions of the absolutely continuous spectrum of the discrete and metric Laplacians are connected by the same relation \er{ms6} as the eigenfunctions of their fiber operators.

4) Eigenfunctions of the discrete and metric Laplacians corresponding to flat bands are connected with eigenfunctions of their fiber operators by the Floquet transforms (for more details see Proposition \ref{TLBS}).
\end{remark}

\subsection{Scattering on metric graphs} We consider the Schr\"odinger operator $H=H_0+Q$ on $L^2(\cG)$, where $H_0=\D_M$ and the potential $Q\in L^2(\cG)\cap L^1(\cG)$ is real. Here $L^{1}(\cG)$ is the space of all functions $f=(f_\be)_{\be\in\cE}$ on $\cG$ equipped with the norm $\|f\|_{L^{1}(\cG)}=\sum_{\be\in \cE}\|f_{\be}\|_{L^{1}(\be)}$. Let $\bB_1$ and $\bB_2$ be the trace and the Hilbert-Schmidt class equipped with the norm
$\|\cdot \|_{\bB_1}$ and $ \|\cdot \|_{\bB_2}$, respectively.

In Theorem \ref{TSD} we show that $QR_0(k)\in \bB_2$ for all $k\in\C_+$. This implies that the operator $H$ is self-adjoint on $\mD(H_0)$. For any $k\in \C_+$ we put
\[\label{SY}
\begin{array}{ll}
R_0(k)=(H_0-k^2)^{-1}, \qqq & R(k)=(H-k^2)^{-1},\\[6pt]
Y_0(k)=|Q|^{1/2}\,R_0(k)\,Q^{1/2}, \qqq & Q^{1/2}=|Q|^{1/2}\sign Q.
\end{array}
\]
Below we show that the operator $Y_0(k)$ belongs to the trace class. Thus, we can define the Fredholm determinant $D$ by
\[\lb{defFD}
D(k)=\det (I+Y_0(k)),\qqq k\in \C_+.
\]
We describe our third  main result about scattering on metric graphs.

\begin{theorem}\lb{TSD}
Let $Q\in L^2(\cG)$ be real. Then $QR_0(k)\in \bB_2$ for all $k\in\C_+$ and
the operator $H=H_0+Q$ is self-adjoint on $\mD(H_0)$. Let, in addition, $Q\in L^1(\cG)$. Then
\[\lb{Sc1}
R(k)-R_0(k), \; Y_0(k)\in \bB_1\qqq \forall\, k\in \C_+,
\]
and the determinant $D(k)=\det (I+Y_0(k))$ is well-defined and is analytic in $\C_+$ and the limit $D(k+i0)$ exists for almost all $k\in\R$.  Furthermore,  the wave operators
\[\lb{Sc2}
W_\pm=\slim\,e^{itH}e^{-itH_0}P_{ac}(H_0) \qqq \textrm{ as }\qqq t\to\pm\infty,
\]
exist and are complete, i.e., the range of $W_\pm$ is equal to $\mH_{ac}(H)$ and $\s_{ac}(H)=\s_{ac}(H_0)$. Moreover, the $S$-operator given by
\[\lb{Sc3}
S =W_+^*W_-
\]
is unitary on $\mH_{ac}(H_0)$ and the corresponding $S$-matrix $S(k)$ (defined by \er{DLH0}) for almost all $k^2\in\s_{ac}(H_0)$ satisfies
\[\lb{Sc4}
S(k)=I_k-2\pi iA(k), \qqq A(k)\in \bB_1,
\]
\[\lb{Sc5}
\det S(k) ={\ol D(k+i0)\/D(k+i0)}\,.
\]
\end{theorem}

\begin{remark}
1) Consider the Schr\"odinger operator $\wt H=-\D+\wt Q$ on $L^2(\R^d)$, $d\ge1$, where $-\D$ is the standard Laplacian in $\R^d$ and $\wt Q\in C_0^{\iy}(\R^d)$ is a real potential. Define the operator $\wt Y_0(k)=|\wt Q|^{1/2}(-\D-k^2)^{-1}\wt Q^{1/2}$, $k\in \C_+$.

Firstly, if $d=1$, then it is well known that the corresponding operator $\wt Y_0(k)\in\bB_1$ and so the Fredholm determinant $\det (I+\wt Y_0(k))$ is well defined. Recall that the Fredholm determinant is the basic function to study trace formulas, spectral shift functions, etc. in the general framework.

Secondly, if $d\ge 2$, then it is well known that the corresponding operator $\wt Y_0(k)$ is not trace class and so the Fredholm determinant $\det (I+\wt Y_0(k))$ is not defined.  Thus, we need an essential modification. For example,  Newton \cite{N77} defined the modified Fredholm determinant by
$$
\wt D(k)=\det\big[(I+\wt Y_0(k))e^{-\wt Y_0(k)}\big],\qqq k\in \C_+,\qq \textrm{ if } \qq d=3.
$$
Note that the case $d>3$ is more complicated.

Thus, in contrast to the Schr\"odinger operator $\wt H$ on $\R^d$ in the case of the Schr\"odinger operator $H=\D_M+Q$ on metric graphs the operator $Y_0(k)$ is trace class and the corresponding Fredholm determinant $D(k)=\det (I+Y_0(k))$ is well defined for all $k\in\C_+$ and for any dimension $d\ge 1$. Then the metric case (for any dimension $d\ge 1$) is very similar to the case of the Schr\"odinger operator $\wt H$ on $\R^1$. This fact is very important for spectral theory of Schr\"odinger operators on periodic metric graphs.

2) In fact the identity \er{Sc5} is the  so-called Birman-Krein \cite{BK62} formula for the $S$-matrix in the case of scattering on metric graphs.

3) Recall that if the scattering operator $S=W_+^*W_-$ is unitary on $\mH_{ac}(H_0)$, then the operators $H_0$ and $S$ commute and thus are simultaneously diagonalizable:
\[\lb{DLH0}
\mH_{ac}(H_0)=P_{ac}(H_0)L^2(\cG)=\int_\s^\oplus \mH_\l d\l,\qqq H_0=\int_\s^\oplus\l I_\l d\l,\qqq
S=\int_\s^\oplus S(\l)d\l;
\]
here $\s=\s_{ac}(H_0)$, $I_\l$ is the identity in the fiber space $\mH_\l$, and $S(\l)$ is the scattering matrix acting in $\mH_\l$ for the pair $H_0,H$.
\end{remark}

The paper is organized as follows. In Section \ref{Sec2} we prove Theorem \ref{T1} about the direct integral decomposition for the metric Laplacian $\D_M$ on periodic graphs. Section \ref{Sec3} is devoted to the eigenfunctions of the fiber metric Laplacian $\D_M(\vt)$. Here we prove Theorem \ref{TT1} about the eigenfunctions of $\D_M(\vt)$ corresponding to the Dirichlet eigenvalues and Theorem \ref{TT2} about the connection between the eigenfunctions of the fiber metric $\D_M(\vt)$ and  discrete $\D(\vt)$ Laplacians. Here we also prove Proposition \ref{TLEF} about the connection between generalized eigenfunctions of the absolutely continuous spectrum of the discrete and metric Laplacians and the eigenfunctions of their fiber operators. Section \ref{Sec4} deals with the periodic metric Laplacian perturbed by real integrable potentials. Here we prove Theorem \ref{TSD} about scattering on metric graphs.
In Sections \ref{Sec5} -- \ref{Sec7} we consider some examples of periodic graphs: the $d$-dimensional lattice, the hexagonal lattice and the stanene lattice and give detailed descriptions of the spectrum of the metric Laplacian on these graphs including descriptions of the eigenfunctions of the fiber operators.

%****************************************************************

\section {\lb{Sec2}Direct integral for metric Laplacians}
\setcounter{equation}{0}

\subsection{Edge indices} We define an {\it edge index}, which was introduced in \cite{KS14a}. The indices are important to study the spectrum of the Laplacians and Schr\"odinger operators on periodic graphs, since the fiber operators are expressed in terms of the indices of the fundamental graph edges (see \er{l2.15} for the discrete fiber Laplacian and \er{di1} -- \er{di2} for the metric one).

For any vertex $v\in\cV$ of the $\G$-periodic graph $\cG$ the following unique representation holds true:
\[
\lb{Dv} v=v_0+[v], \qq \textrm{where}\qq v_0\in\cV_0=\cV\cap\Omega,\qquad [v]\in\G,
\]
$\Omega$ is the fundamental cell of the lattice $\G$ defined by \er{fuce}.
In other words, each vertex $v\in\cV$ can be obtained from a vertex $v_0\in\cV_0$ by the shift by a vector $[v]\in\G$. We call $[v]$ the \emph{integer part of the vertex $v$}.

We recall that $\{a_1,\ldots,a_d\}$ is the basis of the lattice $\G$. For each $x\in\R^d$ we introduce the vector $x_\A\in\R^d$ by
\[\lb{cola}
x_\A=(t_1,\ldots,t_d), \qqq \textrm{where} \qq x=\textstyle\sum\limits_{s=1}^dt_sa_s.
\]
In other words, $x_\A$ is the coordinate vector of $x$ with respect to the basis $\A=\{a_1,\ldots,a_d\}$ of the lattice $\G$.

For any oriented edge $\be=(u,v)\in\cA$ we define the {\bf edge index}
$\t(\be)$ as the vector of the lattice $\Z^d$ given by
\[
\lb{in}
\t(\be)=[v]_\A-[u]_\A\in\Z^d,
\]
where $[v]\in\G$ is the integer part of the vertex $v$, and the vector $[v]_\A\in\Z^d$ is defined by \er{cola}. Edge indices depend on the embedding of $\cG$ into $\R^d$ and on the choice of the basis $a_1,\ldots,a_d$ of the lattice $\G$. We note that edges connecting vertices inside the fundamental cell $\Omega$ have zero indices.

\begin{example}
We consider the periodic graph $\cG$ shown in Fig.\ref{ff.0.11}\emph{a}. The periods $a_1,a_2$ of the graph and the fundamental cell $\Omega$ are also shown in the figure. The index of the edge $(v_1,v_3+a_1)$ is equal to $(1,0)$, since $[v_1]=0$ and $[v_3+a_1]=a_1$. The edge $(v_1,v_4)$ has zero index.
\end{example}

On the $\G$-periodic graph $\cG=(\cV,\cE)$ we define two surjections
\[\lb{sur}
\gf_\cV:\cV\rightarrow \cV_*=\cV/\G, \qqq \gf_\cA:\cA\rightarrow\cA_*=\cA/\G,
\]
which map each vertex $v\in\cV$ and each oriented edge $\be\in\cA$ of $\cG$ to their equivalence classes $\gf_\cV(v)$ and $\gf_\cA(\be)$, respectively, which are a vertex and an oriented edge of the fundamental graph $\cG_*=(\cV_*,\cE_*)$.

For each oriented edge $\be_*\in\cA_*$ of the fundamental graph $\cG_*$ we define the edge index $\t(\be_*)\in\Z^d$ by:
\[
\lb{inf}
\t(\bf e_*)=\t(\be) \qq \textrm{ for some $\be\in\cA$ \; such that }  \; \be_*=\gf_\cA(\be), \qqq \be_*\in\cA_*.
\]
Indices of the fundamental graph edges are induced by indices of the periodic graph edges and uniquely determined by \er{inf}, since
$$
\t(\be+a)=\t(\be),\qqq \forall\, (\be,a)\in\cA \ts\G.
$$

Denote by $\cC$ the cycle space of the fundamental graph $\cG_*$. We define the vector-valued \emph{flux function} $\Phi_\t:\cC\to\R^d$ as follows:
\[\lb{mafla}
\Phi_\t(\mathbf{c})=\sum_{\be\in\mathbf{c}}\t(\be), \qqq \mathbf{c}\in\cC.
\]
It is known (see Propositions 4.1.i and 4.2 in \cite{KS19}) that the image of the flux function $\Phi_{\t}$ is the lattice $\Z^d$:
\[\lb{prao}
\Phi_{\t}(\cC)=\Z^d.
\]

\subsection{Fiber Laplacians.}
The normalized Laplacian $\D$, given by \er{DOL}, on the periodic graph $\cG$ has the standard decomposition into the constant fiber direct integral \er{raz}, where the fiber Laplacian $\D(\vt)$ has the form
\[\lb{l2.15}
\big(\D(\vt)f\big)(v)=-\,\frac1{\vk_v}\sum_{\be=(v,\,u)\in\cA_*} e^{i\lan\t(\be),\,\vt\ran}f(u), \qqq f\in\ell^2(\cV_*),\qqq v\in\cV_*,
\]
see Theorem 2.2 in \cite{KS16b}. We recall that $\vk_v$ is the degree of the vertex $v$, $\lan\cdot\,,\cdot\ran$ denotes the standard inner product in $\R^d$, and $\t(\be)$ is the index of the edge $\be\in\cA_*$ defined by \er{in}, \er{inf}.

\medskip

We identify the vertices of the fundamental
graph $\cG_*=(\cV_*,\cE_*)$ with the corresponding vertices of the periodic graph
$\cG=(\cV,\cE)$ from the set $\cV_0=\cV\cap\Omega$. We also identify an edge $\be_*=(u_*,v_*)\in\cE_*$ of the fundamental graph $\cG_*$ with an index $\t(\be_*)$ with the edge $\be=(u,v+\t(\be_*))\in\cE$ of the periodic graph $\cG$, where $u,v\in\cV_0$ such that $u_*=\gf_\cV(u)$, $v_*=\gf_\cV(v)$, and $\gf_\cV$ is defined in \er{sur}. We note that $\t(\be)=\t(\be_*)$.

Let $a(m)\in\G$, $m\in\Z^d$, be defined by
\[\lb{aofm}
a(m)=\textstyle\sum\limits_{s=1}^dm_sa_s, \qqq m=(m_1,\ldots,m_d)\in\Z^d.
\]

\no {\bf Proof of Theorem \ref{T1}.} Denote by $L_{com}^2(\cG)$ the set of all compactly supported functions $h\in L^2(\cG)$. Standard arguments (see pp. 290--291 in \cite{RS78}) yield that $\mU$ given by \er{5001} is well defined on $L_{com}^2(\cG)$ and has a unique extension to a unitary operator.
For $h\in L_{com}^2(\cG)$ the sum \er{5001} is finite and using the identity  $\cE=\big\{\be_*+a: (\be_*,a)\in\cE_*\ts\G\big\}$ we have
\begin{multline*}
\|\mU h\|^2_{\mH}=\int_{\T^d}\|(\mU h)(\vt,\cdot )\|_{L^2(\cG_*)}^2{d\vt\/(2\pi)^d}\\
=\sum_{\be_*\in\cE_*}\int_0^1\int_{\T^d}\bigg(\sum\limits_{m\in\Z^d}e^{-i\lan m,\,\vt\ran }h_{\be_*+a(m)}(t)\bigg)\bigg(\sum\limits_{m'\in\Z^d}e^{i\lan m',\,\vt\ran }\bar{h}_{\be_*+a(m')}(t)\bigg){d\vt\/(2\pi)^d}\;dt\\
=\sum_{\be_*\in\cE_*}\int_0^1\sum_{m,m'\in\Z^d}\lt(h_{\be_*+a(m)}(t)\bar h_{\be_*+a(m')}(t)\int_{\T^d}e^{-i\lan m-m',\,\vt\ran}{d\vt\/(2\pi)^d}\rt)\;dt\\
=\sum_{\be_*\in\cE_*}\int_0^1\sum_{a\in\G}\big|h_{\be_*+a}(t)\big|^2\;dt
=\sum_{\be\in\cE}\int_0^1|h_\be(t)|^2\,dt=\|h\|_{L^2(\cG)}^2.
\end{multline*}
Thus, $\mU$ is  well defined on $L_{com}^2(\cG)$ and has a unique isometric extension. In order to prove that $\mU$ is onto $\mH$ we compute $\mU^*$. Let
$g=\big(g(\cdot,x_*)\big)_{x_*\in\cG_*}\in\mH$, where $g(\cdot,x_*)
:\T^d\to\C$. We define
\begin{equation}\label{rep}
(\mU^*g)(x)=\int_{\T^d}e^{i\lan m,\vt\ran }g(\vt,x_*){d\vt\/(2\pi)^d}\,,\qquad
 x=x_*+a(m)\in\cG,
\end{equation}
where $(x_*,m)\in\cG_*\ts\Z^d$ are uniquely defined.  A direct computation
gives  that this is indeed the formula for the adjoint of $\mU$.
Moreover, the Parseval identity for Fourier series gives
\begin{multline*}
\|\mU^*g\|^2_{L^2(\cG)}=\sum_{\be\in
\cE}\big\|(\mU^*g)_\be\big\|_{L^2(0,1)}^2=\sum_{\be_*\in
\cE_*}\sum_{a\in\G}\int_0^1\big|(\mU^*g)_{\be_*+a}(t)\big|^2\,dt\\=
\sum_{\be_*\in\cE_*}\sum_{m\in\Z^d}\int_0^1\bigg|\int_{\T^d}e^{i\lan m,\vt\ran
}g_{\be_*}(\vt,t){d\vt \/(2\pi)^d}\bigg|^2\,dt\\
=\sum_{\be_*\in\cE_*}\int_0^1\int_{\T^d}\big|g_{\be_*}(\vt,t)\big|^2{d\vt \/(2\pi)^d}\,dt=\int_{\T^d}\big\|g(\vt,\cdot)\big\|_{L^2(\cG_*)}^2{d\vt \/(2\pi)^d}=\|g\|_\mH^2.
\end{multline*}
Further, for $h\in L_{com}^2(\cG)$ and $x\in\be\in\cE_*$, $x\notin\cV_*$ we obtain
\begin{multline}\label{ext1}
(\mU\D_M h)(\vt,x)=\sum_{m\in\Z^d}e^{-i\lan m,\vt\ran}(\D_M
h)\big(x+a(m)\big)\\=\sum_{m\in\Z^d}e^{-i\lan m,\vt\ran }
h''_{\be+a(m)}(x+a(m))=\D_M(\vt)(\mU h)(\vt,x).
\end{multline}
Let $v\in\cV_*$ and $\be_1,\be_2,\be\in I_*(v)$. Denote by $\d_1=\d(\be_1,v)$, $\d_2=\d(\be_2,v)$ and $\d=\d(\be,v)$, where $\d$ is defined by \er{tev}. Then, using \er{Dom1}, we have
\begin{multline}\label{ext2}
e^{-i\d_1\lan \t(\be_1),\vt\ran}(\mU h)_{\be_1}(\vt,\d_1)=e^{-i\d_1\lan \t(\be_1),\vt\ran}\sum\limits_{m\in\mathbb{Z}^d}e^{-i\lan m,\vt\ran }h_{\be_1+a(m)}(\d_1)\\=e^{-i\d_2\lan \t(\be_2),\vt\ran}\sum\limits_{m\in\mathbb{Z}^d}e^{-i\lan m+\d_1\t(\be_1)-\d_2\t(\be_2),\vt\ran }h_{\be_1+a(m)}(\d_1)\\=e^{-i\d_2\lan \t(\be_2),\vt\ran}\sum\limits_{m\in\mathbb{Z}^d}e^{-i\lan m+\d_1\t(\be_1)-\d_2\t(\be_2),\vt\ran }h_{\be_2+a(m+\d_1\t(\be_1)-\d_2\t(\be_2))}(\d_2)\\=e^{-i\d_2\lan\t(\be_2),\vt\ran}
\sum\limits_{m\in\mathbb{Z}^d}e^{-i\lan m,\vt\ran }h_{\be_2+a(m)}(\d_2)=e^{-i\d_2\lan\t(\be_2),\vt\ran}(\mU h)_{\be_2}(\vt,\d_2),
\end{multline}
and
\begin{multline}\label{ext3}
\sum\limits_{\be\in I_\ast(v)}(-1)^\d e^{-i\d\lan\t
(\be),\,\vt\ran}\,(\mU h)'_{\be}(\vt,\d)=\sum\limits_{\be\in I_*(v)}(-1)^\d e^{-i\d\lan\t
(\be),\,\vt\ran}
\sum\limits_{m\in\mathbb{Z}^d}e^{-i\lan m,\vt\ran }h'_{\be+a(m)}(\d)\\
=\sum\limits_{\be\in I_*(v)}(-1)^\d
\sum\limits_{m\in\mathbb{Z}^d}e^{-i\lan m+\d\t
(\be),\vt\ran }h'_{\be+a(m)}(\d)=
\sum\limits_{\be\in I_*(v)}(-1)^\d
\sum\limits_{m\in\Z^d}e^{-i\lan m,\vt\ran }h'_{\be+a(m-\d\t
(\be))}(\d)\\=\sum\limits_{m\in\mathbb{Z}^d}e^{-i\lan m,\vt\ran }\sum\limits_{\be\in I_*(v)}(-1)^\d h'_{\be+a(m-\d\t
(\be))}(\d)=\sum\limits_{m\in\mathbb{Z}^d}e^{-i\lan m,\vt\ran }\sum\limits_{\be\in I(v+a(m))}(-1)^\d
h'_{\be+a(m)}(\d)=0.
\end{multline}
The identities (\ref{ext1}) -- \er{ext3} yield \er{Mraz}, \er{FBC}, \er{di2}. \qq
\BBox

\medskip

\begin{corollary}\label{cor1}
The Laplacian $\D_M$  on $L^2(\cG)$ has a decomposition into a
constant fiber direct integral of the form \er{Mraz}, where the fiber operator
$\wt\D_M(\vt)$  on $L^2(\cG_*)$ is defined by
\[\lb{di3}
\textstyle (\wt\D_M(\vt) y\,)_{\be}=-\big({d\/dt}+i\,\lan\t
(\be),\vt\ran\big)^2y_\be,\qqq y=(y_\be)_{\be\in\cE_*},
\]
where $(y\,''_\be)_{\be\in\cE_*}\in L^2(\cG_*)$,
$y$ is continuous on $\cG_*$ and satisfies the boundary conditions:
\[\lb{di4}
\textstyle \sum\limits_{\be\in
I_*(v)}(-1)^{\d(\be,v)}\big({d\/dt}+i\,\lan\t
(\be),\vt\ran\big)\,y_\be\big(\d(\be,v)\big)=0 \qqq \forall\, v\in
\cV_*,
\]
or, that is the same,
\[
\lb{di4'}
\begin{aligned}
\textstyle \sum\limits_{\be\in I_\ast(v)}(-1)^{\d(\be,v)}y'_\be\big(\d(\be,v)\big)
=y(v)\,i\sum\limits_{\be\in I_\ast(v)}(-1)^{\d(\be,v)}
\lan\t
(\be),\vt\ran, \qqq \forall v\in \cV_*,\\
y(v)=y_\be\big(\d(\be,v)\big),\qqq \forall\, \be\in I_\ast(v).
\end{aligned}
\]
\end{corollary}

\no \textbf{Proof.} After the gauge transformation
$
\wt y_{\be}(t)=e^{-it\lan\t
({\bf e}),\,\vt\ran}y_{\be}(t)$,  for all $\be\in\cE_\ast,
$
identities \er{di1} -- \er{di2} take the form \er{di3} -- \er{di4}. \qq
\BBox

\section{\lb{Sec3} Eigenfunctions of metric Laplacians} \setcounter{equation}{0}

In this section we prove Theorems \ref{TT1} and \ref{TT2} about eigenfunctions of the fiber metric Laplacian $\D_M(\vt)$ and Proposition \ref{TLEF} about the connection between eigenfunctions of the absolutely continuous spectrum of the discrete and metric Laplacians and the eigenfunctions of their fiber operators.

\no \textbf{Proof of Theorem \ref{TT1}.} \emph{i}) -- \emph{iii}) Let $\vt\in\T^d\sm\{0\}$ and $j\in\N$. We construct all
eigenfunctions $\P_{j}^0=\big(\P_{j,\be}^0\big)_{\be\in\cE_\ast}$
of the operator $\D_M(\vt)$ corresponding to the eigenvalue $(\pi
j)^2$. These eigenfunctions have the form
\begin{equation}\label{fb9}
\textstyle \P_{j,\be}^0(\vt,t)=\P_{j,\be}^0(\vt,0)\,\cos(\pi jt)+
{\P_{j,\be}^0}'(\vt,0)\,{\sin(\pi jt)\/\pi j}\,, \qqq \forall\, t\in [0,1],
\end{equation}
which yields
\[\lb{fbdi8}
\textstyle \P_{j,\be}^0(\vt,1)=(-1)^j\,\P_{j,\be}^0(\vt,0),\qqq
{\P_{j,\be}^0}'(\vt,1)=(-1)^j\,{\P_{j,\be}^0}'(\vt,0).
\]
Let $\bc\in\cC$, where $\cC$ is the cycle space of the fundamental graph $\cG_*$. Then, using the quasi-periodic condition \er{FBC} and the first identity in \er{fbdi8}, we obtain
\[
\lb{cyc}
\P_{j,\be}^0(\vt,0)\big((-1)^{jl_\bc}\,e^{-i\lan\Phi_\t(\mathbf{c}),\,\vt\ran}-1\big)=0,
\qqq \forall\, \be\in\bc,
\]
where $l_\bc$ is the length of the cycle $\bc$, i.e., the number of its edges, and $\Phi_\t(\mathbf{c})$ is defined by \er{mafla}. Since $\Phi_{\t}(\cC)=\Z^d$ and $\vt\neq0$, from \er{cyc} it follows that $\P_{j,\be}^0(\vt,0)=0$ for each $\be\in\bc$, $\bc\in\cC$. Due to the connectivity of the fundamental graph, the quasi-periodic condition \er{FBC} and the first identity in \er{fbdi8} give $ \P_{j,\be}^0(\vt,0)=0$ for all $\be\in\cE_*$. Thus, the eigenfunctions \er{fb9} take the form
\begin{equation}
\label{fb9'} \textstyle \P_{j,\be}^0(\vt,t)= {\P_{j,\be}^0}'(\vt,0)\,{\sin(\pi jt)\/\pi j}\,, \qqq \forall \, \be\in\cE_\ast.
\end{equation}
These functions must satisfy the quasi-periodic condition \er{di2}, which, using the second identity in \er{fbdi8}, can be rewritten in the form
\[
\lb{di2'}
\sum\limits_{\be=(v,u)\in\cE_\ast}{\P_{j,\be}^0}'(\vt,0)-(-1)^j
\sum\limits_{\be=(u,v)\in\cE_\ast}e^{-i\lan\t ({\bf
e}),\,\vt\ran}\,{\P_{j,\be}^0}'(\vt,0)=0, \qqq \forall\, v\in\cV_*.
\]
This is a homogeneous system of $\n=\#\cV_*$ linear equations with $\n_1=\#\cE_*$ variables \linebreak ${\big({\P_{j,\be}^0}'(\vt,0)\big)_{\be\in\cE_*}}$. The $\n\ts\n_1$ coefficient matrix $D_j(\vt)=\{D_{v,\be}^{(j)}(\vt)\}_{v\in\cV_*\atop \be\in \cE_*}$ of the system is given by
\[\lb{MaD}
D_{v,\be}^{(j)}(\vt)=\left\{
\begin{array}{cl}
  (-1)^{j+1}e^{-i\lan\t(\be),\,\vt\ran }, & \textrm{if $v$ is the terminal vertex of $\be$} \\
  1, & \textrm{if $v$ is the initial vertex of $\be$}\\
  1-(-1)^je^{-i\lan\t(\be),\,\vt\ran } \qq &
  \textrm{if $\be$ is a loop in the vertex $v$}\\
  0, & \textrm{otherwise}
\end{array}\right..
\]
The number of linearly independent solutions of the system \er{di2'}
is $\big(\n_1-\rank D_j(\vt)\big)$. We will show that $\rank D_j(\vt)=\n$.

Let $D_{v\,\centerdot}^{(j)}(\vt)$ denote the row of the matrix
$D_j(\vt)$,  corresponding to the vertex $v\in\cV_*$. In order to
show that the rows of this matrix are linearly independent, we
consider their linear combination with coefficients $\a_v$,
$v\in\cV_*$:
\[\lb{lic}
\sum_{v\in\cV_*}\a_v\,D_{v\,\centerdot}^{(j)}(\vt)=0.
\]
From this, using the form \er{MaD} of the matrix $D_j(\vt)$, we
obtain
\[\lb{aa'}
\a_u=\a_v(-1)^{j}e^{-i\lan\t(\be),\,\vt\ran},\qqq \forall\,
\be=(u,v)\in\cE_*.
\]
Let $\bc\in\cC$. Then for each vertex $v$ of the cycle $\bc$, the identities \er{aa'} give
\[\lb{cyc'}
\a_v\big((-1)^{jl_\bc}e^{i\lan \Phi_\t(\mathbf{c}),\,\vt\ran}-1\big)=0.
\]
Since $\vt\neq0$ and $\Phi_{\t}(\cC)=\Z^d$,  we conclude that $\a_v=0$ for each vertex $v$ of each cycle $\bc$. Due to the connectivity of $\cG_*$ and the identity
\er{aa'}, all coefficients $\a_v=0$ in \er{lic}. Thus, the rows of
the matrix $D_j(\vt)$ are linearly independent and $\rank
D_j(\vt)=\n$. Then the number of linearly independent
solutions of the system \er{di2'} is equal to $\n_1-\n$. Using this and the identity \er{Benu}, we obtain that the eigenvalue $(\pi j)^2$ has multiplicity $\b-1$.

From \er{fb9'} it follows that
\[
\|\P_j^0(\vt,\cdot)\|_{L^2(\cG_*)}^2={1\/2(\pi j)^2}\sum_{\be\in\cE_*}
\big|{\P_{j,\be}^0}'(\vt,0)\big|^2.
\]
Choosing the normalized solutions
$\big(X_{j,s,\be}(\vt)\big)_{\be\in\cE_\ast}$,
$s\in\N_{\b-1}$, of the system \er{di2'} or, equivalently, of
the system \er{sq5}, we obtain that the normalized eigenfunctions
$\P_{j,s}^0=\big(\P_{j,s,\be}^0\big)_{\be\in\cE_\ast}$, $s\in
\N_{\b-1}$, have the form \er{eif}.

The identity $X_{j,s,\be}(\vt)=X_{j+2,s,\be}(\vt)$,
$(s,\be)\in\N_{\b-1}\ts\cE_\ast$,  follows from the form of the
system \er{sq5}. The inequality \er{supv} is a direct consequence of
\er{eif}. \qq \BBox

\begin{remark}\lb{De0}
The Floquet operator $\D_M(0)$ is the metric
Laplacian on  the fundamental graph $\cG_*=(\cV_*,\cE_*)$. It is known (see
\cite{B85}, \cite{C97}, \cite{LP08}), that

1) the normalized eigenfunctions $\P_{j,s}^0=\big(\P_{j,s,\be}^0\big)_{\be\in\cE_\ast}$, $s\in
\N_{\b}$ of the operator
$\D_M(0)$ corresponding to the eigenvalue  $(\pi j)^2$, $j$ is even, have the form \er{eif} as $\vt=0$,
where $\big(X_{j,s,\be}(0)\big)_{\be\in\cE_\ast}$ is a normalized
solution of the system of $\n$ equations
\[
\sum_{\substack{\be\in\cE_\ast, \\ \textrm{ starting at
$v$}}}x_\be-\sum_{\substack{\be\in\cE_\ast,\\ \textrm{ ending
at $v$}}}\,x_\be=0, \qqq \forall\, v\in\cV_*,
\]
where $\n=\#\cV_*$,
one equation for each $v\in\cV_*$. This system has rank $\n-1$
and thus there exist $\b$ linearly independent normalized
solutions $\big(X_{j,s,\be}(0)\big)_{\be\in\cE_\ast}$,
$s\in\N_{\b}$, where $\b=\#\cE_*-\#\cV_*+1$ is the Betti number of the fundamental graph $\cG_*$. In this case there exists an additional
normalized  eigenfunction
$\P_{j}^0=\big(\P_{j,\be}^0\big)_{\be\in\cE_\ast}$, having the form
\begin{equation}
\P_{j,\be}^0(0,t)=\sqrt{2\/\n_1}\cos(\pi jt),
\end{equation}
where $\n_1=\#\cE_*$. This additional eigenfunction is related to the eigenfunction of the
discrete Laplacian $\D(0)$  with the eigenvalue $-1$,

2) the normalized eigenfunctions
$\P_{j,s}^0=\big(\P_{j,s,\be}^0\big)_{\be\in\cE_\ast}$ corresponding
to the eigenvalue $(\pi j)^2$, $j$ is odd, also have the form
\er{eif} as $\vt=0$, where $\big(X_{j,s,\be}(0)\big)_{\be\in\cE_*}$
is a normalized solution of the system of $\n$ equations
\[
\sum_{\substack{\be\in\cE_*, \\ \textrm{ starting at
$v$}}}x_\be+\sum_{\substack{\be\in\cE_*, \\ \textrm{ ending
at $v$}}}\,x_\be=0, \qqq \forall\, v\in\cV_*.
\]
If $\cG_*$ is bipartite, then this system has rank $\n-1$. If
$\cG_*$ is not bipartite,  this system has rank $\n$. For a
bipartite graph there exists an additional normalized eigenfunction
$\P_{j}^0=\big(\P_{j,\be}^0\big)_{\be\in\cE_*}$, having the form
\begin{equation}
\P_{j,\be}^0(0,t)=\pm\sqrt{2\/\n_1}\cos(\pi jt),
\end{equation}
where the sign is chosen such that the eigenfunction has alternative
sign in adjacent vertices. This additional eigenfunction is related
to the eigenfunction of the discrete Laplacian $\D(0)$ with the
eigenvalue $1$.
\end{remark}

In the following proposition we describe the connection between
eigenfunctions  and corresponding eigenvalues (except for the
Dirichlet eigenvalues $\s_D$) of the Floquet operators for metric
and discrete Laplacians.

\begin{proposition}\lb{proC}
Let $\vt\in\T^d$. Then the following statements hold true.

i) If $\P(\vt)=\big(\P_{\be}(\vt,t)\big)_{\be\in\cE_*}\in L^2(\cG_*)$
is  an eigenfunction of the operator $\D_M(\vt)$ with an eigenvalue
$E(\vt)\notin\s_D$, then the function
$\p(\vt,\cdot)\in\ell^2(\cV_*)$, defined by
\[\lb{psi}
\p(\vt,v)=e^{-i\d(\be,v)\lan\t
(\be),\,\vt\ran}\P_{\be}(\vt,\d(\be,v)), \qqq \be\in I_*(v),
\]
is an eigenfunction of the operator $\D(\vt)$ with the eigenvalue $\l(\vt)=-\cos\sqrt{E(\vt)}$.

ii) Conversely, if $\p(\vt,\cdot)\in\ell^2(\cV_*)$ is an
eigenfunction   of $\D(\vt)$ with an eigenvalue $\l(\vt)\in(-1,1)$,
then for each $j\in\N$ the function
$\P_j(\vt)=\big(\P_{j,\be}(\vt,t)\big)_{\be\in\cE_*}\in L^2(\cG_*)$,
defined by
\[\lb{ms6'}
\begin{aligned}
\P_{j,\be}(\vt,t)=\textstyle {1\/\sin z_j(\vt)}\,
\big(\p(\vt,u)\sin(z_j(\vt)(1-t))+
\p(\vt,v)e^{i\lan\t(\be),\,\vt\ran}\,\sin(z_j(\vt)\,t)\big), \qqq
\be=(u,v), \\
z_j(\vt)=\ca z(\vt)+\pi j, & \qq j \textrm{ is even} \\
(\pi-z(\vt))+\pi j, & \qq j \textrm{ is odd}
\ac, \qqq z(\vt)=\arccos(-\l(\vt)),
\end{aligned}
\]
is an eigenfunction of $\D_M(\vt)$ with the eigenvalue $E_j(\vt)=z_j^2(\vt)$.
\end{proposition}

\no \textbf{Proof.} \emph{i}) Let $\P(\vt)$ be an eigenfunction of
the operator   $\D_M(\vt)$ with an eigenvalue $E\notin\s_D$,
$E\neq0$, $E\equiv E(\vt)$. The proof for the case $E=0$ is similar.
Then
\begin{equation}\label{9}
\textstyle
\P_{\be}(t)=\P_{\be}(0)\cos(zt)+\P'_{\be}(0)\,{\sin(zt)\/z}\,,  \qqq
z=\sqrt{E},\qqq \P_{\be}(t)\equiv\P_{\be}(\vt,t),
\end{equation}
which yields
\[\lb{di8}
\textstyle \P_{\be}(1)=\P_{\be}(0)\cos z+\P'_{\be}(0)\,{\sin z\/z}\,, \qqq
\P'_{\be}(1)=-\P_{\be}(0)z\sin z+\P'_{\be}(0)\cos z.
\]
From \er{di8} we obtain
\begin{equation}\label{eq3}
\textstyle \P'_{\be}(0)=\frac{z}{\sin
z}\,\big(\P_{\be}(1)-\P_{\be}(0)\cos z\big), \qqq
\P'_{\be}(1)=\frac{z}{\sin z}\,\big(\P_{\be}(1)\cos
z-\P_{\be}(0)\big).
\end{equation}
Using the formulas \er{eq3}, we rewrite the left-hand side of the condition \er{di2} in the form:
\begin{multline*}
\textstyle\sum\limits_{\be\in I_*(v)}(-1)^{\d(\be,v)}e^{-i\d(\be,v)\lan\t
(\be),\,\vt\ran}\,\P'_{\be}\big(\d(\be,v)\big)=
\hspace{-3mm}\sum\limits_{\be=(v,u)\in\cE_*}\hspace{-3mm}\P'_{\be}(0)-
\hspace{-3mm}\sum\limits_{\be=(u,v)\in\cE_*}\hspace{-3mm}e^{-i\lan\t
(\be),\,\vt\ran}\,\P'_{\be}(1)\\
=\textstyle\sum\limits_{\be=(v,u)\in\cE_*}\frac{z}{\sin z}\,
\big(\P_{\be}(1)-\P_{\be}(0)\cos z\big)-\textstyle
\sum\limits_{\be=(u,v)\in\cE_*}e^{-i\lan\t(\be),\,\vt\ran}\,\frac{z}{\sin
z}\,  \big(\P_{\be}(1)\cos z-\P_{\be}(0)\big),
\end{multline*}
for all $v\in\cV_*$. Thus, after dividing by $\frac{z}{\sin z}$\,,
the condition \er{di2} takes the form
\[\lb{ide1}
\textstyle
\sum\limits_{\be=(v,u)\in\cE_*}\big(\P_{\be}(1)-\P_{\be}(0)\,\cos
z\big)
-\sum\limits_{\be=(u,v)\in\cE_*}e^{-i\lan\t(\be),\,\vt\ran}\,\big(\P_{\be}(1)\cos
z-\P_{\be}(0)\big)=0,
\]
for all $v\in\cV_*$. Using the formulas \er{FBC} and \er{psi},
we write the left-hand side of the identity \er{ide1} in the form
\begin{multline*}
\textstyle-\p(\vt,v)\vk_v\cos z+\sum\limits_{\be=(v,u)\in\cE_*}\P_{\be}(1)
+\sum\limits_{\be=(u,v)\in\cE_*}e^{-i\lan\t(\be),\,\vt\ran}\,\P_{\be}(0)\\=
\textstyle-\p(\vt,v)\vk_v\cos z+\sum\limits_{\be=(v,u)\in\cE_*}e^{i\lan\t(\be),\,\vt\ran}\,\p(\vt,u)+
\sum\limits_{\be=(u,v)\in\cE_*}e^{-i\lan\t(\be),\,\vt\ran}\,\p(\vt,u)\\=
\textstyle-\p(\vt,v)\vk_v\cos z+\sum\limits_{\be=(v,u)\in\cA_*}e^{i\lan\t(\be),\,\vt\ran}\,\p(\vt,u).
\end{multline*}
Thus, the condition \er{ide1} takes the form
$$
\p(\vt,v)\cos z=\frac1{\vk_v}\sum\limits_{\be=(v,u)\in\cA_*}e^{i\lan\t
(\be),\,\vt\ran}\,\p(\vt,u),
$$
or, using the formula \er{l2.15},
\[\lb{ms2}
-\textstyle\p(\vt,v)\cos z=\big(\D(\vt)\p(\vt,\cdot)\big)(v).
\]
Thus, $\p(\vt,\cdot)$ is an eigenfunction of the operator $\D(\vt)$
with  the eigenvalue $(-\cos z)$.

\emph{ii}) Conversely, let $\p(\vt,\cdot)$ be an eigenfunction of
the operator   $\D(\vt)$ with the eigenvalue $\l(\vt)\in(-1,1)$. For
the function $\P_j(\vt)=\big(\P_{j,\be}(\vt,t)\big)_{\be\in\cE_*}$
defined by \er{ms6'} we have
\[\lb{fder}
\P'_{j,\be}(t)=\textstyle {z_j\/\sin
z_j}\,\big(-\p(u)\cos(z_j(1-t))+\p(v)
e^{i\lan\t(\be),\,\vt\ran}\,\cos(z_jt)\big),
\]
$$
\P''_{j,\be}(t)=\textstyle -{z_j^2\/\sin
z_j}\,\big(\p(u)\sin(z_j(1-t))+\p(v)
e^{i\lan\t(\be),\,\vt\ran}\,\sin(z_jt)\big),
$$
where for the shortness
$$
\P_{j,\be}(t)=\P_{j,\be}(\vt,t),\qqq z_j=z_j(\vt),\qqq \p(u)=\p(\vt,u).
$$
Thus, on each edge $\be\in\cE_*$ the function $\P_j(\vt)$ satisfies the
equation $-\P''_{j,\be}(t)=z_j^2 \P_{j,\be}(t)$. From \er{ms6'} we
obtain that
$$
\P_{j,\be}(0)=\p(u),\qqq \P_{j,\be}(1)=\p(v)e^{i\lan\t(\be),\,\vt\ran}\,,\qqq \be=(u,v),
$$
which yields
$$
e^{-i\d(\be,v)\lan\t(\be),\,\vt\ran}\P_{j,\be}(\d(\be,v))=\p(v),
\qqq \forall\, v\in\cV_*,\qqq \forall\, \be\in I_*(v),
$$
i.e., the condition \er{FBC} holds true. Similarly, from \er{fder} we have
\[\lb{fd0}
\begin{aligned}
&\P'_{j,\be}(0)=\textstyle {z_j\/\sin z_j}\,\big(-\p(u)\cos z_j+\p(v)e^{i\lan\t(\be),\,\vt\ran}\big),\\
&\P'_{j,\be}(1)=\textstyle {z_j\/\sin
z_j}\,\big(-\p(u)+\p(v)\,e^{i\lan\t(\be),\,\vt\ran}\,\cos z_j\big),
\qqq \be=(u,v).
\end{aligned}
\]
Using \er{fd0}, the left-hand side of the condition \er{di2} takes the form
\begin{multline*}
\textstyle \sum\limits_{\be\in
I_*(v)}(-1)^{\d(\be,v)}e^{-i\d(\be,v)\lan\t(\be),\,\vt\ran}\,
\P'_{j,\be}\big(\d(\be,v)\big)=
\hspace{-3mm}\sum\limits_{\be=(v,u)\in\cE_*}\hspace{-3mm}\P'_{j,\be}(0)-
\hspace{-3mm}\sum\limits_{\be=(u,v)\in\cE_*}\hspace{-3mm}e^{-i\lan\t
(\be),\,\vt\ran}\,\P'_{j,\be}(1)\\
=\textstyle\sum\limits_{\be=(v,u)\in\cE_*}{z_j\/\sin z_j}
\big(\p(u)e^{i\lan\t (\be),\,\vt\ran}-\p(v)\cos z_j\big)+
\sum\limits_{\be=(u,v)\in\cE_*}{z_j\/\sin z_j}
\big(\p(u)e^{-i\lan\t(\be),\,\vt\ran}-\p(v)\cos z_j\big)\\
=\textstyle-{z_j\cos z_j\/\sin z_j}\,\vk_v\p(v)+ {z_j\/\sin
z_j}\sum\limits_{\be=(v,u)\in\cA_*}\,\p(u)\,e^{i\lan\t(\be),\,\vt\ran},\qqq \forall\,v\in\cV_*.
\end{multline*}
This, the definition \er{l2.15} of the operator $\D(\vt)$ and the
fact that $\p(\vt,\cdot)$  is an eigenfunction of $\D(\vt)$ with the
eigenvalue $\l(\vt)=-\cos z_j$ yield
$$
\begin{aligned}
\textstyle \sum\limits_{\be\in
I_*(v)}(-1)^{\d(\be,v)}e^{-i\d(\be,v)\lan\t (\be),\,\vt\ran}\,
\P'_{j,\be}\big(\d(\be,v)\big)=-{z_j\vk_v\/\sin
z_j}\big(\p(\vt,v)\cos z_j+\big(\D(\vt)\p(\vt,\cdot)\big)(v)\big)=0
\end{aligned}
$$
for all $v\in\cV_*$, i.e., the condition \er{di2} also holds true. Thus, $\P_j(\vt)$ is an eigenfunction of the operator $\D_M(\vt)$ with the eigenvalue $z_j^2$. \qq \BBox

\medskip

\no \textbf{Proof of Theorem \ref{TT2}.} For all $\vt\in\T^d\sm\{0\}$
the fiber Laplacian $\D(\vt)$ has $\n$ real
eigenvalues $-1<\l_1(\vt)\leq\ldots\leq\l_{\nu}(\vt)<1$.
Due to Proposition \ref{proC} all eigenvalues $E_{j,n}(\vt)=z^2_{j,n}(\vt)$ of the operator $\D_{MV}(\vt)$ are given by \er{egv1} and the eigenfunction $\P_{j,n}(\vt)$ corresponding to the eigenvalue $E_{j,n}(\vt)$ (defined up to a constant factor) has the form \er{ms6}.

Now we will show that $\P_{j,n}(\vt)$ is a normalized eigenfunction. Let for the shortness
$$
c_{j,n}=\cos z_{j,n}(\vt),\qqq s_{j,n}=\sin z_{j,n}(\vt),\qqq
z_{j,n}=z_{j,n}(\vt),\qqq \p_n(u)=\p_n(\vt,u).
$$
Direct integration yields
\begin{multline}\label{unit}
\|\P_{j,n}(\vt)\|_{L^2(\cG_*)}^2=
\sum_{\be\in\cE_*}\int_0^1 \big|\P_{j,n,\be}(\vt,t)\big|^2dt\\
={1\/z_{j,n}s_{j,n}^2} \sum\limits_{\be=(u,v)\in\cE_*}
\Big[\big(z_{j,n}-c_{j,n}s_{j,n}\big)\big(|\p_n(u)|^2+|\p_n(v)|^2\big)
\\
+\big(s_{j,n}- z_{j,n}c_{j,n}\,\big)\big(e^{-i\lan\t(\be),\vt\ran}\p_n(u)\ol \p_n(v)+e^{i\lan\t (\be),\vt\ran} \ol
\p_n(u)\p_n(v)\big)\Big].
\end{multline}
The definition \er{l2.15} of the operator $\D(\vt)$ and the fact that $\p_n(\vt,\cdot)$ is an eigenfunction of $\D(\vt)$ corresponding to $\l_n(\vt)$ yield
\begin{multline}\label{rrr}
\sum\limits_{\be=(u,v)\in\cE_*}\big(e^{-i\lan\t(\be),\vt\ran}\p_n(u)\ol \p_n(v)+e^{i\lan\t(\be),\vt\ran}\ol \p_n(u)
\p_n(v)\big)
\\
= \sum\limits_{\be=(u,v)\in\cA_*}e^{-i\lan\t(\be),\vt\ran}
\p_n(u)\ol\p_n(v) =\sum\limits_{u\in\cV_*}\p_n(u)\sum\limits_{\be=(u,v)\in\cA_*}e^{-i\lan\t
(\be),\vt\ran}\ol \p_n(v)\\
=-\l_n(\vt)\sum\limits_{u\in\cV_*}\vk_u\p_n(u)\ol\p_n(u)=c_{j,n}
\|\p_n(\cdot)\|^2_{\ell^2(\cV_*)}.
\end{multline}
Substituting the identities
$$
\sum\limits_{(u,v)\in\cE_*} \big(|\p_n(u)|^2+|\p_n(v)|^2\big)
=\sum\limits_{v\in\cV_*}\vk_v|\p_n(v)|^2=\|\p_n(\vt,\cdot)\|^2_{\ell^2(\cV_*)}
$$
and \er{rrr} into \er{unit}, we obtain
$$
\|\P_{j,n}(\vt)\|_{L^2(\cG_*)}^2={(z_{j,n}-c_{j,n} s_{j,n})
+(s_{j,n}-z_{j,n}c_{j,n})c_{j,n} \/z_{j,n}\,s_{j,n}^2}\,
\|\p_n(\vt,\cdot)\|^2_{\ell^2(\cV_*)}=\|\p_n(\vt,\cdot)\|^2_{\ell^2(\cV_*)}=1.
$$
The identity \er{ms8} follows from Propositions \ref{T.As} and \ref{T.As1}. \qq
\BBox

\medskip

The first eigenvalue $\l_1(\vt)$ and the corresponding
normalized eigenfunction $\p_1(\vt,\cdot)$ of the fiber Laplacian $\D(\vt)$ have
asymptotics as $\vt=\ve\omega$, $\ve=|\vt\,|\rightarrow0$, $\omega\in\S^{d-1}$:
\[\label{lam}
\begin{aligned}
\l_1(\vt)=\l_1(0)+\ve^2\m(\omega)+O(\ve^3), \qqq \l_1(0)=-1, \qqq
\m(\omega)=\textstyle{1\/2}\,\ddot\l_1(\ve\omega)\big|_{\ve=0},\\
\p_1(\vt,\cdot)=\p_1(0,\cdot)+\ve\p_1^{(1)}+O(\ve^2)
,\qqq \p_1(0,v)=\vk^{-1/2}, \qqq \vk=\sum_{v\in\cV_*}\vk_v,\\
\p_1^{(1)}=\p_1^{(1)}(\omega,\cdot)=\dot\p_1(\ve\omega,\cdot)\big|_{\ve=0},
\end{aligned}
\]
where $\dot g={\pa g\/\pa \ve}$ and  $\S^{d}$ is the $d$-dimensional
sphere.

\begin{proposition}\lb{T.As}
i) Each eigenfunction $\P_{j,1}(\vt)=
\big(\P_{j,1,\be}(\vt,t)\big)_{\be\in\cE_*}$, $j\in \Z_0$, of the operator
$\D_M(\vt)$, $\vt\in\T^d\sm\{0\}$, defined by \er{ms6}, satisfies as
$\be=(u,v)$ and $\ve=|\vt|\to0$
\[\lb{ass8}
\P_{j,1,\be}(\vt,t)= {(-1)^j\sqrt{2}\/\sqrt{\vk}}\;\cos(\pi j_*t)
-{\sin(\pi j_*t)\/\sqrt{\m(\omega)}}\bigg(\p_1^{(1)}(\omega,u)-\p_1^{(1)}(\omega,v)
-{i\lan\t (\be),\omega\ran\/\sqrt{\vk}}\bigg) +O(\ve),
\]
where $j_*=j$, if $j$ is even, and $j_*=j+1$ if $j$ is odd.

ii) Let $\vt\in\T^d\sm\{0\}$ and $z_1(\vt)\leq{\pi\/2}$\,. Then all eigenfunctions $\P_{j,1}(\vt)$, $j\in\Z_0$, of $\D_M(\vt)$ satisfy
\[
\lb{est}
\|\P_{j,1}(\vt)\|_{L^\iy(\cG_*)}< \sqrt{2}\,\bigg(2+{\pi\/\L}+
M\,{|\vt|\/\sin z_1(\vt)}\bigg),
\]
where
\[\lb{MMM}
M=\max_{\be\in\cA_*}|\t(\be)|+{2\/\L}\max_{u\in \cV_*}{1\/\vk_u}\sum_{\be=(u,v)\in\cA_*}|\t(\be)|,
\]
$\L$ is the distance between $\l_1(0)=-1$ and the set $\s\big(\D(0)\big)\sm\big\{\l_1(0)\big\}$.

iii) Let $\vt\in\T^d\sm\{0\}$ and $z_1(\vt)>{\pi\/2}$\,. Then all eigenfunctions $\P_{j,1}(\vt)$, $j\in\Z_0$, of $\D_M(\vt)$ satisfy
\[
\lb{est'}
\|\P_{j,1}(\vt)\|_{L^\iy(\cG_*)}<{2\sqrt{2}\/\sin z_1^+}<\iy,
\]
where $z_1^+$ is defined in \er{sDM1}.

iv) There exists a constant $C_0>0$ such that
\[\lb{sup}
\sup_{\vt\in\T^d\sm\{0\}}{|\vt|\/\sin z_1(\vt)}=C_0<\iy.
\]
\end{proposition}

\no \textbf{Proof.} Let for the shortness
$$
z_{j,1}=z_{j,1}(\vt),\qqq z_1=z_1(\vt),\qqq \p_1(u)=\p_1(\vt,u).
$$

\emph{i}) For the eigenfunction $\P_{j,1}(\vt)=\big(\P_{j,1,\be}(\vt,t)\big)_{\be\in\cE_*}$,
defined by \er{ms6}, we have
\begin{multline}\label{ass5}
\P_{j,1,\be}(\vt,t)= {\sqrt{2}\/\sin z_1}\,\rt[\p_1(u)\big(\sin
z_{j,1}\cos(z_{j,1}\,t)-\cos z_{j,1}\sin(z_{j,1}\,t)\big)\\+
\p_1(v)\,e^{i\lan\t(\be),\,\vt\ran}\,\sin(z_{j,1}\,t)\rt]=
(-1)^j\sqrt{2}\,\p_1(u)\cos(z_{j,1}\,t)-\sqrt{2}\sin(z_{j,1}\,t)A(\vt),
\end{multline}
where
\[
\lb{as2'} A(\vt)={1\/\sin z_1}\Big(\p_1(u)\cos z_1-
\p_1(v)\,e^{i\lan\t(\be),\,\vt\ran}\Big).
\]

Using \er{egv1} and \er{lam}, we have the following asymptotics
\[\lb{ass3}
\begin{aligned}
& \cos z_1=-\l_1(\vt)=1-\ve^2\m(\omega)+O(\ve^3),
\\
& \sin z_1=\big(1-\l_1^2(\vt)\big)^{1/2}=
\big(2\ve^2\m(\omega)+O(\ve^3)\big)^{1/2}
=\ve\sqrt{2\m(\omega)}\;\big(1+O(\ve)\big).
\end{aligned}
\]
Substituting \er{lam} and \er{ass3} into \er{as2'}, we obtain
\begin{multline}\label{ass6}
A(\vt)={(1+O(\ve))\/\ve\sqrt{2\m(\omega)}\;}\Big[\p_1(0,u)+\ve\p_1^{(1)}(\omega,u)
-\p_1(0,v)(1+\ve\,i\lan\t(\be),\omega\ran)\\-\ve\p_1^{(1)}(\omega,v)+O(\ve^2)\Big]
={1\/\sqrt{2\m(\omega)}}\Big(\p_1^{(1)}(\omega,u)-\p_1^{(1)}(\omega,v)
-{i\/\sqrt{\vk}}\lan\t (\be),\omega\ran\Big)+O(\ve).
\end{multline}
Substituting \er{ass6} into \er{ass5} and using \er{lam} and
$$
\sin(z_{j,1}\,t)=\sin(\pi j_*t)+O(\ve), \qqq \cos(z_{j,1}\,t)=\cos(\pi j_*t)+O(\ve),
$$
we get \er{ass8}.

\emph{ii}) The identity \er{ass5} gives
\[\lb{nac}
|\P_{j,1,\be}(\vt,t)|<\sqrt{2}\big(1+|A(\vt)|\big).
\]
We have the decomposition
$$
\p_1(\vt,v)=\p_1(0,v)+\wt\p_1(\vt,v)=\vk^{-1/2}+\wt\p_1(\vt,v).
$$
We rewrite the difference on the right-hand side of \er{as2'} in the form
\begin{multline*}
\p_1(u)\cos z_1-\p_1(v)\,e^{i\lan\t(\be),\,\vt\ran}=\p_1(u)-\p_1(v)+\p_1(u)(\cos z_1-1)-\p_1(v)(e^{i\lan\t(\be),\,\vt\ran}-1)
\\=\wt\p_1(\vt,u)-\wt\p_1(\vt,v)-2\p_1(u)\sin^2{z_1\/2}- \p_1(v)(e^{i\lan\t(\be),\,\vt\ran}-1),
\end{multline*}
which yields
\[\lb{eeB}
\big|\p_1(u)\cos z_1-\p_1(v)\,e^{i\lan\t(\be),\,\vt\ran}\big|<|\wt\p_1(\vt,u)|+|\wt \p_1(\vt,v)|+2\sin^2{z_1\/2}+ \big|\lan\t(\be),\,\vt\ran\big|.
\]
Now we estimate $\wt\p_1(\vt,\cdot)$. Let
$$
\D(\vt)=\D(0)+\wt\D(\vt),\qqq \l_1(\vt)= \l_1(0)+\wt\l_1(\vt).
$$
Then we can rewrite the equation
$\D(\vt)\p_1(\vt,\cdot)=\l_1(\vt)\p_1(\vt,\cdot)$ in the form
\[\lb{gaas}
\big(\D(0)+\wt\D(\vt)\big)
\big(\p_1(0,\cdot)+\wt\p_1(\vt,\cdot)\big)=
\big(\l_1(0)+\wt\l_1(\vt)\big)
\big(\p_1(0,\cdot)+\wt\p_1(\vt,\cdot)\big)
\]
or, since $\D(0)\p_1(0,\cdot)=\l_1(0)\p_1(0,\cdot)$,
\[\lb{gaas'}
\big(\D(0)-\l_1(0)\1\big)\wt\p_1(\vt,\cdot)=
-\big(\wt\D(\vt)-\wt\l_1(\vt)\1\big)\p_1(\vt,\cdot),
\]
where $\1$ is the identity operator. Let $P$ be the orthogonal projection onto the subspace of $\ell^2(\cV_*)$ orthogonal to $\p_1(0,\cdot)$. Then from \er{gaas'} we obtain
\[\lb{p_1}
\wt\p_1(\vt,\cdot)=-D^{-1}
P\big(\wt\D(\vt)-\wt\l_1(\vt)\1\big)\p_1(\vt,\cdot),\qq \textrm{where} \qq D=P\big(\D(0)-\l_1(0)\1\big).
\]
This yields
\[\lb{sec1}
\begin{aligned}
\|\wt\p_1(\vt,\cdot)\|\leq\big\|D^{-1}P\big\|\cdot
\big\|\big(\wt\D(\vt)-\wt\l_1(\vt)\1\big)\p_1(\vt,\cdot)\big\|
\leq{1\/\L}\,\big\|\wt\D(\vt)-\wt\l_1(\vt)\1\big\|,
\end{aligned}
\]
where $\L$ is the distance between $\l_1(0)=-1$ and $\s\big(\D(0)\big)\sm\big\{\l_1(0)\big\}$. Due to \er{l2.15}, we have
\[
\lb{NN1}
\big\|\wt\D(\vt)-\wt\l_1(\vt)\1\big\|\leq\wt\l_1(\vt)+\max_{u\in
\cV_*}{1\/\vk_u}\sum_{\be=(u,v)\in\cA_*}|\lan\t
(\be),\,\vt\ran|\leq\wt\l_1(\vt)+|\vt|\,T,
\]
where $T=\max\limits_{u\in\cV_*}{1\/\vk_u}\sum\limits_{\be=(u,v)\in\cA_*}|\t(\be)|$.
Substituting \er{NN1} into \er{sec1}, we obtain
\[\lb{np1}
\|\wt\p_1(\vt,\cdot)\|\leq {1\/\L}\big(\,\wt\l_1(\vt)+ |\vt|\,T\big).
\]
 Using \er{eeB}, \er{np1}, for $A(\vt)$ defined by \er{as2'},
we have
\begin{multline}\lb{eB1}
|A(\vt)|\leq{1\/\sin
z_1}\bigg({2\/\L}\big(\wt\l_1(\vt)+|\vt|\,T\big)+2\sin^2{z_1\/2}+
\big|\lan\t (\be),\vt\ran\big|
\bigg)\\\leq\tan{z_1\/2}+{2\/\L}\cdot{\wt\l_1(\vt)\/\sin
z_1}+{|\vt|\/\sin
z_1}\bigg({2\/\L}\,T+\max_{\be\in\cA_*}\big|\t(\be)\big|\bigg).
\end{multline}
Let $z_1\in(0,{\pi\/2}]$. Then $\sin z_1\geq
{2\/\pi}\,z_1$. This and \er{eB1} give
\[\lb{eB2}
|A(\vt)|\leq1+{\pi\/\L}\cdot{\wt\l_1(\vt)\/z_1(\vt)}+
M\,{|\vt|\/\sin z_1(\vt)}\,,
\]
where $M$ is defined by \er{MMM}. We have the simple inequality
$$
\wt\l_1(\vt)=\l_1(\vt)-\l_1(0)=-\cos z_1(\vt)+1=
\int\limits_0^{z_1(\vt)}\sin t\,dt \leq
z_1(\vt),
$$
which yields that $\displaystyle{\wt\l_1(\vt)\/z_1(\vt)}\leq1$. Then the estimate \er{eB2} has the form
\[\lb{eB3}
|A(\vt)|\leq1+{\pi\/\L}+M\,{|\vt|\/\sin z_1(\vt)}\,.
\]
Combining \er{nac} and \er{eB3}, we obtain \er{est}.

\emph{iii}) We have ${\pi\/2}<z_1(\vt)\leq z_1^+<\pi$. This yields $0<\sin z_1^+\leq\sin z_1(\vt)$. We note that $z_1^+\neq\pi$, since we assume that for a bipartite periodic graph $\cG$ its fundamental graph $\cG_*$ is also bipartite and, consequently, has more than 1 vertex. Then the estimate \er{est'} follows directly from the identity \er{ms6}.

\emph{iv}) Using the second formula in \er{ass3}, we obtain
\[\lb{fra1}
{|\vt|\/\sin z_1(\vt)}={1\/\sqrt{2\m(\omega)}}+O\big(|\vt|\big) \qqq \textrm{as}
\qq |\vt|\to0.
\]
The following estimate holds true (see Theorem 1.2 in \cite{KS16b}):
\[\lb{fra2}
\begin{aligned}
\mu(\omega)\geq T_1={1\/\vk\,\n\,d}\left({d-1\/C}\right)^{d-1}>0, \qqq \forall\,\omega\in\S^{d-1},\\
\vk=\sum_{v\in\cV_*}\vk_v, \qqq C=\sum\limits_{j=1}^d|\t (\be_j)|^2,
\end{aligned}
\]
where $\be_1,\ldots,\be_d\in\cA_*$ are fundamental graph edges with linearly independent indices \linebreak $\t(\be_1),\ldots,\t(\be_d)\in\Z^d$.
Substituting the estimate \er{fra2} into \er{fra1}, we get
$$
{|\vt|\/\sin z_1(\vt)}\leq{1\/\sqrt{2T_1}}+O\big(|\vt|\big).
$$
Then there exists a constant $C_1>0$ such that
\[\lb{fra4}
{|\vt|\/\sin z_1(\vt)}<{1\/\sqrt{2T_1}}+C_1,\qq \textrm{for all} \qq |\vt|<R,
\]
where $R$ is some positive number. Now let $|\vt|\geq R$. Then $\sin z_1(\vt)\geq\a$ for some constant $\a>0$ and
\[\lb{fra5}
{|\vt|\/\sin z_1(\vt)}\leq{\pi\sqrt{d}\/\a}<\iy.
\]
Combining \er{fra4} and \er{fra5}, we obtain \er{sup}.
\qq \BBox

\medskip

\begin{proposition}\lb{T.As1}
i) All eigenfunctions $\P_{j,n}(\vt)$, $(j,\vt)\in\Z_0\ts(\T^d\sm\{0\})$, $n=2,\ldots,\n-1$, of $\D_M(\vt)$ satisfy
\[
\lb{ePnj} \|\P_{j,n}(\vt)\|_{L^\iy(\cG_*)}<{2\sqrt{2}\/\a}\,,
\qqq \a=\inf_{n=2,\ldots,\n-1 \atop \vt\in\T^d} \sin z_n(\vt)>0.
\]

ii) Let $\vt\in\T^d\sm\{0\}$ and let $\cG_*$ be non-bipartite. Then all eigenfunctions $\P_{j,\n}(\vt)$, $j\in\Z_0$, of $\D_M(\vt)$ satisfy
\[
\lb{efnb} \|\P_{j,\n}(\vt)\|_{L^\iy(\cG_*)}<{2\sqrt{2}\/\a_\n}\,,
\qqq \a_\n=\inf_{\vt\in\T^d} \sin z_\n(\vt)>0.
\]

iii) Let $\vt\in\T^d\sm\{0\}$ and let $\cG_*$ be bipartite. Then all eigenfunctions $\P_{j,\n}(\vt)$, $j\in\Z_0$, of $\D_M(\vt)$ satisfy
\[
\lb{estn}
\|\P_{j,\n}(\vt)\|_{L^\iy(\cG_*)}< \sqrt{2}\,\bigg(2+{\pi\/\L}+
M\,{|\vt|\/\sin z_1(\vt)}\bigg),
\]
where $M$ and $\L$ are defined in Proposition \ref{T.As}.ii.
\end{proposition}

\no \textbf{Proof.} \emph{i}) This estimate follows from the identity \er{ms6}.

\emph{ii}) Since $\cG$ is non-bipartite, using \er{ms6}, we also have \er{efnb}.

\emph{iii}) Let for the shortness
$$
z_{j,\n}=z_{j,\n}(\vt),\qqq z_\n=z_\n(\vt),\qqq \p_\n(u)=\p_\n(\vt,u).
$$
For the eigenfunction
$\P_{j,\n}(\vt)=\big(\P_{j,\n,\be}(\vt,t)\big)_{\be\in\cE_*}$,
defined by \er{ms6}, we have
\[\lb{ass5b}
\P_{j,\n,\be}(\vt,t)=
(-1)^j\sqrt{2}\,\p_\n(u)\cos(z_{j,\n}t)-\sqrt{2}\sin(z_{j,\n}t)A_\n(\vt),
\]
where
$$
\lb{as2'b} A_\n(\vt)={1\/\sin z_\n}\Big(\p_\n(u)\cos z_\n-
\p_\n(v)\,e^{i\lan\t(\be),\,\vt\ran}\Big).
$$
For the bipartite graph $\cG_*$ with parts $\cV_1$ and $\cV_2$ we have
$$
\sin z_\n=\sin z_1,\qq \cos z_\n=-\cos z_1,\qqq
\p_\n(v)=\ca \, \p_1(v), & \textrm{ if } v\in \cV_1\\
  -\p_1(v), & \textrm{ if } v\in\cV_2
\ac.
$$
This yields that $A_\n(\vt)=A(\vt)$, where $A(\vt)$ is defined by \er{as2'}.
Combining this, \er{ass5b} and \er{eB3}, we obtain \er{estn}. \qq \BBox

\begin{remark}
The estimates \er{est}, \er{est'}, \er{ePnj} -- \er{estn} give that the eigenfunctions $\P_{j,n}(\vt)$, $(j,n,\vt)\in\Z_0\ts\N_\n\ts\T^d\sm\{0\}$, are uniformly bounded with respect to $j,n,\vt$.
\end{remark}

\no \textbf{Proof of Corollary \ref{TsDM}.} \emph{i}) Recall that
$\D_M=\D_{MD}\os \D_{MV}$ and Theorem \ref{TT1} gives
\[
\s(\D_{MD})=\s_{fb}(\D_{MD})= \s_D.
\]
Thus, we need to consider the operator $\D_{MV}$. Since each $\l_n(\cdot)$
is a real and piecewise analytic function on the torus $\T^d$, each
$z_{j,n}(\cdot)$, $(j,n)\in\Z_0\ts\N_\n$, defined by the formula
\er{egv1}, is also a real and piecewise analytic function on $\T^d$.
Then \er{egv1} yields
\[
\s(\D_{MV})= \bigcup_{(j,n)\in\Z_0\ts\N_\n}\s_{j,n}(\D_M),\qqq
\s_{j,n}(\D_M)=[E_{j,n}^-,E_{j,n}^+]=z^2_{j,n}(\T^d).
\]
Since $\l_{n}(\cdot)=\const$, for $\n-r<n\le \n$ and all other eigenvalues $\l_1(\cdot),\dots,\l_{\n-r}(\cdot)$ are not constants (in any small ball), we obtain \er{sDM1}, \er{sDM1a}.

\emph{ii}) -- \emph{iv}) These items are the direct consequences of the previous one. \qq \BBox

\medskip

Recall that we identify an edge $\be_*=(u_*,v_*)\in\cE_\ast$ of the fundamental graph $\cG_*=(\cV_*,\cE_*)$ with an index $\t(\be_*)$ with the edge $\be=(u,v+\t(\be_*))\in\cE$ of the periodic graph $\cG=(\cV,\cE)$, where $u,v\in\cV_0=\cV\cap\Omega$ such that $u_*=\gf_\cV(u)$, $v_*=\gf_\cV(v)$.

\medskip

\no {\bf Proof of Proposition \ref{TLEF}.} \emph{i}) We fix $\vt\in\T^d\sm\{0\}$ and consider the space $X$ of all functions $f:\cV\rightarrow\C$ satisfying the following quasiperiodic condition:
\[\lb{qpco}
f(v+a(m))=e^{i\lan m,\,\vt\ran}f(v), \qqq \forall\,(v,m)\in\cV\ts\Z^d,
\]
where $a(m)$ is defined by \er{aofm}.
Such a function $f$ is uniquely determined by its values at the fundamental graph vertices. Thus, $X$ is isomorphic to $\ell^2(\cV_*)$. The fiber Laplacian $\D(\vt):\ell^2(\cV_*)\rightarrow\ell^2(\cV_*)$ can be also defined as the restriction of the operator $\D$ to the space~$X$.

Let $\p_n(\vt,\cdot)\in\ell^2(\cV_*)$ be the normalized
eigenfunction of $\D(\vt)$ corresponding to the eigenvalue $\l_n(\vt)$.
The function $\wt\p_n(\vt,\cdot)$ defined by \er{gde} satisfies the quasiperiodic condition:
\[\lb{gde1}
\wt\p_n(\vt,v+a(m))=e^{i\lan m,\,\vt\ran}\wt\p_n(\vt,v), \qqq \forall\,(v,m)\in \cV\ts\Z^d,
\]
i.e., $\wt\p_n(\vt,\cdot)\in X$.
For each $v\in \cV$ we have the unique representation $v=v_0+a(m)$ for some $v_0\in \cV_0$ and $m\in\Z^d$ and, due to \er{gde1}, \er{l2.15} and the definition of edge indices, we obtain
\begin{multline*}
(\D\wt\p_n(\vt,\cdot))(v)=(\D(\vt)\wt\p_n(\vt,\cdot))(v)=
e^{i\lan m,\,\vt\ran}(\D(\vt)\wt\p_n(\vt,\cdot))(v_0)\\=-\,\frac{e^{i\lan m,\,\vt\ran}}{\vk_{v_0}}\sum\limits_{(v_0,u)\in\cA}\wt\p_n(\vt,u)=-\,\frac{e^{i\lan m,\,\vt\ran}}{\vk_{v_0}}\sum\limits_{\be=(v_0,u_0)\in\cA_*}e^{i\lan \t(\be),\,\vt\ran}\wt\p_n(\vt,u_0)\\=-\,\frac{e^{i\lan m,\,\vt\ran}}{\vk_{v_0}}\sum\limits_{\be=(v_0,u_0)\in\cA_*}e^{i\lan \t(\be),\,\vt\ran}\p_n(\vt,u_0)=e^{i\lan m,\,\vt\ran}(\D(\vt)\p_n(\vt,\cdot))(v_0)\\=e^{i\lan m,\,\vt\ran}\l_n(\vt)\p_n(\vt,v_0)=\l_n(\vt)\wt\p_n(\vt,v).
\end{multline*}
Thus, $\wt\p_n(\vt,\cdot)\in\ell^2_{loc}(\cV)$ is a generalized
eigenfunction of $\D$ corresponding to the eigenvalue $\l_n(\vt)$.

\emph{ii}) The proof of this item is similar to the proof of the item \emph{i}).

\emph{iii}) If we denote $u=v+a(m)$ in \er{gde1}, where $u\in\cV\ss\R^d$, then for all $v\in \cV$ \er{gde} gives
$$
\wt\p_n(\vt,u)=e^{i\lan u_\A-v_\A,\,\vt\ran}\wt\p_n(\vt,u-m(a))=
e^{i\lan u_\A,\,\vt\ran}e^{-i\lan v_\A+m,\,\vt\ran}\wt\p_n(\vt,u)=e^{i\lan u_\A,\,\vt\ran}\vp_n(\vt,u).
$$
The function $\vp_n(\vt,u)=e^{-i\lan u_\A,\,\vt\ran}\wt\p_n(\vt,u)$ is $\G$-periodic with respect to $u$, since
\begin{multline*}
\vp_n(\vt,u+a(m))=e^{-i\lan u_\A+m,\,\vt\ran}\wt\p_n(\vt,u+a(m))=e^{-i\lan u_\A+m,\,\vt\ran}e^{i\lan m,\,\vt\ran}\wt\p_n(\vt,u)\\=e^{-i\lan u_\A,\,\vt\ran}\wt\p_n(\vt,u)=\vp_n(\vt,u),\qqq (u,m)\in\cV\ts\Z^d.
\end{multline*}
Thus, \er{gdep} has been proved. The proof of the identity \er{gmep} is similar.
\qq \BBox

\medskip

We recall known results (see Theorems 4.5.2 and 4.5.4 in \cite{BK13})  about the connection between true eigenfunctions from $\ell^2(\cV)$ and $L^2(\cG)$
of the discrete and metric Laplacians, respectively, and eigenfunctions of their fiber operators.

\begin{proposition}
\label{TLBS}
Let $\p_n(\vt,\cdot)\in\ell^2(\cV_*)$, $\vt\in\T^d\sm\{0\}$, be the family of the normalized eigenfunctions of the fiber operators $\D(\vt)$ corresponding to the eigenvalue $\l_n\in\s_{fb}(\D)$ for some $n=\nu-r+1,\ldots,\nu$ and let $\P_{j,n}(\vt,\cdot)\in L^2(\cG_*)$ be the family of the normalized eigenfunctions of the fiber operators $\D_M(\vt)$ corresponding to the eigenvalue $z_{j,n}^2\in\s_{fb}(\D_M)$ defined by \er{egv1}, $j=0,1,2,\ldots$\,. Then the following statements hold true.

i) The eigenfunction $\wt\p_n(\cdot)\in\ell^2(\cV)$ of the discrete Laplacian $\D$ corresponding to the eigenvalue $\l_n$ has the form
\[\lb{bsd}
\wt\p_n(v+a(m))=\int_{\T^d}e^{i\lan m,\vt\ran }\p_n(\vt,v){d\vt\/(2\pi)^d}\,, \qqq \forall\,(v,m)\in\cV_*\ts\Z^d,
\]
where $a(m)\in\G$ is defined by \er{aofm}.

ii) The eigenfunction $\wt\P_{j,n}(\cdot)\in L^2(\cG)$ of the metric Laplacian $\D_M$ corresponding to the eigenvalue $z_{j,n}^2$ has the form
\[\lb{bsm}
\wt\P_{j,n}(x+a(m))=\int_{\T^d}e^{i\lan m,\vt\ran }\P_{j,n}(\vt,x){d\vt\/(2\pi)^d}\,, \qqq \forall\,(x,m)\in\cG_*\ts\Z^d.
\]
\end{proposition}

\begin{remark}
1) The right-hand sides of \er{bsd} and \er{bsm} are the inverse operators of the unitary operators (Floquet transforms) $U$ and $\mU$ from \er{raz} and \er{5001}, respectively.

2) The eigenfunctions defined by \er{bsd} and \er{bsm} are compactly supported.
\end{remark}

\section {\lb{Sec4}  Schr\"odinger operators on periodic metric graphs}
\setcounter{equation}{0}

\subsection{Trace class estimates.}
We consider the Schr\"odinger operator $H$ on the periodic graph $\cG$ given by
$$
H=H_0+Q,\qqq H_0=\D_M.
$$
We assume that the potential  $Q=(Q_{\be})_{\be\in \cE}$ is real and
belongs to the space $L^{1}(\cG)\cap L^{2}(\cG)$.
Recall that  $R_0(k)=(H_0-k^2)^{-1}$, $R(k)=(H-k^2)^{-1}$ and
\[
\label{SY'}
\begin{aligned}
Y_0(k) = |Q|^{1/2}\, R_0(k)\, Q^{1/2},\qqq Y(k) = |Q|^{1/2}\,
R(k)\, Q^{1/2}, \quad  Q^{1/2} =|Q|^{1/2}\sign Q,
\end{aligned}
\]
$k\in\C_+$.  Note  that Propositions \ref{T.As} and \ref{T.As1} give
\[\lb{cgam}
C_\cG:=\sup_{(j,n,\vt)\in\Z_0\ts\N_\n\ts\T^d}\|\P_{j,n}(\vt)\|_{L^\iy(\cG_*)}<\iy.
\]

Let $\cR_k(x,x',\vt)$, $x,x'\in\cG_*$, be the kernel of the integral
operator $(\D_{MV}(\vt)-k^2)^{-1}$ on $L^2(\cG_*)$. Due to Theorem
\ref{TT2} we have
\[
\lb{Rk} \cR_k(x,x',\vt)=\sum_{j=0}^\iy\sum_{n=1}^\n {\P_{j,n}(\vt,
x)\ol\P_{j,n}(\vt, x')\/z_{j,n}^2(\vt)-k^2}\,,\qqq x,x'\in\cG_*.
\]

\begin{theorem}
\label{TRQ} i) Let $(Q,k)\in L^{1}(\cG)\ts \C_+$ and let
$\D_M^{1/2}\ge 0$. Then $\big(\D_M^{1/2}-k\big)^{-1}|Q|^{1/2}\in \bB_2$
and $Y_0(k)\in \bB_1$ and they satisfy
\[
\lb{S8}
\begin{aligned}
 \big\|(\D_M^{1/2}-k)^{-1}|Q|^{1/2}\big\|_{\bB_2}^2\le C_0\|Q\|_{L^1(\cG)}\rt(
{2\/(\Im k)^2}+ {1\/\Im k} \rt),
\end{aligned}
\]
\[
\lb{S9} \bea
&\|Y_0(k)\|_{ \bB_1}\le \cC(Q,k),\\
&\textrm{where } \ \cC(Q,k)={ C_0} \|Q\|_{L^1(\cG)} \rt( {2\/(\Im
k)^2}+ {1\/\Im k} \rt)^{1/2}\rt( {2\/|k|^2}+ {1\/|k|}
\rt)^{1/2},\\
& C_0=2(\n C_\cG^2+\n_1-\n), \qqq \n=\#\cV_*, \qqq \n_1=\#\cE_*,
 \ena
\]
\[
\lb{S9x} \|Y_0(k)\|_{ \bB_1}={O(1)\/|k\Im k|^{1\/2}}\qqq as \qqq \Im
k\to \iy.
\]
Moreover, if  $\Re k\le 0$, then
\[
\lb{se2} \big\|(\D_M^{1/2}-k)^{-1}|Q|^{1/2}\big\|_{\bB_2}^2\le
C_0\|Q\|_{L^1(\cG)}\rt({2\/|k|^2}+ {1\/|k|}\rt).
\]
ii) Let, in addition,
\[\lb{leq}
\sup_{x\in\cG_*}\big|Q(x+a)\big|=\g_a, \qq a\in\G, \qqq \sum_{a\in\G}\g_a=\g<\iy.
\]
Then
\[
\lb{S8'}
\begin{aligned}
\big\|(\D_M^{1/2}-k)^{-1}|Q|^{1/2}\big\|_{\bB_2}^2
\leq\g(\n_1+\n)\bigg({2\/(\Im k)^2}+{1\/\Im k}\bigg).
\end{aligned}
\]
Moreover, if $\Re k\le 0$, then
\[
\lb{se2'} \big\|(\D_M^{1/2}-k)^{-1}|Q|^{1/2}\big\|_{\bB_2}^2\le
\g(\n_1+\n)\rt({2\/|k|^2}+ {1\/|k|}
\rt).
\]
\end{theorem}
\no {\bf Proof.} \emph{i}) Due to Theorem \ref{T1}, we have
\[
\lb{S10}
\Big(\mU(\D_M^{1/2}-k)^{-1}h\Big)(\vt,\cdot)=
\Big(\big(\D_{MD}^{1/2}(\vt)-k\big)^{-1}\os
\big(\D_{MV}^{1/2}(\vt)-k\big)^{-1}\Big)(\mU h)(\vt,\cdot),
\]
where the unitary operator $\mU:L^2(\cG)\to\mH$ is defined by
\er{5001}, $h\in L^2(\cG)$ is compactly supported and $\vt\in \T^d$.
We consider the second term on the right-hand side of \er{S10}. Using Theorem \ref{TT2} and denoting $\a=(j,n)\in\Z_0\ts\N_\n$, $\Z_0=\{0,1,2,\ldots\}$, we obtain
\begin{multline*}
\big(\D_{MV}^{1/2}(\vt)-k\big)^{-1}(\mU h)(\vt,x)=\sum_{\a\in\Z_0\ts\N_\n}{1\/z_\a(\vt)-k}\,\cP_\a(\vt)
(\mU h)(\vt,x)\\
=\sum_{\a\in\Z_0\ts\N_\n}{1\/z_\a(\vt)-k}\,\P_\a(\vt,x)\int_{\cG_*}
(\mU h)(\vt,x')\ol{\P}_\a(\vt,x')\,dx'.
\end{multline*}
Using $Q\in L^{1}(\cG)$, the normalized condition $\int_{\cG_*}
\big|\P_\a(\vt,x)\big|^2\,dx=1$ and \er{cgam}, this gives
\begin{multline}\label{scc3}
\big\|\big(\D_{MV}^{1/2}(\cdot)-k\big)^{-1}\mU |Q|^{1/2}\big\|_{\bB_2}^2\\=
{1\/(2\pi)^d}\sum\limits_{a\in\G}\sum_{\a\in\Z_0\ts\N_\n}\int_{\T^d}
{1\/|z_\a(\vt)-k|^2}\int_{\cG_*}\big|Q(x+a)\big|\cdot
\big|\P_\a(\vt,x)\big|^2\,dxd\vt\\
\leq{C_\cG^2\/(2\pi)^d}\sum\limits_{a\in\G}\sum_{\a\in\Z_0\ts\N_\n}
\int_{\T^d}{d\vt\/|z_\a(\vt)-k|^2}
\int_{\cG_*}\big|Q(x+a)\big|\,dx\\
={C_\cG^2\/(2\pi)^d}\,\|Q\|_{L^1(\cG)}\sum_{\a\in\Z_0\ts\N_\n}
\int_{\T^d}{d\vt\/|z_\a(\vt)-k|^2}\,.
\end{multline}
Let $k=p+iq\in\C$, $q\neq0$.  We assume that the following inequality holds true
\[\lb{scc4}
\sum_{j=0}^\iy{1\/|z_{j,n}(\vt)-k|^2}
\leq{4\/q^2}+{2\/|q|}\, \qq \textrm{ for each } \; n\in\N_\n.
\]
Substituting \er{scc4} into \er{scc3}, we obtain
\[\lb{scc5}
\begin{aligned}
\big\|\big(\D_{MV}^{1/2}(\cdot)-k\big)^{-1}\mU |Q|^{1/2}\big\|_{\bB_2}^2
\leq
2\nu\, C_\cG^2\,\|Q\|_{L^1(\cG)}\bigg({2\/q^2}+{1\/|q|}\bigg).
\end{aligned}
\]
Now we prove \er{scc4}. Due to \er{egv1}, we have the following identity
$$
\sum_{j=0}^\iy{1\/|z_{j,n}(\vt)-k|^2}=\sum_{j=0 \atop j\textrm{ is even}}^\iy
{1\/(z_n(\vt)+\pi j-p\big)^2+q^2}+\sum_{j=1 \atop j\textrm{ is odd}}^\iy
{1\/(\pi(j+1)-z_n(\vt)-p\big)^2+q^2}\,.
$$
Let
$$
p-z_n(\vt)=\pi(m+b), \qqq p+z_n(\vt)=\pi(m_1+b_1)
$$
for some $m,m_1\in \Z$, and $b=b(\vt),b_1=b_1(\vt)\in (0,1)$. Then
\begin{multline*}
\sum_{j=0}^\iy{1\/|z_{j,n}(\vt)-k|^2}=
\sum_{j=0 \atop j\textrm{ is even}}^\iy
{1\/\pi^2(j-m-b)^2+q^2}+\sum_{j=1 \atop j\textrm{ is odd}}^\iy
{1\/\pi^2(j+1-m_1-b_1)^2+q^2}\\
\leq\sum_{j=-\iy}^\iy
{1\/\pi^2(j-b)^2+q^2}+\sum_{j=-\iy}^\iy
{1\/\pi^2(j-b_1)^2+q^2}\leq{4\/q^2}+\sum_{j=1}^\iy{4\/\pi^2 j^2+q^2}\\
\leq{4\/q^2}+4\int_0^\iy {dt\/\pi^2t^2+q^2}={4\/q^2}+{2\/|q|}\,.
\end{multline*}

Now we consider the first term on the right-hand side of \er{S10}. Similar arguments yield
\[\lb{Dscc5}
\begin{aligned}
\big\|\big(\D_{MD}^{1/2}(\cdot)-k\big)^{-1}\mU |Q|^{1/2}\big\|_{\bB_2}^2
\leq2(\n_1-\n)\,\|Q\|_{L^1(\cG)}\bigg({2\/q^2}+{1\/|q|}\bigg).
\end{aligned}
\]
Summing \er{scc5} and \er{Dscc5} and using \er{S10}, we get
\begin{multline}\label{Dscc6}
\big\|(\D_M^{1/2}-k)^{-1}|Q|^{1/2}\big\|_{\bB_2}^2=
\big\|\mU(\D_M^{1/2}-k)^{-1}|Q|^{1/2}\big\|_{\bB_2}^2\\
\leq2(\nu C_\cG^2+\n_1-\n) \|Q\|_{L^1(\cG)}\bigg({2\/q^2}+{1\/|q|}\bigg)=C_0\|Q\|_{L^1(\cG)}\rt(
{2\/(\Im k)^2}+ {1\/|\Im k|}\rt)\,.
\end{multline}
Thus, \er{S8} has been proved.

Let, in addition,  $\Re k\le 0$. Then
\[\lb{rek0}
|z_\a(\vt)-k|^2=\big(z_\a(\vt)-p\big)^2+q^2\ge z_\a^2(\vt)+p^2+q^2=
z_\a^2(\vt)+|k|^2
\]
and using the above arguments we obtain \er{se2}.

For any $k\in\C_+$ we have the following factorization
$$
Y_0(k)=|Q|^{1/2}R_0(k)Q^{1/2}=\big(|Q|^{1/2}(\D_M^{1/2}+k)^{-1}\big)
\big((\D_M^{1/2}-k)^{-1}Q^{1/2}\big).
$$
Then
$$
\|Y_0(k)\|_{\bB_1}\leq\big\||Q|^{1/2}(\D_M^{1/2}+k)^{-1}\big\|_{\bB_2}
\big\|(\D_M^{1/2}-k)^{-1}Q^{1/2}\big\|_{\bB_2}
$$
and, applying the estimates \er{Dscc6} and \er{se2}, we obtain
\er{S9}, which yields  asymptotics
\er{S9x}.

\emph{ii}) Let the potential $Q$ satisfy \er{leq}. Then using \er{scc3} and \er{scc4}, we obtain
\begin{multline}\label{scc2'}
\big\|\big(\D_{MV}^{1/2}(\cdot)-k\big)^{-1}\mU |Q|^{1/2}\big\|_{\bB_2}^2\\
={1\/(2\pi)^d}\sum\limits_{a\in\G}\sum_{\a\in\Z_0\ts\N_\n}\int_{\T^d}
{1\/|z_\a(\vt)-k|^2}
\int_{\cG_*}\big|Q(x+a)\big|\cdot\big|{\P}_\a(\vt,x)\big|^2\,dx\,d\vt\\
\leq{1\/(2\pi)^d}\sum\limits_{a\in\G}\g_a\sum_{\a\in\Z_0\ts\N_\n}
\int_{\T^d}{d\vt\/|z_\a(\vt)-k|^2}\leq2\g\,\n\bigg({2\/(\Im k)^2}+{1\/|\Im k|}\bigg).
\end{multline}
Similar arguments yield
\[\lb{Dscc2'}
\big\|\big(\D_{MD}^{1/2}(\cdot)-k\big)^{-1}\mU |Q|^{1/2}\big\|_{\bB_2}^2\leq \g(\n_1-\n)\bigg({2\/(\Im k)^2}+{1\/|\Im k|}\bigg).
\]
Summing \er{scc2'} and \er{Dscc2'} and using \er{S10}, we get \er{S8'}.

Let, in addition,  $\Re k\le 0$. Then using \er{rek0}
 and the above arguments we obtain \er{se2'}.
\qq \BBox

\begin{lemma}
\label{TRQ2} Let $(Q,k)\in L^{1}(\cG)\ts \C_+$. Then
\[
\lb{TY0} \Tr Y_0(k)={1\/(2\pi)^d}\int_{\T^d}\rt(\sum_{j=0}^\iy\sum_{n=1}^\n
{Q_{j,n}(\vt)\/z_{j,n}^2(\vt)-k^2}+\sum_{j=1}^\iy\sum_{s=1}^{\b-1}
{Q_{j,s}^0(\vt)\/(\pi j)^2-k^2}\rt)d\vt,
\]
where $\b=\#\cE_*-\#\cV_*+1$ is the Betti number of the fundamental graph $\cG_*=(\cV_*,\cE_*)$, $\n=\#\cV_*$, and
\[
\begin{aligned}
Q_{j,n}(\vt)=\int_{\cG_*}
\sum\limits_{a\in\G}Q(x+a)\big|\P_{j,n}(\vt,x)\big|^2dx,\\
Q_{j,s}^0(\vt)=\int_{\cG_*}
\sum\limits_{a\in\G}Q(x+a)\big|\P_{j,s}^0(\vt,x)\big|^2dx.
\end{aligned}
\]
\end{lemma}
\no {\bf Proof.}  Due to Theorem \ref{TT2} and \er{Rk} we have
\begin{multline*}
\big((\D_{MV}-k^2)^{-1}f,g\big)\\
=\int_{\cG_*}dx\int_{\cG_*}dx'\int_{\T^d}{d\vt\/(2\pi)^d}
\sum\limits_{m,\,m'\in\Z^d}e^{-i\lan m-m',\vt\ran}\cR_k(x,x',\vt)f(x+a(m))\ol g(x'+a(m')),
\end{multline*}
where $f,g\in L^2(\cG)$ are compactly supported functions, and $a(m)\in\G$ is defined by \er{aofm}. This yields
\begin{multline*}
\Tr\big(|Q|^{1/2}(\D_{MV}-k^2)^{-1}Q^{1/2}\big)=\int_{\cG_*}dx
\int_{\T^d}{d\vt\/(2\pi)^d}\sum\limits_{a\in\G}\cR_k(x,x,\vt)Q(x+a)\\
=\sum_{j=0}^\iy\sum_{n=1}^\n\int_{\cG_*}dx\int_{\T^d}{d\vt\/(2\pi)^d}
\sum\limits_{a\in\G}{Q(x+a)\big|\P_{j,n}(\vt,
x)\big|^2\/z_{j,n}^2(\vt)-k^2} ={1\/(2\pi)^d}\int_{\T^d}\sum_{j=0}^\iy\sum_{n=1}^\n
{Q_{j,n}(\vt)\/z_{j,n}^2(\vt)-k^2}\,d\vt.
\end{multline*}
Similar arguments yield
\begin{multline*}
\Tr\big(|Q|^{1/2}(\D_{MD}-k^2)^{-1} Q^{1/2}\big) \\
=\sum_{j=1}^\iy\sum_{s=1}^{\b-1}\int_{\cG_*}dx
\int_{\T^d}{d\vt\/(2\pi)^d}\sum\limits_{a\in\G}{Q(x+a)\big|\P_{j,s}^0(\vt,
x)\big|^2\/(\pi j)^2-k^2}={1\/(2\pi)^d}\int_{\T^d}\sum_{j=1}^\iy\sum_{s=1}^{\b-1}
{Q_{j,s}^0(\vt)\/(\pi j)^2-k^2}\,d\vt.
\end{multline*}
Substituting these formulas into the identity
$$
\Tr Y_0(k)=\Tr\big(|Q|^{1/2}(\D_{MV}-k^2)^{-1} Q^{1/2}\big)+\Tr
\big(|Q|^{1/2}(\D_{MD}-k^2)^{-1} Q^{1/2}\big),
$$
we obtain \er{TY0}. \qq \BBox

\subsection{Fredholm determinants} In order to discuss Fredholm
determinants we recall well-known facts~\cite{S05}. Let $\mathcal H$ be a
Hilbert space endowed with an inner product $(\, , \, )$ and a norm
$\|\cdot\|$. Let  $\cB$ denote the class of bounded operators.
Let $\bB_1$ be the set of all trace class operators on
$\mathcal H$ equipped with the trace norm $\|\cdot\|_{\bB_1}$.

$\bu$ Let $A, B\in \cB$ and $AB, BA\in \bB_1$. Then
\begin{equation}
\label{AB} {\rm Tr}\, AB={\rm Tr}\, BA,
\end{equation}
\begin{equation}
\label{1+AB} \det (I+ AB)=\det (I+BA).
\end{equation}
$\bu$ Let $A, B\in \bB_1$. Then
\begin{equation}
\label{DA1} |\det (I+ A)|\le e^{\|A\|_{\bB_1}}.
\end{equation}
\begin{equation}
\label{DA1x} |\det (I+ A)-\det (I+ B)|\le \|A-B\|_{\bB_1}
e^{1+\|A\|_{\bB_1}+\|B\|_{\bB_1}}.
\end{equation}
Moreover,  $I+ A$ is invertible if and only if $\det (I+ A)\ne 0$.

$\bu$  Suppose for a domain $\mD \subset {\C}$, a function
$\O(\cdot)-I: \mD\to \bB_1 $ is analytic and is invertible
for any $z\in\mD$. Then $F(z)=\det\O(z)$ is analytic and satisfies
\begin{equation}
\label{F'z}
 F'(z)= F(z){\rm Tr}\,\O^{-1}(z)\O'(z).
\end{equation}

By Theorem \ref{TRQ}, each operator $Y_0(k)\in \bB_1$, $k\in \C_+$, and
we can define the determinant
\[
\begin{aligned}
D(k)=\det (I+Y_0(k)),\qqq  k\in \C_+.
\end{aligned}
\]

\begin{lemma}
\label{TRQ3}
Let $Q\in L^{1}(\cG)$. Then the determinant $D(k)$, $k\in\C_+$, is analytic in $k\in\C_+$. If, in addition, $\cC(Q,k)<1$, where the constant $\cC(Q,k)$ is
defined by \er{S9}, then
\begin{equation}
\lb{Dx1}
\log D(k)=-\sum_{n=1}^{\iy}{1\/n}\Tr\big(-Y_0(k)\big)^n,
\end{equation}
where the series converges absolutely and uniformly, and for any
$N\geq 0$ we have
\[
\lb{Dx2}
 \Big|\log D(k)+\sum _{n=1}^{N}{\Tr (-Y_0(k))^n\/n}\Big|\le
{\|Y_0(k)\|_{\bB_1}^{N+1}\/(N+1)(1-\cC(Q,k))}\,.
\]

\end{lemma}
\no {\bf Proof.} By Theorem \ref{TRQ}, the $\bB_1$-valued function $Y_0(k)$
is analytic in $\C_+$. The series \er{Dx1} is standard, see
\cite{S05}. Under the condition $\|Y_0(k)\|_{\bB_1}\leq\cC(Q,k)$ the identity \er{Dx1}
gives \er{Dx2}. \BBox

\medskip

\no {\bf Proof of Theorem \ref{TSD}.}  Due to \er{S8}, \er{S9} we have
\er{Sc1}. Then it is well known that there exist the wave operators
$W_\pm$ and they are complete, see \cite{RS79}.   Moreover, the
S-operator given by \er{Sc3} is unitary on $\mH_{ac}(H_0)$ and the corresponding
S-matrix $S(k)$ for almost all $ k^2 \in \s_{ac}(H_0)$ satisfies
\er{Sc4} and \er{Sc5}. \BBox

\section {\lb{Sec5} $d$-dimensional lattice}
\setcounter{equation}{0}

We consider the lattice graph $\dL^d=(\cV,\cE)$, where the vertex set and the edge set are given by
\[
\lb{dLg}
\cV=\Z^d,\qquad  \cE=\big\{(m,m+a_1),\ldots,(m,m+a_d), \quad
\forall\,m\in\Z^d\big\},
\]
and $a_1,\ldots,a_d$ is the standard orthonormal basis, see Fig.\ref{ff.0.1}\emph{a}. The "minimal"\, fundamental graph $\dL_*^d=(\cV_*,\cE_*)$ of the lattice $\dL^d$ consists of one vertex $v$ and $d$ loop edges $\be_1=\ldots=\be_d=(v,v)$ with indices $a_1,\ldots,a_d$, see Fig.\ref{ff.0.1}\emph{b}. For each $\vt=(\vt_1,\ldots,\vt_d)\in\T^d$ the fiber operator $\D(\vt)$ degenerates to the scalar function
\[
\textstyle \D(\vt)=-{1\/d}\,(\cos\vt_1+\ldots+\cos\vt_d),
\]
which yields that the unique eigenvalue of $\D(\vt)$ and the corresponding normalized eigenfunction are given by
\[\lb{sql2}
\l(\vt)=\D(\vt),\qqq \p(\vt,v)=(2d)^{-1/2}.
\]
The fiber operator $\D_M(\vt)$ acts on functions
$y=(y_\be)_{\be\in\cE_*}\in L^2(\dL_*^d)$ by
$(\D_M(\vt)y)_\be=-y''_\be$, where $(y''_\be)_{\be\in\cE_*}\in L^2(\dL_*^d)$ and
$y$ satisfies the quasi-periodic conditions:
\[\lb{sq1}
y_{\be_1}(0)=\ldots=y_{\be_d}(0)=e^{-i\vt_1}y_{\be_1}(1)=\ldots=
e^{-i\vt_d}y_{\be_d}(1),
\]
\[\lb{sq3}
y'_{\be_1}(0)+\ldots+y'_{\be_d}(0)-e^{-i\vt_1}y'_{\be_1}(1)-\ldots-
e^{-i\vt_d}y'_{\be_d}(1)=0.
\]
Indeed, substituting the indices $a_1,\ldots,a_d$ of the fundamental graph edges  into the formulas \er{FBC}, \er{di2}, we obtain the conditions \er{sq1} -- \er{sq3}.

\setlength{\unitlength}{1.0mm}
\begin{figure}[h]
\centering

\unitlength 1.0mm % = 2.845pt
\linethickness{0.4pt}
\ifx\plotpoint\undefined\newsavebox{\plotpoint}\fi % GNUPLOT compatibility

\begin{picture}(140,60)(0,0)

%  убическа€ решетка

\put(10,10){\line(1,0){40.00}}
\put(10,30){\line(1,0){40.00}}
\put(10,50){\line(1,0){40.00}}

\bezier{60}(17,33.5)(37,33.5)(57,33.5)
\bezier{60}(17,13.5)(37,13.5)(57,13.5)

\put(17,53.5){\line(1,0){40.00}}
\put(10,10){\line(0,1){40.00}}
\put(30,10){\line(0,1){40.00}}
\put(50,10){\line(0,1){40.00}}
\put(57,13.5){\line(0,1){40.00}}
\bezier{12}(10,10)(13.5,11.75)(17,13.5)

\bezier{60}(17,13.5)(17,33.5)(17,53.5)
\bezier{60}(37,13.5)(37,33.5)(37,53.5)
\bezier{12}(30,10)(33.5,11.75)(37,13.5)
\put(50,10){\line(2,1){7.00}}
\bezier{12}(10,30)(13.5,31.75)(17,33.5)
\bezier{12}(30,30)(33.5,31.75)(37,33.5)

\put(50,30){\line(2,1){7.00}}
\put(10,50){\line(2,1){7.00}}
\put(30,50){\line(2,1){7.00}}
\put(50,50){\line(2,1){7.00}}
\put(10,10){\circle{1.0}}
\put(30,10){\circle{1.0}}
\put(50,10){\circle*{1.0}}
\put(10,30){\circle{1.0}}
\put(30,30){\circle{1.0}}
\put(50,30){\circle{1.0}}
\put(10,50){\circle{1.0}}
\put(30,50){\circle{1.0}}
\put(50,50){\circle{1.0}}

\put(17,53.5){\circle{1.0}}
\put(37,53.5){\circle{1.0}}
\put(57,53.5){\circle{1.0}}
\put(17,33.5){\circle{1.0}}
\put(37,33.5){\circle{1.0}}
\put(57,33.5){\circle{1.0}}
\put(17,13.5){\circle{1.0}}
\put(37,13.5){\circle{1.0}}
\put(57,13.5){\circle{1.0}}

\put(50,10){\vector(0,1){20.00}}
\put(50,10){\vector(-1,0){20.00}}
\put(50,10){\vector(2,1){7.00}}

\put(49.0,7.0){$v$}
\put(54,10.0){$a_1$}
\put(35.0,7.5){$a_2$}
\put(46.0,26){$a_3$}

\put(-5,5){(\emph{a})}
%*********************************
\bezier{30}(87,13.5)(97,13.5)(107,13.5)
\bezier{30}(87,13.5)(87,23.5)(87,33.5)
\bezier{12}(87,13.5)(83.5,11.75)(80,10)

\put(100,10){\vector(2,1){7.00}}
\bezier{12}(100,30)(103.5,31.75)(107,33.5)
\bezier{12}(80,30)(83.5,31.75)(87,33.5)
\bezier{30}(80,10)(80,20)(80,30)
\bezier{30}(107,13.5)(107,23.5)(107,33.5)
\put(100,10){\vector(-1,0){20.00}}
\bezier{30}(87,33.5)(97,33.5)(107,33.5)

\put(100,10){\vector(0,1){20.00}}
\bezier{30}(80,30)(90,30)(100,30)
\put(80,10){\circle{1}}
\put(80,30){\circle{1}}
\put(100,30){\circle{1}}
\put(87,13.5){\circle{1}}
\put(107,13.5){\circle{1}}
\put(107,33.5){\circle{1}}
\put(87,33.5){\circle{1}}

\put(100,10){\circle*{1}}

\put(100.0,7.0){$v$}
\put(79.0,7.0){$v$}
\put(101.0,28.5){$v$}
\put(77.0,28.5){$v$}

\put(84.5,33.5){$v$}
\put(108.0,33.5){$v$}
\put(108.0,13.5){$v$}

\put(103.0,9.5){$\be_1$}
\put(89.0,7.0){$\be_2$}
\put(96.0,21.0){$\be_3$}

\put(65,5){(\emph{b})}
\end{picture}

\vspace{-0.5cm} \caption{\footnotesize \emph{a}) The lattice
$\dL^3$;\quad \emph{b})  the fundamental graph $\dL^3_*$.}
\label{ff.0.1}
\end{figure}

\begin{proposition}\lb{T.Gr}
Let $j\in\N$. Then for each $\vt\in\T^d\sm\{0\}$ the operator
$\D_M(\vt)$ has the eigenvalue $(\pi j)^2$ of
multiplicity $d-1$ and the corresponding normalized
eigenfunctions $\P_{j,s}^0=(\P_{j,s,\be}^0)_{\be\in\cE_*}, \
s=1,\ldots,d-1$,  have the form
\[
\lb{X1}
\begin{aligned}
& \P_{j,s,\be_p}^0(\vt,t)={\sqrt{2}\;J_p\sin(\pi
jt)\/\sqrt{1+|C_{j,s}(\vt)|^2}}\, ,\qqq \ca J_d=1, \qq J_s=-C_{j,s}(\vt)\\
J_p=0,\qq \forall\, p\neq s,d\ac,\qqq t\in[0,1],\\
&  \textrm{ where } \qq C_{j,s}(\vt)=
{1-(-1)^je^{-i\vt_d}\/1-(-1)^je^{-i\vt_s}}\,.
\end{aligned}
\]
\end{proposition}

\no \textbf{Proof.}  The Betti number $\b$ of the fundamental graph $\dL^d_*$, defined by \er{Benu}, is equal to $d$. Then, due to Theorem \ref{TT1}, for all
$\vt\in\T^d\sm\{0\}$ the eigenvalue $(\pi j)^2$ has multiplicity
$d-1$ and all normalized eigenfunctions corresponding to the
eigenvalue $(\pi j)^2$, have the form
\[\lb{co1}
\P_{j,s}^0=\big(\P_{j,s,\be}^0\big)_{\be\in\cE_*},\qqq
\P_{j,s,\be}^0(\vt,t)=X_{j,s,\be}(\vt)\sqrt{2}\, \sin(\pi jt)\,,\qq
s=1,\ldots,d-1,
\]
where $\big(X_{j,s,\be}(\vt)\big)_{\be\in\cE_*}$ is a normalized solution of the equation
$$
\sum\limits_{p=1}^dx_p\big(1-(-1)^je^{-i\vt_p}\big)=0.
$$
This equation has $d-1$ linearly independent solutions
\[\lb{co2}
\begin{aligned}
\textstyle X_{j,1,\be_1}(\vt)=
-\, C_{j,1}\,X_{j,1,\be_d}(\vt), \qqq X_{j,1,\be_2}(\vt)=0,\\
 \ldots\,, \qqq X_{j,1,\be_{d-1}}(\vt)=0, \qqq X_{j,1,\be_d}(\vt);\\
\textstyle X_{j,2,\be_1}(\vt)=0, \qqq X_{j,2,\be_2}(\vt)
=-C_{j,2}\,X_{j,2,\be_d}(\vt),
\\ \ldots\,, \qqq X_{j,2,\be_{d-1}}(\vt)=0, \qqq X_{j,2,\be_d}(\vt);\\
\ldots\hspace{5cm}\\
\textstyle X_{j,d-1,\be_1}(\vt)=0, \qqq X_{j,d-1,\be_2}(\vt)=0,\qqq \ldots\,,\\
X_{j,d-1,\be_{d-1}}(\vt)=-C_{j,d-1}
X_{j,d-1,\be_d}(\vt), \qqq X_{j,d-1,\be_d}(\vt),
\end{aligned}
\]
where $C_{j,s}=C_{j,s}(\vt)$ are defined in \er{X1}.
Each constant $X_{j,s,\be_d}(\vt)$,  $s\in\N_{d-1}$, is
determined by the condition $\|\P_{j,s}^0(\vt,\cdot)\|_{L^2(\dL_*^d)}=1$.
Then \er{co1}, \er{co2} give
\begin{multline*}
\textstyle 1=\|\P_{j,s}^0(\vt,\cdot)\|_{L^2(\dL_*^d)}^2=
\sum\limits_{p=1}^d\int\limits_0^1 \big|\P_{j,s,\be_p}^0(\vt,t)\big|^2dt=
\int\limits_0^1 \big|\P_{j,s,\be_s}^0(\vt,t)\big|^2dt+\int\limits_0^1
\big|\P_{j,s,\be_d}^0(\vt,t)\big|^2dt
\\
=2\Big(\big|X_{j,s,\be_s}(\vt)\big|^2+ \big|X_{j,s,\be_d}(\vt)
\big|^2\Big)\int\limits_0^1\sin^2(\pi jt)\,dt
=\big(1+|C_{j,s}|^2\big)\big|X_{j,s,\be_d}(\vt) \big|^2,
\end{multline*}
which yields that the constant $X_{j,s,\be_d}(\vt)$ may be chosen in the
form $X_{j,s,\be_d}(\vt)={1\/\sqrt{1+|C_{j,s}|^2}}>0$\,. Substituting this and \er{co2} into \er{co1}, we obtain \er{X1}.
\BBox

\medskip

From  \er{sql2}  and  Theorem \ref{TT2} we deduce that
the operator $\D_{MV}(\vt)$ on the metric graph $\dL_*^d$, defined by \er{1S}, has the spectral representation
$$
\D_{MV}(\vt)=\sum_{j=0}^\iy z_j^2(\vt)\cP_j(\vt),\qq
\textrm{ for all } \; \vt\in\T^d\sm\{0\},
$$
where its eigenvalues $z_j^2(\vt)$ and the corresponding normalized
eigenfunctions \linebreak $\P_j(\vt)=\big(\P_{j,\be}(\vt,t)\big)_{\be\in\cE_*}$ are given by:
\[
\lb{egv1.s}
\begin{aligned}
& z_j(\vt)=\ca
  z(\vt)+\pi j, & \qq j \textrm{ is even} \\
  (\pi-z(\vt))+\pi j, & \qq j \textrm{ is odd}
\ac, \qqq z(\vt)=\arccos\big(-\l(\vt)\big)\in [0,\pi],\\[6pt]
&\P_{j,\be_s}(\vt,t)=\textstyle {1\/\sqrt{d}\;\sin z(\vt)}\,
\big(\sin(z_{j}(\vt)\,(1-t))+e^{i\vt_s}\,\sin(z_{j}(\vt)\,t)\big),\qqq
t\in [0,1].
\end{aligned}
\]

Since $\l(\vt)=-{1\/d}\,(\cos\vt_1+\ldots+\cos\vt_d)$, we have
$$
[\l^-,\l^+]=[-1,1], \qqq  z^-=\arccos(-\l^-)=0,\qqq
z^+=\arccos(-\l^+)=\pi.
$$
Recall that $\Z_0=\{0,1,2,\ldots\,\}$. Then, due to the formulas
\er{sDM1}, \er{sDM1a}, the Laplacian $\D_M$ on $L^2(\dL^d)$ has the
spectrum given by
\[\lb{sDM1.s}
\begin{array}{ll}
\s(\D_M)=\s_{ac}(\D_M)\cup \s_{fb}(\D_M),\qq &
\s_{fb}(\D_M)= \s_{fb}(\D_{MD})=\s_D,\\[6pt]
\s_{ac}(\D_M)=\bigcup_{j\in\Z_0}\s_j(\D_M)=[0,+\iy),\qq & \s_j(\D_M)=\big[(\pi j)^2,(\pi+\pi j)^2\big].
\end{array}
\]

\section {\lb{Sec6}Graphene lattice}
\setcounter{equation}{0}

We consider the hexagonal lattice $\bG$ shown in Fig.\ref{ff.0.3}\emph{a}. The periods $a_1,a_2$ of the lattice $\bG$ are also shown in the figure. The fundamental graph $\bG_*=(\cV_*,\cE_*)$, where $\cV_*=\{v_1,v_2\}$, consists of two vertices and three \emph{multiple} edges $\be_1,\be_2,\be_3$, having the form $(v_1,v_2)$ (Fig.\ref{ff.0.3}\emph{b}), with the indices $\t(\be_1)=(1,0)$, $\t(\be_2)=(0,1)$, $\t(\be_3)=(0,0)$.

The fiber Laplacian $\D(\vt)$, $\vt=(\vt_1,\vt_2)\in\T^2$, in the standard orthonormal basis of the space $\ell^2(\cV_*)$ has the form
\[
\lb{Dx}
\textstyle \D(\vt)=\ma  0 & \f(\vt)\\
                             \bar\f(\vt) & 0\am, \qqq
                              \f(\vt )=-{1\/3}\,(1+e^{-i\vt_1}+e^{-i\vt_2}).
\]
Then the eigenvalues of $\D(\vt)$ are given by
\[
\lb{zen}
\l_1(\vt)=-\l_2(\vt)=-|\f(\vt)|.
\]

\setlength{\unitlength}{1.0mm}
\begin{figure}[h]
\centering

\unitlength 1mm % = 2.845pt
\linethickness{0.4pt}
\ifx\plotpoint\undefined\newsavebox{\plotpoint}\fi % GNUPLOT compatibility
\begin{picture}(100,45)(0,0)

\put(5,10){(\emph{a})}

% √ексагональна€ решетка
\put(14,10){\circle{1}}
\put(28,10){\circle{1}}
\put(34,10){\circle{1}}
\put(48,10){\circle{1}}

\put(18,16){\circle{1}}
\put(24,16){\circle*{1}}
\put(38,16){\circle{1}}
\put(44,16){\circle{1}}

\put(14,22){\circle{1}}
\put(28,22){\circle*{1}}
\put(34,22){\circle{1}}
\put(48,22){\circle{1}}

\put(18,28){\circle{1}}
\put(24,28){\circle{1}}
\put(38,28){\circle{1}}
\put(44,28){\circle{1}}

\put(14,34){\circle{1}}
\put(28,34){\circle{1}}
\put(34,34){\circle{1}}
\put(48,34){\circle{1}}

\put(18,40){\circle{1}}
\put(24,40){\circle{1}}
\put(38,40){\circle{1}}
\put(44,40){\circle{1}}

\put(28,10){\line(1,0){6.00}}
\put(18,16){\line(1,0){6.00}}
\put(38,16){\line(1,0){6.00}}

\put(28,22){\line(1,0){6.00}}
\put(18,28){\line(1,0){6.00}}
\put(38,28){\line(1,0){6.00}}

\put(28,34){\line(1,0){6.00}}
\put(18,40){\line(1,0){6.00}}
\put(38,40){\line(1,0){6.00}}

\put(14,10){\line(2,3){4.00}}
\put(34,10){\line(2,3){4.00}}
\put(24,16){\line(2,3){4.00}}
\put(44,16){\line(2,3){4.00}}

\put(14,22){\line(2,3){4.00}}
\put(34,22){\line(2,3){4.00}}
\put(24,28){\line(2,3){4.00}}
\put(44,28){\line(2,3){4.00}}

\put(14,34){\line(2,3){4.00}}
\put(34,34){\line(2,3){4.00}}

\put(28,10){\line(-2,3){4.00}}
\put(48,10){\line(-2,3){4.00}}
\put(38,16){\line(-2,3){4.00}}
\put(18,16){\line(-2,3){4.00}}

\put(28,22){\line(-2,3){4.00}}
\put(48,22){\line(-2,3){4.00}}
\put(38,28){\line(-2,3){4.00}}
\put(18,28){\line(-2,3){4.00}}

\put(28,34){\line(-2,3){4.00}}
\put(48,34){\line(-2,3){4.00}}

\put(30,18){$\scriptstyle a_1$}
\put(20.5,22){$\scriptstyle a_2$}

\put(24,16){\vector(0,1){12.0}}
\put(33,21.3){\vector(3,2){0.5}}

%\drawline(24,16)(34,22)
\qbezier(24,16)(29,19)(34,22)

\put(24.8,21.5){$\scriptstyle v_1$}
\put(35,21.5){$\scriptstyle v_2+a_1$}
\put(25.0,28.0){$\scriptstyle v_2+a_2$}
\put(21.8,13.5){$\scriptstyle v_2$}

%***************************
\put(75,10){\circle*{1}}
\put(83,21){\circle*{1}}

\put(75,30){\circle{1}}
\put(95,40){\circle{1}}
\put(95,20){\circle{1}}

\put(75,10){\vector(0,1){20.0}}
\put(75,10){\vector(2,1){20.0}}

\multiput(95,20)(0,7){3}{\line(0,1){4}}
\put(75,30){\line(2,1){4.0}}
\put(82,33.5){\line(2,1){4.0}}
\put(89,37){\line(2,1){4.0}}

\qbezier(83,21)(89,20.5)(95,20)
\qbezier(83,21)(79,15.5)(75,10)
\qbezier(83,21)(79,25.5)(75,30)

\put(71,8.0){$v_2$}
\put(83,22){$v_1$}
\put(96,19.0){$v_2$}
\put(71.0,31.0){$v_2$}
\put(93.5,42.0){$v_2$}
\put(85,13){$a_1$}
\put(70,20.0){$a_2$}

\put(76,17.2){$\be_3$}
\put(89,21.2){$\be_1$}
\put(78,26.6){$\be_2$}

\put(64,10){(\emph{b})}
\end{picture}

\vspace{-0.5cm} \caption{\footnotesize  \emph{a}) Graphene $\bG$; \quad \emph{b}) the fundamental graph $\bG_*$ of graphene.} \label{ff.0.3}
\end{figure}

\begin{proposition}\lb{T.Gr'}
i) The fiber operator $\D_M(\vt)$ acts on
$y=(y_\be)_{\be\in\cE_*}\in L^2(\bG_*)$ by
$(\D_M(\vt)y)_\be=-y''_\be$, where $(y''_\be)_{\be\in\cE_*}\in L^2(\bG_*)$ and $y$ satisfies the quasi-periodic
conditions:
\[
\lb{FBC'}
y_{\be_1}(0)=y_{\be_2}(0)=y_{\be_3}(0),\qqq
e^{-i\vt_1}y_{\be_1}(1)=e^{-i\vt_2}y_{\be_2}(1)=y_{\be_3}(1),
\]
\[
\lb{dii2''}
\begin{array}{c}
y_{\be_1}'(0)+y_{\be_2}'(0)+y_{\be_3}'(0)=0,\\[6pt]
e^{-i\vt_1}y_{\be_1}'(1)+e^{-i\vt_2}y_{\be_2}'(1)+y_{\be_3}'(1)=0.
\end{array}
\]

ii) Let $j\in \N$. Then for each $\vt\in\T^2\sm\{0\}$ the operator $\D_M(\vt)$ has the simple eigenvalue $(\pi j)^2$ and the corresponding normalized eigenfunction $\P_j^0$ has the form
\[\lb{ve2}
\textstyle \P_j^0=\big(\P_{j,\be}^0\big)_{\be\in\cE_*},\qqq
\P_{j,\be}^0(\vt,t)=X_{\be}(\vt)\sqrt{2}\,\sin(\pi jt)\,, \qq t\in [0,1],
\]
where
\[\lb{ve3}
X_{\be_1}(\vt)=-X_{\be_3}(\vt)
{1-e^{-i\vt_2}\/e^{-i\vt_1}-e^{-i\vt_2}}\,,\qqq
X_{\be_2}(\vt)=X_{\be_3}(\vt)
{1-e^{-i\vt_1}\/e^{-i\vt_1}-e^{-i\vt_2}}\,,
\]
\[\lb{ve4}
X_{\be_3}(\vt)=
\bigg({1-\cos(\vt_1-\vt_2)\/3-\cos\vt_1-\cos\vt_2-\cos(\vt_1-\vt_2)}\bigg)^{1/2}\,.
\]
\end{proposition}

\no \textbf{Proof.} \emph{i}) Substituting the indices $\t(\be_1)=(1,0)$, $\t(\be_2)=(0,1)$, $\t(\be_3)=(0,0)$ of the fundamental graph edges  into the formulas \er{FBC}, \er{di2}, we obtain the conditions \er{FBC'} -- \er{dii2''}.

\emph{ii}) The Betti number $\b$ of the fundamental graph $\bG_*$ is equal to $2$. Then, due to Theorem \ref{TT1}, for all $\vt\in\T^2\sm\{0\}$ the eigenvalue $(\pi j)^2$ is simple, and
the unique normalized eigenfunction, corresponding to the
eigenvalue $(\pi j)^2$, has the form \er{ve2}, where
$\big(X_{\be}(\vt)\big)_{\be\in\cE_*}$ is the normalized solution
of the system of two equations
\[\lb{ve5}
\begin{array}{c}
x_{\be_1}+x_{\be_2}+x_{\be_3}=0, \\[6pt]
e^{-i\vt_1}\,x_{\be_1}+e^{-i\vt_2}\,x_{\be_2}+x_{\be_3}=0,
\end{array}
\]
which yields \er{ve3}. The constant $X_{\be_3}(\vt)$ is found using the condition $\big\|\P_j^0(\vt,\cdot)\big\|_{L^2(\bG_*)}=1$.
Using \er{ve2}, \er{ve3}, we have
\begin{multline*}
1=\big\|\P_{j}^0(\vt,\cdot)\big\|_{L^2(\bG_*)}^2=
\sum\limits_{s=1}^3\int\limits_0^1
\big|\P_{j,\be_s}(\vt,t)\big|^2dt=
2\Big(\sum\limits_{s=1}^3\big|X_{\be_s}(\vt)\big|^2\Big)\int\limits_0^1
\sin^2(\pi jt)\;dt\\
=\Big(1+\Big|{1-e^{-i\vt_1}\/e^{-i\vt_1}-e^{-i\vt_2}}\Big|^2+
\Big|{1-e^{-i\vt_2}\/e^{-i\vt_1}-e^{-i\vt_2}}\Big|^2\Big)\big|X_{\be_3}(\vt)
\big|^2
=\Big(1+{2-\cos\vt_1-\cos\vt_2\/1-\cos(\vt_1-\vt_2)}\Big)
\big|X_{\be_3}(\vt)\big|^2,
\end{multline*}
which yields that the constant $X_{\be_3}(\vt)$ may be chosen in the
form \er{ve4}. \qq \BBox

\begin{proposition}\label{TT2.h}
The operator $\D_{MV}(\vt)$ on
$\bG_*$, defined by  \er{1S}, has the form
\[
\lb{DVS.h} \D_{MV}(\vt)=\sum_{j=0}^\iy z_{j,1}^2(\vt)\cP_{j,1}(\vt)
+\sum_{j=0}^\iy z_{j,2}^2(\vt)\cP_{j,2}(\vt), \qq \textrm{for all} \qq \vt\in\T^2\sm\{0\},
\]
where its eigenvalues $z_{j,n}^2(\vt)$ and
the corresponding normalized eigenfunctions $\P_{j,n}(\vt)$ \linebreak $=\big(\P_{j,n,\be}(\vt,t)\big)_{\be\in\cE_*}$ have the form:
\[\lb{egv1.h}
\begin{aligned}
z_{j,n}(\vt)=\ca
  z_n(\vt)+\pi j, &  j \textrm{ is even} \\
  (\pi-z_n(\vt))+\pi j, &  j \textrm{ is odd}
\ac, \qq
z_n(\vt)=\textstyle\arccos\big({(-1)^{n+1}}\,|\f(\vt)|\big)\in[0,\pi],
\end{aligned}
\]
\[
\lb{ms6.h}
\begin{aligned}
\P_{j,n,\be_s}(\vt,t)=\textstyle {1\/\sqrt{3}\;\sin z_n(\vt)}\,
\Big({(-1)^n\f(\vt)\/|\f(\vt)|}\,\sin(z_{j,n}(\vt)\,(1-t))+ e^{i\lan\t
({\be_s}),\,\vt\ran}\,\sin(z_{j,n}(\vt)\,t)\Big),
\end{aligned}
\]
$\f(\vt)$ is defined in \er{Dx}.
\end{proposition}
\no {\bf Proof.}
Using \er{DVS}, \er{egv1} and \er{zen}, we obtain \er{DVS.h}, \er{egv1.h}. From \er{Dx} it follows that the normalized eigenfunction $\p_n(\vt,\cdot)$, $n=1,2$, corresponding
to the eigenvalue $\l_n(\vt)$ satisfies
\[\lb{hel2}
\begin{array}{c}
-\l_n(\vt)\p_n(\vt,v_1)+\f(\vt)\p_n(\vt,v_2)=0, \\[6pt]
3\big(|\p_n(\vt,v_1)|^2+|\p_n(\vt,v_2)|^2\big)=1,
\end{array}
\]
which yields
\[\lb{hel3}
\p_n(\vt,v_1)={\f(\vt)\/\l_n(\vt)}\,\p_n(\vt,v_2)=
{\f(\vt)\/(-1)^n|\f(\vt)|}\,\p_n(\vt,v_2)
\]
and we may choose $\p_n(\vt,v_2)=\frac1{\sqrt{6}}$.
Substituting this and \er{hel3} into \er{ms6} of Theorem
\ref{TT2}, we obtain \er{ms6.h}. \qq \BBox

\begin{corollary}\label{TsDM.h}
The Laplacian $\D_M$ on $L^2(\bG)$ has the
spectrum given by
\[\lb{sDM1.h}
\begin{aligned}
&\s(\D_M)=\s_{ac}(\D_M)\cup \s_{fb}(\D_M),\\
&\s_{ac}(\D_M)=[0,+\iy),\qqq \s_{fb}(\D_M)=\s_{fb}(\D_{MD})=\s_D.
\end{aligned}
\]
\end{corollary}
\no {\bf Proof.}
The eigenvalues of $\D(\vt)$ are given by
$\l_n(\vt)={(-1)^n\/3}\,|1+e^{-i\vt_1}+e^{-i\vt_2}|$,
$n=1,2$,
which yields
$$
[\l_1^-,\l_1^+] =\l_1(\T^2)=[-1,0],\qqq [\l_2^-,\l_2^+]
=\l_2(\T^2)=[0,1]
$$
and
$$
\textstyle z_1^-=\arccos(-\l_1^-)=0,\qqq z_1^+=\arccos(-\l_1^+)={\pi\/2}\,,
$$
$$
\textstyle z_2^-=\arccos(-\l_2^-)={\pi\/2}\,,\qqq z_2^+=\arccos(-\l_2^+)=\pi.
$$
Then, using the formulas \er{sDM1}, \er{sDM1a}, we obtain
$$
\begin{aligned}
&\s_{j,1}(\D_M)=[E_{j,1}^-,E_{j,1}^+]=\ca
  \big[(\pi j)^2,\big({\pi\/2}+\pi j\big)^2\big], & \qq j \textrm{ is even} \\
  \big[\big({\pi\/2}+\pi j\big)^2,(\pi+\pi j)^2\big], & \qq j \textrm{ is odd}
\ac,
\\
&\s_{j,2}(\D_M)=[E_{j,2}^-,E_{j,2}^+]=\ca
  \big[\big({\pi\/2}+\pi j\big)^2,\big(\pi+\pi j\big)^2\big], & \qq j \textrm{ is even} \\
  \big[(\pi j)^2,\big({\pi\/2}+\pi j\big)^2\big], & \qq j \textrm{ is odd}
\ac,\qqq j\in\Z_0=\{0,1,2,\ldots\,\},\\
&\s_{ac}(\D_M)=\bigcup_{j\in\Z_0}\bigcup_{n=1}^2\s_{j,n}(\D_M)=[0,+\iy),\qq
\s_{fb}(\D_{MV})=\varnothing,\qq \s_{fb}(\D_M)=
\s_{fb}(\D_{MD})=\s_D.
\end{aligned}
$$
Thus, \er{sDM1.h} has been proved.  \BBox

\begin{remark}
The spectrum of the Schr\"odinger operator $H=\D_M+Q$ with a periodic potential $Q$ on the hexagonal lattice was studied in \cite{KP07}.
\end{remark}

\section {\lb{Sec7}Stanene  lattice}
\setcounter{equation}{0}

Stanene is a topological insulator, theoretically predicted by Prof. Shoucheng Zhang's group at Stanford, which may display dissipationless current at its edges above room temperature \cite{Z13}. It is composed of tin atoms arranged in a single layer, in a manner similar to graphene.
The addition of fluorine atoms to the tin lattice could extend the
critical temperature up to $100^\circ$C. This would make it
practical for use in integrated circuits to make smaller, faster and
more energy efficient computers. Stanene has a band gap, it is a
semiconducting material. That makes it useful as material for use in
a transistor, which must have a component that turns on and off. For
more details see \cite{Z13} and references therein.

\setlength{\unitlength}{1.0mm}
\begin{figure}[h]
\centering

\unitlength 1mm % = 2.845pt
\linethickness{0.4pt}
\ifx\plotpoint\undefined\newsavebox{\plotpoint}\fi % GNUPLOT compatibility
\begin{picture}(100,45)(0,0)

\put(3,5){(\emph{a})}

% √ексагональна€ решетка
\put(14,10){\circle{1}}
\put(14,10){\line(0,-1){4.00}}
\put(14,6){\circle{1}}
\put(14,22){\line(0,-1){4.00}}
\put(14,18){\circle{1}}
\put(14,34){\line(0,-1){4.00}}
\put(14,30){\circle{1}}

\put(28,10){\circle{1}}
\put(28,10){\line(0,1){4.00}}
\put(28,14){\circle{1}}
\put(28,22){\line(0,1){4.00}}
\put(28,26){\circle*{1}}
\put(28,34){\line(0,1){4.00}}
\put(28,38){\circle{1}}

\put(34,10){\circle{1}}
\put(34,10){\line(0,-1){4.00}}
\put(34,6){\circle{1}}
\put(34,22){\line(0,-1){4.00}}
\put(34,18){\circle{1}}
\put(34,34){\line(0,-1){4.00}}
\put(34,30){\circle{1}}

\put(48,10){\circle{1}}
\put(48,10){\line(0,1){4.00}}
\put(48,14){\circle{1}}
\put(48,22){\line(0,1){4.00}}
\put(48,26){\circle{1}}
\put(48,34){\line(0,1){4.00}}
\put(48,38){\circle{1}}

\put(18,16){\circle{1}}
\put(18,16){\line(0,1){4.00}}
\put(18,20){\circle{1}}
\put(18,28){\line(0,1){4.00}}
\put(18,32){\circle{1}}
\put(18,40){\line(0,1){4.00}}
\put(18,44){\circle{1}}

\put(24,16){\circle*{1}}
\put(24,16){\line(0,-1){4.00}}
\put(24,12){\circle*{1}}
\put(24,28){\line(0,-1){4.00}}
\put(24,24){\circle{1}}
\put(24,40){\line(0,-1){4.00}}
\put(24,36){\circle{1}}

\put(38,16){\circle{1}}
\put(38,16){\line(0,1){4.00}}
\put(38,20){\circle{1}}
\put(38,28){\line(0,1){4.00}}
\put(38,32){\circle{1}}
\put(38,40){\line(0,1){4.00}}
\put(38,44){\circle{1}}

\put(44,16){\circle{1}}
\put(44,16){\line(0,-1){4.00}}
\put(44,12){\circle{1}}
\put(44,28){\line(0,-1){4.00}}
\put(44,24){\circle{1}}
\put(44,40){\line(0,-1){4.00}}
\put(44,36){\circle{1}}

\put(14,22){\circle{1}}
\put(28,22){\circle*{1}}

\put(34,22){\circle{1}}
\put(48,22){\circle{1}}

\put(18,28){\circle{1}}
\put(24,28){\circle{1}}
\put(38,28){\circle{1}}
\put(44,28){\circle{1}}

\put(14,34){\circle{1}}
\put(28,34){\circle{1}}
\put(34,34){\circle{1}}
\put(48,34){\circle{1}}

\put(18,40){\circle{1}}
\put(24,40){\circle{1}}
\put(38,40){\circle{1}}
\put(44,40){\circle{1}}

\put(28,10){\line(1,0){6.00}}
\put(18,16){\line(1,0){6.00}}
\put(38,16){\line(1,0){6.00}}

\put(28,22){\line(1,0){6.00}}
\put(18,28){\line(1,0){6.00}}
\put(38,28){\line(1,0){6.00}}

\put(28,34){\line(1,0){6.00}}
\put(18,40){\line(1,0){6.00}}
\put(38,40){\line(1,0){6.00}}

\put(14,10){\line(2,3){4.00}}
\put(34,10){\line(2,3){4.00}}
\put(24,16){\line(2,3){4.00}}
\put(44,16){\line(2,3){4.00}}

\put(14,22){\line(2,3){4.00}}
\put(34,22){\line(2,3){4.00}}
\put(24,28){\line(2,3){4.00}}
\put(44,28){\line(2,3){4.00}}

\put(14,34){\line(2,3){4.00}}
\put(34,34){\line(2,3){4.00}}

\put(28,10){\line(-2,3){4.00}}
\put(48,10){\line(-2,3){4.00}}
\put(38,16){\line(-2,3){4.00}}
\put(18,16){\line(-2,3){4.00}}

\put(28,22){\line(-2,3){4.00}}
\put(48,22){\line(-2,3){4.00}}
\put(38,28){\line(-2,3){4.00}}
\put(18,28){\line(-2,3){4.00}}

\put(28,34){\line(-2,3){4.00}}
\put(48,34){\line(-2,3){4.00}}

\put(27.8,19.5){$\scriptstyle v_1$}
\put(28.8,26.5){$\scriptstyle v_3$}
\put(22.0,17.5){$\scriptstyle v_2$}
\put(21.0,10.5){$\scriptstyle v_4$}

%***************************
\put(76.5,12){\vector(-2,-3){1.00}}
\put(76.5,28.2){\vector(-2,3){1.00}}
\put(93.0,20){\vector(4,-1){1.00}}

\put(75,10){\circle*{1}}
\put(75,4){\vector(0,1){6.00}}
\put(75,4){\circle*{1}}

\put(83,21){\circle*{1}}
\put(83,21){\vector(0,1){6.00}}
\put(83,27){\circle*{1}}

\put(75,30){\circle{1}}
\put(95,40){\circle{1}}
\put(95,20){\circle{1}}

\put(75,10){\vector(0,1){20.0}}
\put(75,10){\vector(2,1){20.0}}

\multiput(95,20)(0,7){3}{\line(0,1){4}}
\put(75,30){\line(2,1){4.0}}
\put(82,33.5){\line(2,1){4.0}}
\put(89,37){\line(2,1){4.0}}

\qbezier(83,21)(89,20.5)(95,20)
\qbezier(83,21)(79,15.5)(75,10)
\qbezier(83,21)(79,25.5)(75,30)

\put(70.5,9.0){$v_2$}
\put(70.5,4.0){$v_4$}
\put(83,18){$v_1$}
\put(83,28){$v_3$}
\put(96,19.0){$v_2$}
\put(71.0,31.0){$v_2$}
\put(93.5,42.0){$v_2$}
\put(85,12.5){$a_1$}
\put(70.5,20.0){$a_2$}

\put(77,17.2){$\scriptstyle\be_3$}
\put(89,21.2){$\scriptstyle\be_1$}
\put(78,27.0){$\scriptstyle\be_2$}
\put(83.5,23.5){$\scriptstyle\be_4$}
\put(75.5,6.2){$\scriptstyle\be_5$}

\put(62,5){(\emph{b})}
\end{picture}

\vspace{-0.5cm} \caption{\footnotesize  \emph{a}) Stanene $\bS$; \quad \emph{b}) the fundamental graph $\bS_*$ of stanene.} \label{ff.0.3'}
\end{figure}

A vertex of degree one is called an end vertex. An edge incident to an end vertex is called a pendant edge.
The stanene lattice $\bS$ is obtained from the hexagonal
lattice $\bG$ by adding a pendant edge at each vertex of $\bG$ (see
Fig.\ref{ff.0.3'}a). We choose the orientations of the edges as shown in Fig.\ref{ff.0.3'}b. We consider the discrete and metric Laplacians on
the stanene lattice $\bS$. The fundamental graph $\bS_*$ consists of 4 vertices
$v_1,v_2,v_3,v_4$ with degrees $\vk_{v_1}=\vk_{v_2}=4$, $\vk_{v_3}=\vk_{v_4}=1$ and 5 oriented edges
$$
\be_1=(v_1,v_2),\qq \be_2=(v_1,v_2),\qq \be_3=(v_1,v_2),\qq \be_4=(v_1,v_3),\qq \be_5=(v_4,v_2).
$$
The indices of the fundamental graph edges are given by
\[\lb{sta1}
\t(\be_1)=(1,0),\qqq \t(\be_2)=(0,1),\qqq
\t(\be_3)=\t(\be_4)=\t(\be_5)=(0,0).
\]
For each $\vt=(\vt_1,\vt_2)\in\T^2$ the fiber Laplacian $\D(\vt)$
in the standard orthonormal basis of the space $\ell^2(\cV_*)$ has the form
\[
\lb{el}
\D(\vt)=-
{1\/2}\begin{pmatrix}
0&b(\vt)&1&0\\
\overline b(\vt)&0&0&1 \\
1&0&0&0 \\
0&1&0&0 \\
\end{pmatrix},\qqq b(\vt)={1\/2}\big(1+e^{-i\vt_1}+e^{-i\vt_2}\big).
\]
A direct calculation gives
\[\lb{egv1'}
\textstyle\det\big(\D(\vt)-\l
\1_4\big)=\l^4-\big(\frac{1}{2}+\frac{|b(\vt)|^2}{4}\big)
\l^2+{1\/16}\,,
\]
where $\1_4$ is the identity $4\ts4$ matrix.
The eigenvalues of the matrix $\D(\vt)$ are given by
\[
\lb{ev}
\l_{1,2,3,4}(\vt)=\pm\frac{|b(\vt)|}4\pm\frac{\sqrt{|b(\vt)|^2+4}}{4}\,.
\]

\begin{proposition}\lb{T.StD}
The spectrum of the discrete Laplacian $\D$ on the stanene lattice
$\bS$ has the form
\[
\lb{sSL}
\begin{aligned}
&\textstyle\s(\D)=\s_{ac}(\D)=\bigcup\limits_{n=1}^4\s_n(\D)=
\big[-1;-{1\/4}\,\big]\cup\big[{1\/4}\,;1\big],\\
&\textstyle \s_1(\D)=[-1,-{1\/2}],\qq
\s_2(\D)=[-{1\/2},-{1\/4}],\qq \s_3(\D)=[{1\/4},{1\/2}],\qq
\s_4(\D)=[{1\/2},1].
\end{aligned}
\]
\end{proposition}

\no {\bf Proof.} We need the following statement (see Theorem 4.3.8
in \cite{HJ85}).

\emph{Let $B=\ma   A & y \\
y^* & a \am$ be a self-adjoint $(\nu+1)\ts(\nu+1)$ matrix
for some self-adjoint $\nu\ts \nu$ matrix $A$, some  real number $a$
and some vector $y\in\C^{\nu}$. Denote by $\l_1(A)\leq\ldots\leq\l_\n(A)$ the eigenvalues of $A$, arranged in
increasing order, counting multiplicities.
Then
$$
\l_1(B)\leq\l_1(A)\leq\l_2(B)\leq\ldots\leq\l_\n(B)\leq\l_\n(A)\leq\l_{\nu+1}(B).
$$}

Applying this statement to the matrix $\D(\vt)$, we obtain
\[
\lb{razz}
\textstyle \l_1(\vt)\le-{1\/2}\le\l_2(\vt)\le0\le\l_3(\vt)\le{1\/2}\le\l_4(\vt),\qqq \forall\,\vt\in\T^2.
\]
Since
$$
\textstyle \l_1(0)=-1,\qqq
\l_1\big({2\pi\/3}\,,-{2\pi\/3}\big)=-{1\/2}\,,
$$
using \er{razz} and the fact that the spectrum of a bipartite graph is symmetric with respect to the point 0, we have
\[\lb{egv2}
\textstyle \s_1=\l_1(\T^2)=[-1,-{1\/2}],\qqq \s_4=\l_4(\T^2)=-\l_1(\T^2)=[{1\/2},1].
\]
From \er{egv1'} and \er{razz} it follows that
\[
\textstyle\l_1(\vt)\l_2(\vt)={1\/4}\,.
\]
This and \er{egv2} yield
\[\lb{sb23}
\textstyle
\s_2=\l_2(\T^2)=\big[-{1\/2}\,,-{1\/4}\big],\qqq
\s_3=\l_3(\T^2)=-\l_2(\T^2)=\big[{1\/4}\,,{1\/2}\big].
\]
Thus, \er{sSL} has been proved.
\qq \BBox

\begin{proposition}\lb{T.StM}
i) The fiber operator $\D_M(\vt)$ acts on
$y=(y_\be)_{\be\in\cE_*}\in L^2(\bS_*)$ by
$(\D_M(\vt)y)_\be=-y''_\be$, where $(y''_\be)_{\be\in\cE_*}\in L^2(\bS_*)$ and
$y$ satisfies the quasi-periodic conditions:
\[
\lb{FBC''}
\begin{array}{c}
y_{\be_1}(0)=y_{\be_2}(0)=y_{\be_3}(0)=y_{\be_4}(0),\\[6pt]
e^{-i\vt_1}y_{\be_1}(1)=e^{-i\vt_2}y_{\be_2}(1)=y_{\be_3}(1)=y_{\be_5}(1),
\end{array}
\]
\[
\lb{di2''}
\begin{array}{c}
y_{\be_1}'(0)+y_{\be_2}'(0)+y_{\be_3}'(0)+y_{\be_4}'(0)=0,\\[6pt]
e^{-i\vt_1}y_{\be_1}'(1)+e^{-i\vt_2}y_{\be_2}'(1)+y_{\be_3}'(1)+
y_{\be_5}'(1)=0,\\[6pt]
y_{\be_4}'(1)=0,\qqq y_{\be_5}'(0)=0.
\end{array}
\]

ii) Let $j\in\N$. Then for each $\vt\in\T^2\sm\{0\}$ the operator $\D_M(\vt)$ has the simple eigenvalue $(\pi j)^2$ and the
corresponding normalized eigenfunction
$\P_j^0=\big(\P_{j,\be}^0\big)_{\be\in\cE_*}$ has the form
\[
\lb{sta3} \textstyle \P_{j,\be_s}^0(\vt,t)=X_{\be_s}(\vt)\sqrt{2}
\,{\sin(\pi jt)}\,, \qq s=1,2,3,\qqq \P_{j,\be_4}^0(\vt,t)
=\P_{j,\be_5}^0(\vt,t)\equiv0,
\]
where the constants $X_{\be_s}(\vt)$, $s=1,2,3$, are defined by
\er{ve3}, \er{ve4}.
\end{proposition}
\no \textbf{Proof.} \emph{i}) Substituting the indices \er{sta1} of the
fundamental graph edges into  the formulas \er{FBC}, \er{di2}, we
obtain the conditions \er{FBC''} -- \er{di2''}.

\emph{ii}) The Betti number $\b$ of the fundamental graph $\bS_*$ is equal to $2$. Then, due to
Theorem \ref{TT1}, for all $\vt\in\T^2\sm\{0\}$ the eigenvalue
$(\pi j)^2$ is simple and
the unique normalized eigenfunction, corresponding to the
eigenvalues $(\pi j)^2$, has the form
\[
\textstyle \P_{j,\be}^0(\vt,t)=X_{\be}(\vt)\sqrt{2} \,{\sin(\pi jt)}\,,
\]
where $\big(X_{\be}(\vt)\big)_{\be\in\cE_*}$ is the normalized
solution of the system of 4 equations
\[
\begin{array}{c}
x_{\be_1}+x_{\be_2}+x_{\be_3}+x_{\be_4}=0, \\[6pt]
e^{-i\vt_1}x_{\be_1}+ e^{-i\vt_2}\,x_{\be_2}+
x_{\be_3}+x_{\be_5}=0, \qqq x_{\be_4}=0, \qqq
x_{\be_5}=0,
\end{array}
\]
which yields the system \er{ve5} for  $X_{\be_s}$, $s=1,2,3$. The unique normalized solution of \er{ve5} is given by \er{ve3}, \er{ve4}. Thus, \er{sta3} has been proved. \qq
\BBox

\begin{proposition}\label{TT2.st}
The operator $\D_{MV}(\vt)$ on  $\bS_*$, defined by \er{1S}, has
the form
\[
\lb{DVS.st} \D_{MV}(\vt)=\sum_{j=0}^\iy\sum_{n=1}^4
z_{j,n}^2(\vt)\cP_{j,n}(\vt), \qq \textrm{for all} \qq \vt\in\T^2\sm\{0\},
\]
where its eigenvalues $z_{j,n}^2(\vt)$
and the corresponding  normalized eigenfunctions $\P_{j,n}(\vt)$ \linebreak
$=\big(\P_{j,n,\be}(\vt,t)\big)_{\be\in\cE_*}$ have
the form:
\[\lb{egv1.h'}
\begin{aligned}
 z_{j,n}(\vt)=\ca
  z_n(\vt)+\pi j, & \; j \textrm{ is even} \\
  (\pi-z_n(\vt))+\pi j, & \; j \textrm{ is odd}
\ac, \qq
z_n(\vt)=\textstyle\arccos\big(-\l_n(\vt)\big)\in[0,\pi],\\
\end{aligned}
\]
\[\lb{ms6.h'}
\begin{aligned}
&\P_{j,n,\be_s}(\vt,t)=\textstyle{\sqrt{2}\/\sin z_n(\vt)}
\big(\p_n(\vt,v_1)\sin(z_{j,n}(\vt)\,(1-t))+\p_n(\vt,v_2)e^{i\lan\t
(\be_s),\,\vt\ran}\,\sin(z_{j,n}(\vt)\,t)\big),\\
& \hspace{130mm} s=1,2,3,\\
&\P_{j,n,\be_4}(\vt,t)=\textstyle {\sqrt{2}\/\sin z_{n}(\vt)}\,
\big(\p_n(\vt,v_1)\sin(z_{j,n}(\vt)\,(1-t))+
\p_n(\vt,v_3)\sin(z_{j,n}(\vt)\,t)\big),\\
&\P_{j,n,\be_5}(\vt,t)=\textstyle {\sqrt{2}\/\sin z_{n}(\vt)}\,
\big(\p_n(\vt,v_4)\sin(z_{j,n}(\vt)\,(1-t))+
\p_n(\vt,v_2)\sin(z_{j,n}(\vt)\,t)\big),
\end{aligned}
\]
\[\lb{hel3'}
\begin{aligned}
&\textstyle\p_n(\vt,v_1)=\Big({2\l^2_n(\vt)\/16\l^2_n(\vt)+1}\Big)^{1/2},\hspace{11mm}
\p_n(\vt,v_2)=-{1\/b(\vt)}\,\big(2\l_n(\vt)-{1\/2\l_n(\vt)}\big)
\p_n(\vt,v_1),\\
&\textstyle\p_n(\vt,v_3)=-{1\/2\l_n(\vt)}\,\p_n(\vt,v_1), \qq
\p_n(\vt,v_4)={1\/2\l_n(\vt) b(\vt)}\,\big(2\l_n(\vt)-{1\/2\l_n(\vt)}\big)\p_n(\vt,v_1),
\end{aligned}
\]
where $b(\vt)$ is defined in \er{el}.
\end{proposition}
\no {\bf Proof.}
Using \er{DVS}, \er{egv1}, \er{ms6}, we obtain \er{DVS.st}, \er{egv1.h'}, \er{ms6.h'}.

For each $\vt\in\T^2$ the fiber Laplacian $\D(\vt)$ has the form \er{el}
and the eigenvalues of $\Delta(\vt)$ are given by \er{ev}.
Then the eigenfunction $\p_n(\vt,\cdot)$, corresponding to the eigenvalue $\l_n(\vt)$ satisfies
\[\lb{hel2'}
\begin{array}{l}
\l_n(\vt)\p_n(\vt,v_1)+{1\/2}\,b(\vt)\,\p_n(\vt,v_2)+{1\/2}\,\p_n(\vt,v_3)=0,
\qqq
\l_n(\vt)\p_n(\vt,v_3)+{1\/2}\,\p_n(\vt,v_1)=0, \\[6pt]
\l_n(\vt)\p_n(\vt,v_2)+{1\/2}\,\bar b(\vt)\,\p_n(\vt,v_1)+{1\/2}\,\p_n(\vt,v_4)=0, \qqq
\l_n(\vt)\p_n(\vt,v_4)+{1\/2}\,\p_n(\vt,v_2)=0,
\end{array}
\]
which are normalized by $\sum\limits_{s=1}^4\vk_{v_s}\big|\p_n(\vt,v_s)\big|^2=1$.
Using \er{egv1'}, we obtain that
\[\lb{hel4'}
\textstyle|b(\vt)|=\Big|2\l_n(\vt)-{1\/2\l_n(\vt)}\Big|\,.
\]
Solving the system of the equations \er{hel2'} and applying \er{hel4'}, we have
\er{hel3'}. \qq \BBox

\begin{corollary}\label{TsDM.h'}
The Laplacian $\D_M$ on $L^2(\bS)$ has the
spectrum given by
\[\lb{sDM1'.h}
\begin{aligned}
&\s(\D_M)=\s_{ac}(\D_M)\cup \s_{fb}(\D_M),\qqq
\s_{fb}(\D_M)=\s_{fb}(\D_{MD})=\s_D,\\
&\textstyle\s_{ac}(\D_M)=\big[0,r_+^2\big]\textstyle\cup_{j\in\Z_0}
\big[\big(r_-+\pi j\big)^2,\big(r_++\pi(j+1)\big)^2\big],
\end{aligned}
\]
where  $r_\pm =\arccos\big(\pm {1\/4}\big)$.
\end{corollary}
\no {\bf Proof.} Due to Proposition \ref{T.StD}, the bands $[\l_n^-,\l_n^+]=\l_n(\T^2)$ satisfy
$$
\textstyle [\l_1^-,\l_1^+]=[-1,-{1\/2}],\qqq
[\l_2^-,\l_2^+]=[-{1\/2},-{1\/4}],\qqq [\l_3^-,\l_3^+]=[{1\/4},{1\/2}],\qqq
[\l_4^-,\l_4^+]=[{1\/2},1],
$$
and
$$
\begin{array}{ll}
\textstyle z_1^-=\arccos(-\l_1^-)=0,\qqq & z_1^+=\arccos(-\l_1^+)={\pi\/3}\,,\\
 [6pt]
\textstyle z_2^-=\arccos(-\l_2^-)={\pi\/3}\,,\qqq & z_2^+=\arccos(-\l_2^+)=
\arccos\big({1\/4}\big), \\ [6pt]
\textstyle z_3^-=\arccos(-\l_3^-)=\arccos\big(-{1\/4}\big)\,,\qqq & z_3^+=
\arccos(-\l_3^+)={2\pi\/3}\,, \\ [6pt]
\textstyle z_4^-=\arccos(-\l_4^-)={2\pi\/3}\,,\qqq & z_4^+=
\arccos(-\l_4^+)=\pi.
\end{array}
$$
Then, using the formulas \er{sDM1}, we obtain
$$
\begin{aligned}
&\s_{j,1}(\D_M)=[E_{j,1}^-,E_{j,1}^+]=\ca
  \big[(\pi j)^2,\big({\pi\/3}+\pi j\big)^2\big], & \qq j \textrm{ is even} \\
  \big[\big({2\pi\/3}+\pi j\big)^2,(\pi+\pi j)^2\big], & \qq j \textrm{ is odd}
\ac,
\\
&\s_{j,2}(\D_M)=[E_{j,2}^-,E_{j,2}^+]=\ca
  \big[\big({\pi\/3}+\pi j\big)^2,\big(\arccos\big({1\/4}\big)+\pi j\big)^2\big],
   & \qq j \textrm{ is even} \\
  \big[(\arccos\big(-{1\/4}\big)+\pi j)^2,\big({2\pi\/3}+\pi j\big)^2\big],
   & \qq j \textrm{ is odd}
\ac,\\
&\s_{j,3}(\D_M)=[E_{j,3}^-,E_{j,3}^+]=\ca
  \big[\big(\arccos\big(-{1\/4}\big)+\pi j\big)^2,\big({2\pi\/3}+\pi j\big)^2\big], & \qq j \textrm{ is even} \\
  \big[\big({\pi\/3}+\pi j\big)^2,(\arccos\big({1\/4}\big)+\pi j)^2\big], &
  \qq j \textrm{ is odd}
\ac,\\
&\s_{j,4}(\D_M)=[E_{j,4}^-,E_{j,4}^+]=\ca
  \big[({2\pi\/3}+\pi j)^2,\big(\pi+\pi j\big)^2\big], & \qq j \textrm{ is even} \\
  \big[\big(\pi j\big)^2,({\pi\/3}+\pi j)^2\big], & \qq j \textrm{ is odd}
\ac,
\end{aligned}
$$
for all $j\in\Z_0=\{0,1,2,\ldots\,\}$. Substituting these identities
into the formula
$$
\textstyle\s_{ac}(\D_M)=\bigcup\limits_{j\in\Z_0}\bigcup\limits_{n=1}^4
\s_{j,n}(\D_M),
$$
we obtain the third identity in \er{sDM1'.h}. Since $\s_{fb}(\D_{MV})=\varnothing$, we have the second formula in \er{sDM1'.h}.
\qq \BBox

\

\lb{Ack}
\no\textbf{Acknowledgments.} \footnotesize Evgeny Korotyaev was
supported by the RSF grant  No. 18-11-00032. Natalia Saburova was
supported by the RFBR grant No. 19-01-00094.

\lb{Ref}

\end{document}